\documentclass[11pt,a4paper]{article}
\usepackage{amsfonts}
\usepackage{amsmath,amssymb}
\usepackage{mathrsfs}
\usepackage{hyperref}
\usepackage{graphicx}
\usepackage{float}
\usepackage{indentfirst}
\newtheorem{thm}{Theorem}[section]

\newtheorem{lem}[thm]{Lemma}

\newtheorem{rem}[thm]{Remark}

\numberwithin{equation}{section}\allowdisplaybreaks
\def\i{\mathbf{i}}
\def\endproof{\hfill $\Box$}
\def\s{\,\,\,\,\,}
\begin{document}

\title{\large\bf Global Existence for The Massive Dirac Equations with small initial datum on Tori}

\author{\normalsize \bf $Zonglin\quad Jia^a$\footnote{Corresponding author.} \\
\footnotesize
\it a. Department of Mathematics and Physics, North China Electric Power University, Beijing 102206,  PR China.\\
\footnotesize
\it Emails: 50902525@ncepu.edu.cn\\
}
\date{} \maketitle
\begin{minipage}{13.5cm}
\footnotesize \bf Abstract. \rm In the article we obtain almost global existence for Dirac Equations with high regularity and small initial datum on Tori. Besides, the global existence with low regularity and small initial datum is gotten. The approaches are mainly Gagliardo-Nirenberg-Moser estimates and Bernstein-Type Lemma.
\\

\\

\end{minipage}

\section{Introduction}\label{section1}
The Dirac equation is a basic equation in quantum mechanics. Mathematicians' researches on this equation often focus on the case that the domain is a Euclidean space. This paper considers the case where the domain is a flat torus.

At first, let us introduce the $d$-dimensional flat torus $\mathbb{T}^d$. Its precise definition is $(\mathbb{R}/2\pi\mathbb{Z})^d$. Because $\mathbb{R}/2\pi\mathbb{Z}$ is isometric to the circle $\mathbb{S}^1$, $\mathbb{T}^d$ can be regarded as
$$\underbrace{\mathbb{S}^1\times\mathbb{S}^1\times\cdots\times\mathbb{S}^1}\limits_{\mbox{$d$-times}}.$$
Namely,
$$
\mathbb{T}^d:=\{(e^{\i x^1},e^{\i x^2},\cdots,e^{\i x^d})|x^j\in\mathbb{R},\quad j=1,2,\cdots,d\}
$$
with the imaginary unit $\i$.

In order to accurately describe the Dirac equation, we must first make some symbolic explanations. Sometimes we write $x^0$ as $t$ meaning the time variable belonging to $\mathbb{R}$. $\mu\in\{0,1,2,\cdots,d\}$ and $j\in\{1,2,\cdots,d\}$ are indices. $I_n$ means the unit matrix of order $n$. We study the Dirac equations of the following form
\begin{eqnarray}\label{1}
-\i\sum\limits_{\mu=0}^d\gamma^{\mu}\frac{\partial\psi}{\partial x^{\mu}}+m\psi=F(\psi),\s\s\s\s\psi(0,\cdot)=\psi_0
\end{eqnarray}
where $\gamma^{\mu}$ are Dirac matrices. The mass $m$ is a positive constant and the unknown spinor $\psi$ is a mapping from $\mathbb{R}\times\mathbb{T}^d$ to $\mathbb{C}^{d_0}$, while the known function $F$ maps $\mathbb{C}^{d_0}$ to $\mathbb{C}^{d_0}$ with $d_0:=2^{\lfloor (d+1)/2\rfloor}$ and $\lfloor x\rfloor$ being the integer part of a real number $x$.\\
\textbf{Previous Results}

In terms of the number of equations, the Dirac equation can be divided into a single equation and coupled equations of other physical systems, for example The Dirac-Klein-Gordon equations, Maxwell-Dirac equations, etc. Its structure can be divided into linear equation, quasilinear equation and nonlinear equation. In view of the above situations, mathematicians have made in-depth research and exploration. The ground state solution of the Dirac equation(i.e. the time-independent solution) has been discussed in many literatures. For a single nonlinear equation, Bartsch and Ding proved in \cite{2} that, for any real number $p$ not less than 2, the Dirac equation has a nontrivial minimum energy solution in the $L^p$-integrable Sobolev space of order 1 when the nonlinear part satisfies some smooth condition, growth condition and periodic one. In \cite{10} Ding considers semi-classical ground state solutions of quasilinear cases, where the quasilinear part is a power function of the norm of the solution. The author proves that the equation has a minimum energy solution for all small positive parameters. And when the parameter goes to 0, the minimum energy solution of this family tends to the minimum energy solution of the related limit problem. In addition, in the process of convergence, this family of solutions will focus on the global maximum point of the potential function in a certain sense. In \cite{12}, Ding and Liu consider the case where the quasilinear part is a general function and get a conclusion similar to \cite{10}. Compared with \cite{12}, the equation considered in \cite{13} is slightly different: the quasilinear part of the equation is added with a linear function without constant terms. \cite{14} improved the conclusion of \cite{13} and its condition was weakened as: the quasilinear term is either superlinear or asymptotically linear at positive infinity; The potential function of a linear term has a local maximum but not a global maximum. In \cite{3}, Bartsch and Xu improved the result of \cite{14}, that is, the potential function of the linear term is not required to have a maximum value point, but only a critical point. For equations coupled with other physical systems, Ding and Xu studied the semi-classical ground state solutions of quasilinear Maxwell-Dirac equations and nonlinear Dirac-Klein-Gordon equations in \cite{15} and \cite{16} respectively. The quasilinear part in \cite{15} is the same as that in \cite{12}, while \cite{16} contains the nonlinear term. They extend the results of \cite{12} to equations with electromagnetic fields and equations coupled with state functions of elementary particles respectively.

Mathematicians have also made a lot of explorations on the equation of variable including time. The scattering and non-scattering problems of hartree-type nonlinear Dirac systems with critical regularity were studied by Cho, Hong and Lee in \cite{5}. The nonlinear term of this equation is not a function but a functional, and is the convolution of a function called the Yukawa potential (under some conditions it is called the Coulomb potential) with the quadratic form of the solution. In this paper, they first prove that under the Yukawa potential, if the system admits small initial data and angular regularity, the solution scatters and has global well-posedness. Through the process, the authors observe that only a small amount of angular regularity is required to obtain global well-posedness. In addition, \cite{5} studies a class of solutions with given Coulomb potential, which are not scattering. In addition to the Dirac type equations with 3-spatial variables, the case of 2 elements has also been studied. In \cite{19}, Lee discussed two kinds of systems: one is massless honeycomb potential Dirac equation; The other is the Hartree Dirac equation. The potential of the former is a diagonal matrix of order 2, whose elements are complex quadratic forms of solutions. The nonlinear term of the latter is the convolution of the Coulomb potential with the square of the norm of the solution. In \cite{19} The results are as follows: first, these two kinds of systems have local well-posedness in square integrable fractional Sobolev space, where the order of Sobolev space is greater than a constant; second, for any square integrable Sobolev space of order less than the constant, the flow maps determined by these two kinds of systems, if them exist, are not differentiable at the origin of order 3. In the same case of 2 elements, \cite{6} considers the same Yukawa potential as \cite{5}. It is concluded that the Dirac equation is unique and exists globally with small initial data in a square integrable function space, and the solution scatterers. In order to make up for the deficiency of \cite{6}, \cite{7} studied the system of 2 elements Coulomb potential. Compared with the classical Coulomb potential, the potential function here is more general, but the conclusion is stronger: this equation has the global existence and uniqueness with small initial data in square integrable fractional Sobolev space, and the solutions scatter.

Some scholars also consider Dirac equations on manifolds. Ding and Li discussed in \cite{11} the boundary value problem of nonlinear Dirac equation on a compact Spin manifold of general dimension, which is a space-like hypersurface of a Lorentz manifold and a Riemannian manifold. The Dirac operator here is a self-adjoint elliptic operator of first order. The authors prove that the boundary value problem has a nontrivial weak solution. Further, if the nonlinear term is odd, then the problem has infinitely many weak solutions. The above result can be regarded as a discussion of the ground state solutions of Dirac equation on a general manifold. For special Lorentz manifolds, such as $(1+ N)$ dimensional Minkowski Spaces, some scholars also discuss the Dirac system above. In \cite{22}, Tzvetkov proved the uniqueness and global existence with small initial datum of a class of semilinear massless Dirac equations by contraction mapping principle, and obtained the decay properties. In \cite{20}, Li and Zang studied the asymptotic properties of linear and nonlinear Dirac equations without mass using vector field method: the solution of linear equation has negative power decay behavior with respect to spatial-time variables; Nonlinear equations (nonlinear part admits null structure) have global existence of small initial data and decay behavior of negative power, and the power is the same as linear equations.\\
\textbf{Our Strategy}

Noting that the Dirac equation is closely related to the wave equation or Klein-Gordon equation, we take time derivative to both ends of the Dirac equation to obtain a wave equation or Klein-Gordon equation, if we assume that the solution is $C^2$-continuously differentiable(From Sobolev Embedding Theorem it follows that the initial data $\psi_0$ must be in $H^s(\mathbb{T}^d,\mathbb{C}^{d_0})$ for $s\geq 2+d/2$). As is well known, people have done very in-depth research on these two types of equations. Related literatures can refer to \cite{DS,FHZ}. Using their methods, we briefly proved that when the regularity of the initial data is greater than $(9/2)d+9$, the Dirac equation admit a unique almost global solution. Therefore, the focus of this article is on the Dirac equation with low regularity.

Because Fourier analysis cannot be established on general compact manifolds, we have to use the periodicity of the flat torus to deal with the problem through the Fourier series. We define an inverse transformation and then give some frequency cut-off operators following the Littlewood-Paley theory. The idea here comes from \cite{DS}. By rewriting the Dirac operator, we turn the equation to be considered into its equivalent integral equation, and let the right hand side of integral equation to be an operator. This paper will employ the contraction mapping principle to get the global existence. The proof is divided into two steps: the first step is to ensure the operator is from a space into itself; the second one is to guarantee the operator is a contraction.

As for the first step, we decompose the estimation into a linear part and a non-linear part. In order to treat the linear part, we have to use a fact that the Besov Space is equivalent to the Sobolev Space in some special cases. For nonlinear term, Lemma \ref{lem1} tells us that we just need to prove an inequality (\ref{6}). Under the premise that $F$ is analytic near the origin, exploiting the Taylor Series(its coefficients are denoted by $\mathbf{c_p}$) leads to decomposing the nonlinear term into many single terms. Next, use frequency cut-off operators to decompose the unity for single terms and there will be a lot of summations for multiple indicators which will be ``canceled" by power series. From Plancherel formula it follows that we are able to rewrite these single terms as forms of convolutions. The frequency cut-off operators ensure that each item in the convolutions cannot be too far apart(This will play a key role in the subsequent process). In the next, we use Young's Inequality to make these single terms into many products of $L_t^pL_x^q$-norms, where the indices are to be determined. We shall employ Lemma \ref{lem2} to transform $L_t^pL_x^q$-norms into some new norms denoted by $||\cdot||_{S_j^{\pm}}$. However, that lemma requires $1/p+1/q=1/2$. In order to meet the assumption, we utilize Gagliardo-Nirenberg-Moser estimates to change $q$ into $\tilde{q}$, but the price to pay is to improve the regularity of the unknown functions. Thanks to Bernstein-Type Lemma, i.e. Lemma 2.1 in chapter 2 of \cite{BCD} we can let spatial derivatives disappear and some powers of a certain constant $C$ will arise. Choose indices of $L^p_tL^q_x$-norms appropriately to ``cancel" summations for multiple indicators said before by negative powers of $C$. So we get a power series with respect to $||\cdot||_{S^{\pm}_j}$-norms of the unknown functions. At last we condition $\mathbf{c_p}$ in order to guarantee that the power series is convergent.

The approach for the second step is similar, but the computation is more complicated.

Our main results are the following

\begin{thm}\label{thm1}
Given $s>(9/2)d+9$, there are $\epsilon_0>0$ and $c>0$ such that for any $\epsilon\in(0,\epsilon_0]$ and $\psi_0\in H^s(\mathbb{T}^d,\mathbb{C}^{d_0})$ with $||\psi_0||_{H^s}\leq\epsilon$, (\ref{1}) admits a unique solution
$$\psi\in C^0((-T_{\epsilon},T_{\epsilon}),H^s(\mathbb{T}^d,\mathbb{C}^{d_0}))\bigcap C^1((-T_{\epsilon},T_{\epsilon}),H^{s-1}(\mathbb{T}^d,\mathbb{C}^{d_0}))$$
with $T_{\epsilon}\geq c\epsilon^{-2}$, provided $F$ is smooth and vanishes at least at order 2 at 0.
\end{thm}

\begin{thm}\label{thm2}
Assume $F$ meets the conditions of Theorem \ref{thm1}. Moreover, $F$ is odd with respect to $\psi$, i.e. $F(-\psi)=-F(\psi)$. Then there is a zero measure subset $\mathcal{N}$ of $(0,\infty)$ and for any $m\in(0,\infty)\setminus\mathcal{N}$, any $N\in\mathbb{N}$, there exists an $s_0\geq(9/2)d+9$ such that for any $s>s_0$, there are $\epsilon_0>0$, $c>0$, $C>0$ and for any $\epsilon\in(0,\epsilon_0]$, any $\mathbb{C}^{d_0}$-valued function $\psi_0$ satisfying $||\psi_0||_{H^s}\leq\epsilon$, the following holds: the solution $\psi$ to (\ref{1}) exists at least over a time interval $(-c\epsilon^{-N},c\epsilon^{-N})$ with the uniform bound
$$
||\psi(t,\cdot)||_{H^s}\leq C\epsilon
$$
holding on that interval.
\end{thm}

Before illustrating the most important theorem of this article, we give a definition called Condition (A):
\subsection{A Hypothesis}
We say an function, which is analytic in a small neighborhood of $0$,
\begin{eqnarray*}
F(\psi)=\sum\limits_{|\mathbf{p}|=0}^{\infty}\mathbf{c}_{\mathbf{p}}\psi^{\mathbf{p}}
\end{eqnarray*}
(where $\mathbf{p}:=(p_1,p_2,\cdots,p_{d_0})$, $\mathbf{c}_{\mathbf{p}}\in\mathbb{C}^{d_0}$ and $\psi^{\mathbf{p}}:=\psi_1^{p_1}\cdots \psi_{d_0}^{p_{d_0}}$) meets Condition (A) if and only if
$$\limsup\limits_{p\rightarrow\infty}\max\{B_p,A_p\}^{1/p}<C^{-d/2-1/4}2^{-d/2}/3,$$
where the constant $C$ is from Bernstein-Type Lemma, i.e. Lemma 2.1 in chapter 2 of \cite{BCD}.

Let us now explain the above notation in detail.\\
\textbf{Case for $B_p$}

Set
\begin{eqnarray*}
C_{\vec{i}\vec{j}\vec{l}\vec{v}}:=c_{\mathbf{m}\mathbf{n}}:=\mathbf{c_p}E_{\mathbf{mn}}
\end{eqnarray*}
with
$$
\vec{i}:=(\vec{i}_1,\vec{i}_2,\cdots,\vec{i}_{d_0})\s\mbox{and}\s\vec{j}:=(\vec{j}_1,\vec{j}_2,\cdots,\vec{j}_{d_0}),
$$

$$
\vec{l}:=(\vec{l}_1,\vec{l}_2,\cdots,\vec{l}_{d_0})\s\mbox{and}\s\vec{v}:=(\vec{v}_1,\vec{v}_2,\cdots,\vec{v}_{d_0}),
$$

$$
\vec{i}_k:=(i^1_k,i^2_k,\cdots,i^{m_k}_k)\in\mathbb{Z}^{m_k}\s\mbox{and}\s\vec{j}_k:=(j^1_k,j^2_k,\cdots,j^{n_k}_k)\in\mathbb{Z}^{n_k},
$$

$$
\vec{l}_k:=(l^1_k,l^2_k,\cdots,l^{m_k}_k)\in\mathbb{Z}^{m_k}\s\mbox{and}\s\vec{v}_k:=(v^1_k,v^2_k,\cdots,v^{n_k}_k)\in\mathbb{Z}^{n_k}
$$
for $1\leq k\leq d_0$. Here, the coefficient is from this expression
$$
(u^1+u^2)^{\mathbf{p}}:=\sum\limits_{\mathbf{m}+\mathbf{n}=\mathbf{p}}E_{\mathbf{mn}}(u^1)^{\mathbf{m}}(u^2)^{\mathbf{n}}\s\s\mbox{with}\s E_{\mathbf{mn}}\in\mathbb{N}, u^1,u^2\in\mathbb{C}^{d_0}
$$
for
$$
\mathbf{m}:=(m_1,m_2,\cdots,m_{d_0})\s\mbox{and}\s\mathbf{n}:=(n_1,n_2,\cdots,n_{d_0}).
$$

Define
\begin{eqnarray*}
\sum_1:=\sum\limits_{\vec{i}_1=i-1}^{i+1}\sum\limits_{\vec{i}_2=i-1}^{i+1}\cdots\sum\limits_{\vec{i}_{d_0}=i-1}^{i+1}\sum\limits_{\vec{j}_1=i-1}^{i+1}\sum\limits_{\vec{j}_2=i-1}^{i+1}\cdots\sum\limits_{\vec{j}_{d_0}=i-1}^{i+1}
\end{eqnarray*}
with
\begin{eqnarray*}
\sum\limits_{\vec{i}_k=i-1}^{i+1}:=\sum\limits_{i^1_k=i-1}^{i+1}\sum\limits_{i^2_k=i-1}^{i+1}\cdots\sum\limits_{i^{m_k}_k=i-1}^{i+1}\s\mbox{and}\s\sum\limits_{\vec{j}_k=i-1}^{i+1}:=\sum\limits_{j^1_k=i-1}^{i+1}\sum\limits_{j^2_k=i-1}^{i+1}\cdots\sum\limits_{j^{n_k}_k=i-1}^{i+1}.
\end{eqnarray*}

Besides, we define
$$
B^1_{\vec{l}\vec{v}}:=\sup\limits_{i\in\mathbb{Z}}\sum_1|C_{\vec{i}\vec{j}\vec{l}\vec{v}}|,\s B^2_{\vec{l}\vec{v}}:=\sup\limits_{i\in\mathbb{Z}}\sum_1|\gamma^0C_{\vec{i}\vec{j}\vec{l}\vec{v}}|\s\mbox{and}\s B^3_{\vec{l}\vec{v}}:=\sup\limits_{i\in\mathbb{Z}}\sum_1\sum\limits_{\dot{k}=1}^d|\gamma^{\dot{k}} C_{\vec{i}\vec{j}\vec{l}\vec{v}}|.
$$
For $q=1,2,3$, we set
$$
B^q_{\mathbf{mn}}:=\sup\limits_{\vec{l}_1}\sup\limits_{\vec{l}_2}\cdots\sup\limits_{\vec{l}_{d_0}}\sup\limits_{\vec{v}_1}\sup\limits_{\vec{v}_2}\cdots\sup\limits_{\vec{v}_{d_0}}B^q_{\vec{l}\vec{v}}
$$
with
$$
\sup\limits_{\vec{l}_k}:=\sup\limits_{l^1_k}\cdots\sup\limits_{l^{m_k}_k}\s\mbox{and}\s\sup\limits_{\vec{v}_k}:=\sup\limits_{v^1_k}\cdots\sup\limits_{v^{n_k}_k}.
$$
At last, define
$$
B_{r}:=\sup\limits_{q=1,2,3}\sup\limits_{|\mathbf{m}|+|\mathbf{n}|=r}B^q_{\mathbf{mn}}.
$$
\textbf{Case for $A_p$}

Let
\begin{eqnarray*} C^i_{\varsigma\overrightarrow{wvfb\tau\iota\sigma\varrho}\zeta\chi}:=C^i_{\varsigma\overrightarrow{wvfb\tau\iota\sigma\varrho}\vartheta\varepsilon}:=a_{\mathbf{k}_i\mathbf{l}_i}b_{\mathbf{r}_i\mathbf{s}_i}p_i\frac{c_{\mathbf{m}_i\mathbf{n}_i}}{|\mathbf{m}_i|+1}\mathbf{c_p}
\end{eqnarray*}
with
$$
\overrightarrow{wvfb\tau\iota\sigma\varrho}:=(\vec{w},\vec{v},\vec{f},\vec{b},\vec{\tau},\vec{\iota},\vec{\sigma},\vec{\varrho}),
$$
where
$$
\vec{w}:=(\vec{w}_1,\vec{w}_2,\cdots,\vec{w}_{d_0})\s\mbox{with}\s \vec{w}_t:=(w_t^1,w_t^2,\cdots,w_t^{k_i^t}),
$$
$$
\vec{v}:=(\vec{v}_1,\vec{v}_2,\cdots,\vec{v}_{d_0})\s\mbox{with}\s \vec{v}_t:=(v_t^1,v_t^2,\cdots,v_t^{k_i^t}),
$$
$$
\vec{f}:=(\vec{f}_1,\vec{f}_2,\cdots,\vec{f}_{d_0})\s\mbox{with}\s \vec{f}_t:=(f_t^1,f_t^2,\cdots,f_t^{l_i^t}),
$$
$$
\vec{b}:=(\vec{b}_1,\vec{b}_2,\cdots,\vec{b}_{d_0})\s\mbox{with}\s \vec{b}_t:=(b_t^1,b_t^2,\cdots,b_t^{l_i^t}),
$$
$$
\vec{\tau}:=(\vec{\tau}_1,\vec{\tau}_2,\cdots,\vec{\tau}_{d_0})\s\mbox{with}\s \vec{\tau}_t:=(\tau_t^1,\tau_t^2,\cdots,\tau_t^{r_i^t}),
$$
$$
\vec{\iota}:=(\vec{\iota}_1,\vec{\iota}_2,\cdots,\vec{\iota}_{d_0})\s\mbox{with}\s \vec{\iota}_t:=(\iota_t^1,\iota_t^2,\cdots,\iota_t^{r_i^t}),
$$
$$
\vec{\sigma}:=(\vec{\sigma}_1,\vec{\sigma}_2,\cdots,\vec{\sigma}_{d_0})\s\mbox{with}\s \vec{\sigma}_t:=(\sigma_t^1,\sigma_t^2,\cdots,\sigma_t^{s_i^t})
$$
and
$$
\vec{\varrho}:=(\vec{\varrho}_1,\vec{\varrho}_2,\cdots,\vec{\varrho}_{d_0})\s\mbox{with}\s \vec{\varrho}_t:=(\varrho_t^1,\varrho_t^2,\cdots,\varrho_t^{s_i^t}).
$$
Here, the coefficients are from these expressions
$$
(u^1+u^2)^{(\mathbf{p}-e_i)^+}:=\sum\limits_{\mathbf{m}_i+\mathbf{n}_i=(\mathbf{p}-e_i)^+}c_{\mathbf{m}_i\mathbf{n}_i}(u^1)^{\mathbf{m}_i}(u^2)^{\mathbf{n}_i}\s\s\mbox{with}\s c_{\mathbf{m}_i\mathbf{n}_i}\in\mathbb{N}
$$
$$
(u^1+u^2)^{\mathbf{m}_i}:=\sum\limits_{\mathbf{k}_i+\mathbf{l}_i=\mathbf{m}_i}a_{\mathbf{k}_i\mathbf{l}_i}(u^1)^{\mathbf{k}_i}(u^2)^{\mathbf{l}_i}\s\s\mbox{with}\s a_{\mathbf{k}_i\mathbf{l}_i}\in\mathbb{N}
$$
and
$$
(u^1+u^2)^{\mathbf{n}_i}:=\sum\limits_{\mathbf{r}_i+\mathbf{s}_i=\mathbf{n}_i}b_{\mathbf{r}_i\mathbf{s}_i}(u^1)^{\mathbf{r}_i}(u^2)^{\mathbf{s}_i}\s\s\mbox{with}\s b_{\mathbf{r}_i\mathbf{s}_i}\in\mathbb{N}
$$
where $(\mathbf{p}-e_i)^+$ means a $d_0$-dimensional vector whose the $i$-th component is $\max\{p_i-1,0\}$ and the others are $p_k$ for $k\not=i$.

Besides, we appoint
\begin{eqnarray*}
\sum_2:=\sum\limits_{\vec{w}=\varsigma-1}^{\varsigma+1}\sum\limits_{\vec{f}=\varsigma-1}^{\varsigma+1}\sum\limits_{\vec{\tau}=\varsigma-1}^{\varsigma+1}\sum\limits_{\vec{\sigma}=\varsigma-1}^{\varsigma+1}\sum\limits_{\vartheta=\varsigma-1}^{\varsigma+1},
\end{eqnarray*}
where
\begin{eqnarray*}
\sum\limits_{\vec{w}=\varsigma-1}^{\varsigma+1}:=\sum\limits_{\vec{w}_1=\varsigma-1}^{\varsigma+1}\sum\limits_{\vec{w}_2=\varsigma-1}^{\varsigma+1}\cdots\sum\limits_{\vec{w}_{d_0}=\varsigma-1}^{\varsigma+1}\s\mbox{with}\s\sum\limits_{\vec{w}_t=\varsigma-1}^{\varsigma+1}:=\sum\limits_{w^1_t=\varsigma-1}^{\varsigma+1}\sum\limits_{w^2_t=\varsigma-1}^{\varsigma+1}\cdots\sum\limits_{w^{k_i^t}_t=\varsigma-1}^{\varsigma+1}
\end{eqnarray*}
and
\begin{eqnarray*}
\sum\limits_{\vec{f}=\varsigma-1}^{\varsigma+1}:=\sum\limits_{\vec{f}_1=\varsigma-1}^{\varsigma+1}\sum\limits_{\vec{f}_2=\varsigma-1}^{\varsigma+1}\cdots\sum\limits_{\vec{f}_{d_0}=\varsigma-1}^{\varsigma+1}\s\mbox{with}\s\sum\limits_{\vec{f}_t=\varsigma-1}^{\varsigma+1}:=\sum\limits_{f^1_t=\varsigma-1}^{\varsigma+1}\sum\limits_{f^2_t=\varsigma-1}^{\varsigma+1}\cdots\sum\limits_{f^{l_i^t}_t=\varsigma-1}^{\varsigma+1}
\end{eqnarray*}
and
\begin{eqnarray*}
\sum\limits_{\vec{\tau}=\varsigma-1}^{\varsigma+1}:=\sum\limits_{\vec{\tau}_1=\varsigma-1}^{\varsigma+1}\sum\limits_{\vec{\tau}_2=\varsigma-1}^{\varsigma+1}\cdots\sum\limits_{\vec{\tau}_{d_0}=\varsigma-1}^{\varsigma+1}\s\mbox{with}\s\sum\limits_{\vec{\tau}_t=\varsigma-1}^{\varsigma+1}:=\sum\limits_{\tau^1_t=\varsigma-1}^{\varsigma+1}\sum\limits_{\tau^2_t=\varsigma-1}^{\varsigma+1}\cdots\sum\limits_{\tau^{r_i^t}_t=\varsigma-1}^{\varsigma+1}
\end{eqnarray*}
and
\begin{eqnarray*}
\sum\limits_{\vec{\sigma}=\varsigma-1}^{\varsigma+1}:=\sum\limits_{\vec{\sigma}_1=\varsigma-1}^{\varsigma+1}\sum\limits_{\vec{\sigma}_2=\varsigma-1}^{\varsigma+1}\cdots\sum\limits_{\vec{\sigma}_{d_0}=\varsigma-1}^{\varsigma+1}\s\mbox{with}\s\sum\limits_{\vec{\sigma}_t=\varsigma-1}^{\varsigma+1}:=\sum\limits_{\sigma^1_t=\varsigma-1}^{\varsigma+1}\sum\limits_{\sigma^2_t=\varsigma-1}^{\varsigma+1}\cdots\sum\limits_{\sigma^{s_i^t}_t=\varsigma-1}^{\varsigma+1}.
\end{eqnarray*}
Moreover, we define
$$
B^i_{\overrightarrow{vb\iota\varrho}\varepsilon}:=\sup\limits_{\varsigma\in\mathbb{Z}}\sum_2\left\{\Big|C^i_{\varsigma\overrightarrow{wvfb\tau\iota\sigma\varrho}\vartheta\varepsilon}\Big|+\Big|\gamma^0C^i_{\varsigma\overrightarrow{wvfb\tau\iota\sigma\varrho}\vartheta\varepsilon}\Big|+\sum\limits_{\dot{k}=1}^d\Big|\gamma^{\dot{k}} C^i_{\varsigma\overrightarrow{wvfb\tau\iota\sigma\varrho}\vartheta\varepsilon}\Big|\right\}
$$
and
$$
B^i_{\mathbf{klrs}}:=\sup\limits_{\vec{v}}\sup\limits_{\vec{b}}\sup\limits_{\vec{\iota}}\sup\limits_{\vec{\varrho}}\sup\limits_{\varepsilon}B^i_{\overrightarrow{vb\iota\varrho}\varepsilon}
$$
with
$$
\sup\limits_{\vec{v}}:=\sup\limits_{\vec{v}_1}\sup\limits_{\vec{v}_2}\cdots\sup\limits_{\vec{v}_{d_0}}\s\mbox{and}\s \sup\limits_{\vec{b}}:=\sup\limits_{\vec{b}_1}\sup\limits_{\vec{b}_2}\cdots\sup\limits_{\vec{b}_{d_0}}
$$
and
$$
\sup\limits_{\vec{\varrho}}:=\sup\limits_{\vec{\varrho}_1}\sup\limits_{\vec{\varrho}_2}\cdots\sup\limits_{\vec{\varrho}_{d_0}}\s\mbox{and}\s \sup\limits_{\vec{\iota}}:=\sup\limits_{\vec{\iota}_1}\sup\limits_{\vec{\iota}_2}\cdots\sup\limits_{\vec{\iota}_{d_0}},
$$
where
$$
\sup\limits_{\vec{v}_t}:=\sup\limits_{v^1_t}\sup\limits_{v^2_t}\cdots\sup\limits_{v^{k_i^t}_t}\s\mbox{and}\s \sup\limits_{\vec{b}_t}:=\sup\limits_{b^1_t}\sup\limits_{b^2_t}\cdots\sup\limits_{b^{l_i^t}_t}
$$
and
$$
\sup\limits_{\vec{\iota}_t}:=\sup\limits_{\iota^1_t}\sup\limits_{\iota^2_t}\cdots\sup\limits_{\iota^{r_i^t}_t}\s\mbox{and}\s \sup\limits_{\vec{\varrho}_t}:=\sup\limits_{\varrho^1_t}\sup\limits_{\varrho^2_t}\cdots\sup\limits_{\varrho^{s_i^t}_t}.
$$
At last, set
$$
B^i_r:=\sup\limits_{|\mathbf{k}|+|\mathbf{l}|+|\mathbf{r}|+|\mathbf{s}|=r}B^i_{\mathbf{klrs}}\s\mbox{and}\s A_r:=\max\{B^i_r|1\leq i\leq d_0\}.
$$

\begin{thm}\label{thm3}
Given $d\geq 9$, there exist a small positive real number $\epsilon_0$ such that if for any $\epsilon\in(0,\epsilon_0]$ and any
$$
||\psi_0||_{H^{d/2}(\mathbb{T}^d,\mathbb{C}^{d_0})}\leq\epsilon,
$$
the Cauchy problem $(\ref{1})$ admits a unique global solution in the space $Y^{d/2}$, provided $F$ is analytic in a small neighborhood of $0$, $F(0)=0$ and meets the Condition (A).
\end{thm}

\begin{rem}
The definition of $Y^{d/2}$ is given in Section \ref{section4}.
\end{rem}
\textbf{The Organization}

The present article is organized as follow: Section \ref{section2} gives some preliminaries and notations about function spaces and frequency cut-off operators relying upon Fourier Series; we get almost global existence for high regularity in Section \ref{section3}; the main part of this article, solving the problem of low regularity, focuses on Section \ref{section4}; in the Appendix, we prove the projection operator $\Pi_{\pm}(D)$ is bounded.

\section{Preliminaries and Notations}\label{section2}
Set $A\lesssim B$ to denote $A\leq eB$ for a absolute constant $e>0$.

For a multi-index $\mathbf{p}:=(p_1,p_2,\cdots,p_{d_0})$, we set $|\mathbf{p}|:=\sum\limits_{k=1}^{d_0}p_k$. Besides we also define
\begin{eqnarray*}
\mathcal{F}_{x\mapsto\xi}u(\xi):=\frac{1}{(2\pi)^d}\int_{\mathbb{T}^d}e^{-\i x\cdot\xi}u(x)\,dx\s\s\mbox{for}\s\s\xi\in\mathbb{Z}^d\s\mbox{and}\s u\in L^1(\mathbb{T}^d),
\end{eqnarray*}
\begin{eqnarray*}
\mathcal{F}_{(t,x)\mapsto(\tau,\xi)}u(\tau,\xi):=\frac{1}{(2\pi)^{d+1}}\int_{\mathbb{R}}\int_{\mathbb{T}^d}e^{-\i x\cdot\xi-\i t\tau}u(t,x)\,dtdx
\end{eqnarray*}
for
\begin{eqnarray*}
(\tau,\xi)\in\mathbb{R}\times\mathbb{Z}^d\s\mbox{and}\s u\in L^1(\mathbb{R}\times\mathbb{T}^d),
\end{eqnarray*}

\begin{eqnarray*}
\mathcal{F}^{-1}_{\xi\mapsto x}\phi(x):=\sum\limits_{\xi\in\mathbb{Z}^d}e^{\i x\cdot\xi}\phi(\xi)\s\s\mbox{for}\s\s x\in\mathbb{T}^d\s\mbox{and}\s \phi\in l^1(\mathbb{Z}^d):=\left\{\phi\Bigg|\s\sum\limits_{\xi\in\mathbb{Z}^d}|\phi(\xi)|<\infty\right\}
\end{eqnarray*}
and
\begin{eqnarray*}
\mathcal{F}^{-1}_{(\tau,\xi)\mapsto(t,x)}\phi(t,x):=\sum\limits_{\xi\in\mathbb{Z}^d}e^{\i x\cdot\xi}\int_{\mathbb{R}}e^{\i t\tau}\phi(\tau,\xi)\,d\tau
\end{eqnarray*}
for
$$
(t,x)\in(\mathbb{R},\mathbb{T}^d)\s\mbox{and}\s \phi\in L^1(\mathbb{R},l^1(\mathbb{Z}^d)):=\left\{\phi\Bigg|\s\sum\limits_{\xi\in\mathbb{Z}^d}\int_{\mathbb{R}}|\phi(\tau,\xi)|\,d\tau<\infty\right\}.
$$
Set
$$
(\phi_1\ast\phi_2)(\xi):=\sum\limits_{\zeta\in\mathbb{Z}^d}\phi_1(\xi-\zeta)\phi_2(\zeta)
$$
and
$$
(\phi_1\ast\phi_2)(\tau,\xi):=\sum\limits_{\zeta\in\mathbb{Z}^d}\int_{\mathbb{R}}\phi_1(\tau-t,\xi-\zeta)\phi_2(t,\zeta)\,dt.
$$
Fix a cut-off function $\rho\in C_0^{\infty}(\mathbb{R})$ with compact support $[-2,2]$, $0\leq\rho\leq1$ and $\rho(x)\equiv1$ on $[-1,1]$. Define $\varphi(\xi):=\rho(|\xi|/2)-\rho(|\xi|)$,  $\tilde{\varphi}_j(\xi):=\varphi(2^{-j+1}\xi)+\varphi(2^{-j}\xi)+\varphi(2^{-j-1}\xi)$.

Denote
$$P_j:=\mathcal{F}^{-1}_{\xi\mapsto x}\circ\varphi(2^{-j}\cdot)\circ\mathcal{F}_{x\mapsto\xi}\s\mbox{and}\s\tilde{P}_j:=\mathcal{F}^{-1}_{\xi\mapsto x}\circ\tilde{\varphi}_j\circ\mathcal{F}_{x\mapsto\xi}.$$
Let $\mathcal{K}_l$ be a collection of spherical caps of diameter $2^{-l}$ which provide a symmetric and finitely overlapping cover of the unit sphere $\mathbb{S}^{d-1}$. $\dot{\Gamma}_{\kappa}$ means the cone generated by $\kappa\in\mathcal{K}_l$ and the origin.

Further, let $\eta_{\kappa}$ be smooth partition of unity subordinate to the covering of $\mathbb{R}^d\setminus\{0\}$ with the cones $\dot{\Gamma}_{\kappa}$ such that each $\eta_{\kappa}$ is supported in $2\dot{\Gamma}_{\kappa}$. We define
$$P_{\kappa}:=\mathcal{F}^{-1}_{\xi\mapsto x}\circ\eta_{\kappa}\circ\mathcal{F}_{x\mapsto\xi}.$$

Set
$$
Q_j^{\pm}u(t,x):=\mathcal{F}^{-1}_{(\tau,\xi)\mapsto(t,x)}\{[\rho(2^{-j-1}|\tau\pm\langle\xi\rangle|)-\rho(2^{-j}|\tau\pm\langle\xi\rangle|)]\mathcal{F}_{(t,x)\mapsto(\tau,\xi)}u\}.
$$

The next matrix will be used widely
$$
\Pi_{\pm}(\xi):=\frac{1}{2}\left\{I_{d_0}\pm\left(\sum\limits_{j=1}^d\alpha^j\xi_j+\beta\right)/\langle\xi\rangle\right\}
$$
with $\langle\xi\rangle:=(1+|\xi|^2)^{1/2}$, $\alpha^j:=\gamma^0\gamma^j$ and $\beta:=\gamma^0$, where $I_{d_0}$ is the unit matrix of order $d_0$.

Given any positive integer $k$, let $\Xi_k:=2^k\cdot\mathbb{Z}^d$. Define $\gamma^{(0)}:\mathbb{R}\rightarrow[0,1]$ to be an even smooth function supported in $[-1,1]$ and $\gamma^{(1)}$ to be
\begin{eqnarray*}
\gamma^{(1)}(t):=\frac{\gamma^{(0)}(t)}{\sum\limits_{z\in\mathbb{Z}}\gamma^{(0)}(t-z)}.
\end{eqnarray*}
It is not difficult to see that $\gamma^{(1)}$ is supported in $[-1,1]$ and
\begin{eqnarray*}
\sum\limits_{z\in\mathbb{Z}}\gamma^{(1)}(t-z)=1.
\end{eqnarray*}
Let $\gamma:\mathbb{R}^d\rightarrow[0,1]$, $\gamma(\xi):=\gamma^{(1)}(\xi_1)\cdots\gamma^{(1)}(\xi_d)$. For $n\in\Xi_k$, let
$$
\gamma_{k,n}(\xi):=\gamma((\xi-n)/2^k)\s\s\mbox{and}\s\s\Gamma_{k,n}:=\mathcal{F}^{-1}_{\xi\mapsto x}\circ\gamma_{k,n}(\xi)\circ\mathcal{F}_{x\mapsto\xi}.
$$

\subsection{Function spaces}
For $1\leq p\leq\infty$, $b'\in\mathbb{R}$, we define
$$
||f||_{\dot{X}^{\pm,b',p}}:=\Big|\Big|\left(2^{b'j}\Big|\Big|Q_j^{\pm}f\Big|\Big|_{L_t^2L_x^2}\right)_{j\in\mathbb{Z}}\Big|\Big|_{l^p}
$$
and
$$
||f||_{L^p_tL_x^q[l,k']}:=\sum\limits_{\kappa\in\mathcal{K}_l}\sum\limits_{n\in\Xi_{k'}}||\Gamma_{k',n}P_{\kappa}f||_{L_t^pL_x^q}.
$$
Moreover, we set
\begin{eqnarray*}
||f||_{S^{\pm}_j}:&=&||f||_{L_t^{\infty}L_x^2}+||f||_{\dot{X}^{\pm,1/2,\infty}}\\
&&+\sup\limits_{\substack{0\leq k'\leq j\\(d+2)j/(2d-2)\leq l\leq j}}\left\{2^{-(k'+j)/a}||f||_{L_t^aL_x^b[l,k']}+2^{-(k'+j)/b}||f||_{L_t^bL_x^a[l,k']}\right\}
\end{eqnarray*}
with $a,b$ to be determined. Define
\begin{eqnarray*}
||f||_{S^{\pm,\sigma}}:=\sum\limits_{j\geq0}2^{\sigma j}||P_jf||_{S^{\pm}_{j+1}}.
\end{eqnarray*}

\begin{lem}\label{lem2}
Given $k\in\mathbb{Z}$ and $j,k'\in\mathbb{N}$, we have
\begin{eqnarray*}
\sup\limits_{0\leq l\leq j}\left\{2^{-\frac{k'+j}{a}}||Q^{\pm}_kf||_{L^a_tL^b_x[l,k']}\right\}\lesssim||f||_{S^{\pm}_j}
\end{eqnarray*}
with $1/a+1/b=1/2$.
\end{lem}
\textbf{Proof.} The details are similar to the proof of Part $\textit{iii})$ of Lemma 3.3 in \cite{BH}.
\endproof

\section{Proof of Theorem \ref{thm1} and Theorem \ref{thm2}}\label{section3}
Assume $\psi$ is a second-order continuous differentiable solution to (\ref{1}). Applying $\partial/\partial_t$ to both sides of (\ref{1}) leads to
\begin{eqnarray}\label{2}
\frac{\partial^2\psi}{\partial t^2}-\Delta\psi+m^2\psi=G(\psi,\partial\psi)
\end{eqnarray}
with
$$
G(\psi,\partial\psi):=mF+\i\sum\limits_{j=1}^d\gamma^j\frac{\partial F}{\partial\psi}\frac{\partial\psi}{\partial x^j}-\i\sum\limits_{j=1}^d\gamma^0\frac{\partial F}{\partial\psi}\gamma^0\gamma^j\frac{\partial\psi}{\partial x^j}+m\gamma^0\frac{\partial F}{\partial\psi}\gamma^0\psi-\gamma^0\frac{\partial F}{\partial\psi}\gamma^0 F
$$
and
\begin{eqnarray}\label{3}
\psi(0,\cdot)=\psi_0,\s\s\frac{\partial\psi}{\partial t}(0,\cdot)=-\gamma^0\sum\limits_{j=1}^d\gamma^j\frac{\partial\psi_0}{\partial x^j}-\i m\gamma^0\psi_0+\i\gamma^0 F(\psi_0).
\end{eqnarray}

On the other hand, we take into account of the classical solution $\psi$ to (\ref{2}) with the initial datum (\ref{3}). Define
$$
V(t):=\int_{\mathbb{T}^3}\Big|\i\sum\limits_{\mu=0}^d\gamma^{\mu}\frac{\partial\psi}{\partial x^{\mu}}-m\psi+F(\psi)\Big|^2.
$$
It is easy to see $V(0)=0$ and $V'(t)\equiv0$ since of (\ref{2}). Hence $V(t)\equiv0$ as long as the solution $\psi$ exists.

In conclusion, we have the following theorem
\begin{thm}
A second-order continuous differentiable function solves the Dirac equation (\ref{1}) if and only if it solves the Klein-Gordon type equation (\ref{2}).
\end{thm}

However, Delort and Szeftel have given the existence and uniqueness for Klein-Gordon equation with high regularity on the standard $d$-dimensional torus $\mathbb{T}^d$.
\begin{thm}(\cite{DS})
Given $s>(9/2)d+8$, there are $\epsilon_0>0$ and $c>0$ such that for any pair $(u_0, u_1)$ in the unit ball of $H^{s+1}(\mathbb{T}^d)\times H^s(\mathbb{T}^d)$ and any $\epsilon\in(0,\epsilon_0]$, the following system
\begin{eqnarray*}
\left\{ \begin{aligned}
&D_t^2u-\sum\limits_{j,k=1}^dg^{jk}(u,\partial_tu,\partial u)D_jD_ku-u=\tilde{G}(u,\partial_tu,\partial u) \\
&u|_{t=0}=\epsilon u_0,\s\s\partial_tu|_{t=0}=\epsilon u_1
                          \end{aligned} \right.
\end{eqnarray*}
admits a unique solution
$$
u\in C^0((-T_{\epsilon},T_{\epsilon}),H^{s+1}(\mathbb{T}^d))\bigcap C^1((-T_{\epsilon},T_{\epsilon}),H^s(\mathbb{T}^d))
$$
with $g^{jk}:=\delta^{jk}+a^{jk}$, $T_{\epsilon}\geq c\epsilon^{-2}$, $D_tu=-\i\partial u/\partial t$, $D_ju=-\i\partial u/\partial x^j$ and $D:=(D_1,D_2,\cdots,D_d)$, where $a^{jk}$ is a real symmetric matrix of smooth functions defined on a neighborhood of 0 in $\mathbb{R}\times\mathbb{R}\times\mathbb{R}^d$ and $\tilde{G}$ is a smooth real function defined near 0 in $\mathbb{R}^{d+2}$ vanishing at least at order 2 at 0.
\end{thm}

In our case, $a^{jk}$ is 0. Since $m$ can be canceled by scaling, Theorem \ref{thm1} holds.

Furthermore, we have $G(-\psi,-\partial\psi)=-G(\psi,\partial\psi)$ by the assumption of Theorem \ref{thm2} that $F$ is odd. Hence, from Theorem 2.1 in \cite{FHZ} it follows that Theorem \ref{thm2} is true.
\endproof

\begin{rem}
In \cite{FHZ} Fang, Han and Zhang deal with the case $d=1$ for semi-linear Klein-Gordon equation. But their methods are also suitable for $d>1$.
\end{rem}

\section{Proof for Theorem \ref{thm3}}\label{section4}
Canceling $m$ by scaling and employing the notations of (3.1) in \cite{BH2} we rewrite (\ref{1}) as
\begin{eqnarray*}
-\i\left(\frac{\partial\psi}{\partial t}+\sum\limits_{j=1}^d\alpha^j\frac{\partial\psi}{\partial x^j}+\i\beta\psi\right)=\beta F(\psi).
\end{eqnarray*}
Define
\begin{eqnarray*}
\Pi_{\pm}(D):=\mathcal{F}^{-1}_{\xi\mapsto x}\circ\Pi_{\pm}(\xi)\circ\mathcal{F}_{x\mapsto\xi}\s\mbox{and}\s\langle D\rangle:=\mathcal{F}^{-1}_{\xi\mapsto x}\circ\langle\xi\rangle\circ\mathcal{F}_{x\mapsto\xi}
\end{eqnarray*}
and $\psi_{\pm}=\Pi_{\pm}(D)\psi$. Then $\psi$ can be rewritten as $\psi=\psi_++\psi_-$. For simplicity, we use $\Pi_{\pm}(D)\psi_k$ to denote the $k$-th component of $\Pi_{\pm}(D)\psi$.

Following the idea of (3.4) in \cite{BH2} we rewrite the Dirac equations as
\begin{eqnarray*}
\left\{ \begin{aligned}
         &\i\frac{\partial \psi_+}{\partial t}+\langle D\rangle\psi_+=-\Pi_+(D)(\beta F(\psi_++\psi_-)),\s\s\psi_+(0)=\Pi_+(D)\psi_0  \\
         &\i\frac{\partial \psi_-}{\partial t}-\langle D\rangle\psi_-=-\Pi_-(D)(\beta F(\psi_++\psi_-)),\s\s\psi_-(0)=\Pi_-(D)\psi_0.
                          \end{aligned} \right.
\end{eqnarray*}

In the next, we are going to employ contraction mapping principle. Define an operator $\mathcal{T}$ from
$$Y^{d/2}:=\{(\psi,\varphi)\in S^{+,d/2}\times S^{-,d/2}|\s||\psi||_{S^{+,d/2}}+||\varphi||_{S^{-,d/2}}\leq\epsilon\leq 1\}$$
to itself(this claim need to be checked)
\begin{eqnarray*}
\left\{ \begin{aligned}
         \mathcal{T}\psi(t):&=e^{-\i t\langle D\rangle}\Pi_+(D)\psi_0+\i\int_0^t e^{\i(t-s)\langle D\rangle}\Pi_+(D)\{\beta F(\Pi_+(D)\psi+\Pi_-(D)\varphi)\}\\
                             :&=e^{-\i t\langle D\rangle}\Pi_+(D)\psi_0+\i\int_0^t e^{\i(t-s)\langle D\rangle}f_+(s)\,ds\s\s\mbox{for}\s\psi\in S^{+,d/2}\\
         \mathcal{T}\varphi(t):&=e^{ \i t\langle D\rangle}\Pi_-(D)\psi_0+\i\int_0^t e^{-\i(t-s)\langle D\rangle}\Pi_-(D)\{\beta F(\Pi_+(D)\psi+\Pi_-(D)\varphi)\}\\
                             :&=e^{ \i t\langle D\rangle}\Pi_-(D)\psi_0+\i\int_0^t e^{-\i(t-s)\langle D\rangle}f_-(s)\,ds\s\s\mbox{for}\s\varphi\in S^{-,d/2}.
                          \end{aligned} \right.
\end{eqnarray*}

We need the coming lemma.
\begin{lem}\label{lem1}
For any $j\in\mathbb{N}$, $u_0=\tilde{P}_ju_0\in L^2(\mathbb{T}^d,\mathbb{C}^{d_0})$ and $f=\tilde{P}_jf\in L^1_t(\mathbb{R},L^2(\mathbb{T}^d,\mathbb{C}^{d_0}))$, set
\begin{eqnarray*}
u(t):=e^{\pm\i t\langle D\rangle}u_0+\i\int_0^t e^{\mp\i(t-s)\langle D\rangle}f(s)\,ds.
\end{eqnarray*}
Then $u=\tilde{P}_ju$ is the unique solution of
$$
\i\partial_t u\pm\langle D\rangle u=f,
$$
and $u\in C(\mathbb{R},L^2(\mathbb{T}^d,\mathbb{C}^{d_0}))$ and
$$
||u||_{S^{\pm}_j}\lesssim ||u_0||_{L^2(\mathbb{T}^d,\mathbb{C}^{d_0})}+\sup\limits_{g\in G^{\pm}_j}\Big|\int_{\mathbb{R}}\int_{\mathbb{T}^d}\langle f,g\rangle_{\mathbb{C}^{d_0}}\,dxdt\Big|
$$
where $G^{\pm}_j$ is defined as the set of all $g=\tilde{P}_jg\in L_t^{\infty}(\mathbb{R},L^2(\mathbb{T}^d,\mathbb{C}^{d_0}))$ such that $||g||_{S^{\pm}_j}=1$.
\end{lem}
\textbf{Proof.} The method is the same as Lemma 3.4 in \cite{BH}, since that proof does not involve the geometric properties of domain.
\endproof

Our goal is to estimate $||\mathcal{T}\psi||_{S^{+,d/2}}+||\mathcal{T}\varphi||_{S^{-,d/2}}$. Indeed, we have
\begin{eqnarray*}
||\mathcal{T}\psi||_{S^{+,d/2}}\leq||e^{-\i t\langle D\rangle}\Pi_+(D)\psi_0||_{S^{+,d/2}}+\Big|\Big|\int_0^t e^{\i(t-s)\langle D\rangle}f_+(s)\,ds\Big|\Big|_{S^{+,d/2}}
\end{eqnarray*}
and
\begin{eqnarray*}
||\mathcal{T}\varphi||_{S^{-,d/2}}\leq||e^{\i t\langle D\rangle}\Pi_-(D)\psi_0||_{S^{-,d/2}}+\Big|\Big|\int_0^t e^{-\i(t-s)\langle D\rangle}f_-(s)\,ds\Big|\Big|_{S^{-,d/2}}.
\end{eqnarray*}
\subsection{Step 1: $\mathcal{T}$ is from $Y^{d/2}$ to itself}
\textbf{linear estimations}\\

From the Lemma \ref{lem1} it follows that
\begin{eqnarray*}
&&\sum\limits_{j\geq0}2^{dj/2}||e^{-\i t\langle D\rangle}P_j\Pi_+(D)\psi_0||_{S^+_{j+1}}\lesssim\sum\limits_{j\geq0}2^{dj/2}||P_j\Pi_+(D)\psi_0||_{L^2(\mathbb{T}^d,\mathbb{C}^{d_0})}\\
&\lesssim&||\Pi_+(D)\psi_0||_{B^{d/2}_{2,2}(\mathbb{T}^d,\mathbb{C}^{d_0})}\lesssim||\Pi_+(D)\psi_0||_{H^{d/2}(\mathbb{T}^d,\mathbb{C}^{d_0})}\lesssim||\psi_0||_{H^{d/2}(\mathbb{T}^d,\mathbb{C}^{d_0})}
\end{eqnarray*}
where $B^{d/2}_{2,2}(\mathbb{T}^d,\mathbb{C}^{d_0})$ is Besov space and we have used the fact $B^s_{2,2}=H^s$ which is referred to Exercise 5 in Section 9 of Chapter 13 of \cite{T}. Hence, we have
\begin{eqnarray*}
||e^{-\i t\langle D\rangle}\Pi_+(D)\psi_0||_{S^{+,d/2}}\lesssim||\psi_0||_{H^{d/2}(\mathbb{T}^d,\mathbb{C}^{d_0})}.
\end{eqnarray*}
Similar idea gives
\begin{eqnarray*}
||e^{\i t\langle D\rangle}\Pi_-(D)\psi_0||_{S^{-,d/2}}\lesssim||\psi_0||_{H^{d/2}(\mathbb{T}^d,\mathbb{C}^{d_0})}.
\end{eqnarray*}

\textbf{Nonlinear estimations}\\

The nonlinear term can be controlled by
\begin{eqnarray*}
\Big|\Big|\int_0^t e^{\pm\i(t-s)\langle D\rangle}f_{\pm}(s)\,ds\Big|\Big|_{S^{\pm,d/2}}\leq\sum\limits_{j=0}^{\infty}2^{dj/2}\Big|\Big|\int_0^t P_je^{\pm\i(t-s)\langle D\rangle}f_{\pm}(s)\,ds\Big|\Big|_{S_{j+1}^{\pm}}.
\end{eqnarray*}
Lemma \ref{lem1} implies
\begin{eqnarray*}
&&\Big|\Big|\int_0^t P_je^{\pm\i(t-s)\langle D\rangle}f_{\pm}(s)\,ds\Big|\Big|_{S_{j+1}^{\pm}}\lesssim\sup\limits_{g\in G^{\pm}_{j+1}}\Big|\int_{\mathbb{R}}\int_{\mathbb{T}^d}\langle P_jf_{\pm},g\rangle_{\mathbb{C}^{d_0}}\,dxdt\Big|\\
&\leq&\sup\limits_{\phi\not=0}\Big|\int_{\mathbb{R}}\int_{\mathbb{T}^d}\langle P_jf_{\pm},\tilde{P}_{j+1}\phi\rangle_{\mathbb{C}^{d_0}}\,dxdt\Big|/||\tilde{P}_{j+1}\phi||_{S^{\pm}_{j+1}}.
\end{eqnarray*}
Then it suffices to prove
\begin{eqnarray}\label{6}
&&\Big|\int_{\mathbb{R}}\int_{\mathbb{T}^d}\langle F(\Pi_+(D)\psi+\Pi_-(D)\varphi),\beta\Pi_{\pm}(D)P_j\phi\rangle_{\mathbb{C}^{d_0}}\,dxdt\Big|\\
&\lesssim& h(||\Pi_+(D)\psi||_{S^{+,d/2}},||\Pi_-(D)\varphi||_{S^{-,d/2}})2^{-yj}||\tilde{P}_{j+1}\phi||_{S^{\mp}_{j+1}}\nonumber
\end{eqnarray}
for some smooth function $h(t,\tau)$ which vanishes at $(0,0)$ and $y>d/2$ is to be determined, where we have used the fact
$$P_j\tilde{P}_{j+1}=\tilde{P}_{j+1}P_j=P_j.$$
Indeed, we only have to obtain
\begin{eqnarray*}
&&\sum\limits_{i\in\mathbb{Z}}\Bigg|\int_{\mathbb{R}}\int_{\mathbb{T}^d}\left\langle F\left(\Pi_+(D)\psi+\Pi_-(D)\varphi\right),\beta Q^{\mp}_i\Pi_{\pm}(D)P_j\phi\right\rangle_{\mathbb{C}^{d_0}}\,dxdt\Bigg|\\
&\lesssim& h(||\Pi_+(D)\psi||_{S^{+,d/2}},||\Pi_-(D)\varphi||_{S^{-,d/2}})2^{-yj}||\tilde{P}_{j+1}\phi||_{S^{\mp}_{j+1}}.
\end{eqnarray*}
As soon as we have (\ref{6}), combining the Appendix leads to that $\mathcal{T}$ maps $Y^{d/2}$ into itself.

Because of
\begin{eqnarray}\label{7}
F(\psi)=\sum\limits_{|\mathbf{p}|=1}^{\infty}\mathbf{c}_{\mathbf{p}}\psi^{\mathbf{p}}
\end{eqnarray}
we have
\begin{eqnarray*}
&&F(\Pi_+(D)\psi+\Pi_-(D)\varphi)=\sum\limits_{|\mathbf{p}|=1}^{\infty}\sum\limits_{\mathbf{m}+\mathbf{n}=\mathbf{p}}c_{\mathbf{m}\mathbf{n}}\{\Pi_+(D)\psi\}^{\mathbf{m}}\{\Pi_-(D)\varphi\}^{\mathbf{n}}\\
&=&\sum\limits_{|\mathbf{p}|=1}^{\infty}\sum\limits_{\mathbf{m}+\mathbf{n}=\mathbf{p}}c_{\mathbf{m}\mathbf{n}}\left\{\sum\limits_{l,i\in\mathbb{Z}}Q^+_iP_l\Pi_+(D)\psi\right\}^{\mathbf{m}}\left\{\sum\limits_{v,j\in\mathbb{Z}}Q^-_jP_v\Pi_-(D)\varphi\right\}^{\mathbf{n}}\\
&=&\sum\limits_{|\mathbf{p}|=1}^{\infty}\sum\limits_{\mathbf{m}+\mathbf{n}=\mathbf{p}}\sum\limits_{\substack{\vec{i}:=(\vec{i}_1,\vec{i}_2,\cdots,\vec{i}_{d_0})\\ \vec{j}:=(\vec{j}_1,\vec{j}_2,\cdots,\vec{j}_{d_0})}}\sum\limits_{\substack{\vec{l}:=(\vec{l}_1,\vec{l}_2,\cdots,\vec{l}_{d_0})\\ \vec{v}:=(\vec{v}_1,\vec{v}_2,\cdots,\vec{v}_{d_0})}}C_{\vec{i}\vec{j}\vec{l}\vec{v}}\cdot Q^+_{\vec{i}_1}P_{\vec{l}_1}\Pi_+(D)\psi_1\cdots Q^+_{\vec{i}_{d_0}}P_{\vec{l}_{d_0}}\Pi_+(D)\psi_{d_0}\\
&&\cdot Q^-_{\vec{j}_1}P_{\vec{v}_1}\Pi_-(D)\varphi_1\cdots Q^-_{\vec{j}_{d_0}}P_{\vec{v}_{d_0}}\Pi_-(D)\varphi_{d_0}
\end{eqnarray*}
with
$$
\mathbf{m}:=(m_1,m_2,\cdots,m_{d_0})\s\mbox{and}\s\mathbf{n}:=(n_1,n_2,\cdots,n_{d_0}),
$$

$$
\vec{i}_k:=(i^1_k,i^2_k,\cdots,i^{m_k}_k)\in\mathbb{Z}^{m_k}\s\mbox{and}\s\vec{j}_k:=(j^1_k,j^2_k,\cdots,j^{n_k}_k)\in\mathbb{Z}^{n_k},
$$

$$
\vec{l}_k:=(l^1_k,l^2_k,\cdots,l^{m_k}_k)\in\mathbb{Z}^{m_k}\s\mbox{and}\s\vec{v}_k:=(v^1_k,v^2_k,\cdots,v^{n_k}_k)\in\mathbb{Z}^{n_k},
$$

$$
Q^+_{\vec{i}_k}P_{\vec{l}_k}\Pi_+(D)\psi_k:=Q^+_{i^1_k}P_{l^1_k}\Pi_+(D)\psi_k\cdot Q^+_{i^2_k}P_{l^2_k}\Pi_+(D)\psi_k\cdots Q^+_{i^{m_k}_k}P_{l^{m_k}_k}\Pi_+(D)\psi_k
$$
and
$$
Q^-_{\vec{j}_k}P_{\vec{v}_k}\Pi_-(D)\varphi_k:=Q^-_{j^1_k}P_{v^1_k}\Pi_-(D)\varphi_k\cdot Q^-_{j^2_k}P_{v^2_k}\Pi_-(D)\varphi_k\cdots Q^-_{j^{n_k}_k}P_{v^{n_k}_k}\Pi_-(D)\varphi_k
$$
for $1\leq k\leq d_0$.

From Plancherel formula it follows that
\begin{eqnarray*}
&&\int_{\mathbb{R}}\int_{\mathbb{T}^d}\left\{Q^+_{\vec{i}_1}P_{\vec{l}_1}\Pi_+(D)\psi_1\cdots Q^+_{\vec{i}_{d_0}}P_{\vec{l}_{d_0}}\Pi_+(D)\psi_{d_0}\cdot Q^-_{\vec{j}_1}P_{\vec{v}_1}\Pi_-(D)\varphi_1\cdots Q^-_{\vec{j}_{d_0}}P_{\vec{v}_{d_0}}\Pi_-(D)\varphi_{d_0}\right\}\\
&&\left\langle C_{\vec{i}\vec{j}\vec{l}\vec{v}},\beta Q^{\mp}_i\Pi_{\pm}(D)P_j\phi\right\rangle_{\mathbb{C}^{d_0}}\\
&&=\sum\limits_{\xi\in\mathbb{Z}^d}\int_{\mathbb{R}}\mathfrak{F}_{\vec{i}_1\vec{l}_1}\ast\cdots\ast\mathfrak{F}_{\vec{i}_{d_0}\vec{l}_{d_0}}\ast\mathfrak{F}_{\vec{j}_1\vec{v}_1}\ast\cdots\ast\mathfrak{F}_{\vec{j}_{d_0}\vec{v}_{d_0}}\ast\langle C_{\vec{i}\vec{j}\vec{l}\vec{v}},\beta\mathfrak{F}\rangle_{\mathbb{C}^{d_0}}
\end{eqnarray*}
where
$$
\mathfrak{F}_{\vec{i}_k\vec{l}_k}:=\mathcal{F}_{(t,x)\mapsto(\tau,\xi)}\left(Q^+_{\vec{i}_k}P_{\vec{l}_k}\Pi_+(D)\psi_k\right),\s\mathfrak{F}_{\vec{j}_k\vec{v}_k}:=\mathcal{F}_{(t,x)\mapsto(\tau,\xi)}\left(Q^-_{\vec{j}_k}P_{\vec{v}_k}\Pi_-(D)\varphi_k\right)
$$
and
$$
\mathfrak{F}:=\mathcal{F}_{(t,x)\mapsto(\tau,\xi)}\left(Q^{\mp}_i\Pi_{\pm}(D)P_j\phi\right).
$$
Indeed, we easily have
\begin{eqnarray*}
\mathfrak{F}_{\vec{i}_k\vec{l}_k}=\mathcal{F}_{(t,x)\mapsto(\tau,\xi)}\left(Q^+_{i^1_k}P_{l^1_k}\Pi_+(D)\psi_k\right)\ast\cdots\ast\mathcal{F}_{(t,x)\mapsto(\tau,\xi)}\left(Q^+_{i^{m_k}_k}P_{i^{m_k}_k}\Pi_+(D)\psi_k\right),
\end{eqnarray*}
\begin{eqnarray*}
\mbox{supp}\left\{\mathcal{F}_{(t,x)\mapsto(\tau,\xi)}\left(Q^+_{i^r_k}P_{l^r_k}\Pi_+(D)\psi_k\right)\right\}\subset\{(\tau,\xi)|2^{i^r_k}\leq|\tau\pm\langle\xi\rangle|\leq2^{2+i^r_k},\s 2^{l^r_k}\leq|\xi|\leq2^{2+l^r_k}\}
\end{eqnarray*}
and
\begin{eqnarray*}
\mbox{supp}\mathfrak{F}\subset\{(\tau,\xi)|2^i\leq|\tau\mp\langle\xi\rangle|\leq2^{2+i},\s 2^j\leq|\xi|\leq2^{j+2}\}
\end{eqnarray*}
for $1\leq r\leq m_k$. The accurate expression of $\mathfrak{F}_{\vec{i}_k\vec{l}_k}$ is
\begin{eqnarray*}
\mathfrak{F}_{\vec{i}_k\vec{l}_k}=\int_{\sum\limits_{r=1}^{m_k}\tau_r=\tau}\sum\limits_{\sum\limits_{r=1}^{m_k}\xi_r=\xi}\prod\limits_{r=1}^{m_k}\varphi(2^{-i^r_k}\{\tau_r\mp\langle\xi_r\rangle\})\varphi(2^{-l^r_k}\xi_r)\Pi_+(\xi_r)\left\{\mathcal{F}_{(t,x)\mapsto(\tau,\xi)}\psi_k\right\}(\tau_r,\xi_r)
\end{eqnarray*}
which implies that if $|i^r_k-i^s_k|\geq2$ or $|l^r_k-l^s_k|\geq2$ for $r\not=s$ or $|i^r_k-l^s_k|\geq2$, $\mathfrak{F}_{\vec{i}_k\vec{l}_k}=0$. This also means that if
$$
||(\vec{i}_1,\vec{i}_2,\cdots,\vec{i}_{d_0},\vec{j}_1,\vec{j}_2\cdots,\vec{j}_{d_0},i)||\geq2
$$
or
$$
||(\vec{l}_1,\vec{l}_2,\cdots,\vec{l}_{d_0},\vec{v}_1,\vec{v}_2,\cdots,\vec{v}_{d_0},j)||\geq2,
$$
then
\begin{eqnarray}\label{10}
&&\int_{\mathbb{R}}\int_{\mathbb{T}^d}Q^+_{\vec{i}_1}P_{\vec{l}_1}\Pi_+(D)\psi_1\cdots Q^+_{\vec{i}_{d_0}}P_{\vec{l}_{d_0}}\Pi_+(D)\psi_{d_0}\cdot Q^-_{\vec{j}_1}P_{\vec{v}_1}\Pi_-(D)\varphi_1\cdots Q^-_{\vec{j}_{d_0}}P_{\vec{v}_{d_0}}\Pi_-(D)\varphi_{d_0}\nonumber\\
&&\left\langle C_{\vec{i}\vec{j}\vec{l}\vec{v}},\beta Q^{\mp}_i\Pi_{\pm}(D)P_j\phi\right\rangle_{\mathbb{C}^{d_0}}=0
\end{eqnarray}
where
\begin{eqnarray*}
&&||(\vec{i}_1,\vec{i}_2,\cdots,\vec{i}_{d_0},\vec{j}_1,\vec{j}_2\cdots,\vec{j}_{d_0},i)||:=\max\{|i_k^{\dot{m}}-i_t^{\tilde{m}}|,|i_k^{\dot{m}}-i|,|j_k^{\dot{n}}-j_t^{\tilde{n}}|,|j_k^{\dot{n}}-i|,|i_k^{\dot{m}}-j_t^{\tilde{n}}|,\\
&&|1\leq k,t\leq d_0,1\leq\tilde{n}\leq n_t,1\leq \dot{n}\leq n_k,1\leq\tilde{m}\leq m_t,1\leq \dot{m}\leq m_k\}
\end{eqnarray*}
and
\begin{eqnarray*}
&&||(\vec{l}_1,\vec{l}_2,\cdots,\vec{l}_{d_0},\vec{v}_1,\vec{v}_2,\cdots,\vec{v}_{d_0},j)||:=\max\{|l_k^{\dot{p}}-l_t^{\tilde{p}}|,|v_k^{\dot{q}}-v_t^{\tilde{q}}|,|l_k^{\dot{p}}-v_t^{\tilde{q}}|,|l_k^{\dot{p}}-j|,|v_k^{\dot{q}}-j|\\
&&|1\leq k,t\leq d_0,1\leq\tilde{q}\leq n_t,1\leq \dot{q}\leq n_k,1\leq\tilde{p}\leq m_t,1\leq \dot{p}\leq m_k\}.
\end{eqnarray*}
From Young's Inequality it follows that
\begin{eqnarray*}
&&\Bigg|\int_{\mathbb{R}}\int_{\mathbb{T}^d}Q^+_{\vec{i}_1}P_{\vec{l}_1}\Pi_+(D)\psi_1\cdots Q^+_{\vec{i}_{d_0}}P_{\vec{l}_{d_0}}\Pi_+(D)\psi_{d_0}\cdot Q^-_{\vec{j}_1}P_{\vec{v}_1}\Pi_-(D)\varphi_1\cdots Q^-_{\vec{j}_{d_0}}P_{\vec{v}_{d_0}}\Pi_-(D)\varphi_{d_0}\\
&\leq&\Big|\Big|Q^+_{\vec{i}_1}P_{\vec{l}_1}\Pi_+(D)\psi_1\cdots Q^+_{\vec{i}_{d_0}}P_{\vec{l}_{d_0}}\Pi_+(D)\psi_{d_0}\cdot Q^-_{\vec{j}_1}P_{\vec{v}_1}\Pi_-(D)\varphi_1\cdots Q^-_{\vec{j}_{d_0}}P_{\vec{v}_{d_0}}\Pi_-(D)\varphi_{d_0}\Big|\Big|_{L_t^{p_{\vec{i}\vec{j}\vec{l}\vec{v}ij}}L_x^{\tilde{p}_{\vec{i}\vec{j}\vec{l}\vec{v}ij}}}\\
&&\cdot\Big|\Big|\left\langle C_{\vec{i}\vec{j}\vec{l}\vec{v}},\beta Q^{\mp}_i\Pi_{\pm}(D)P_j\phi\right\rangle_{\mathbb{C}^{d_0}}\Big|\Big|_{L_t^{q_{\vec{i}\vec{j}\vec{l}\vec{v}ij}}L_x^{\tilde{q}_{\vec{i}\vec{j}\vec{l}\vec{v}ij}}}
\end{eqnarray*}
with
$$
1/p_{\vec{i}\vec{j}\vec{l}\vec{v}ij}+1/q_{\vec{i}\vec{j}\vec{l}\vec{v}ij}=1\s\s\mbox{and}\s\s 1/\tilde{p}_{\vec{i}\vec{j}\vec{l}\vec{v}ij}+1/\tilde{q}_{\vec{i}\vec{j}\vec{l}\vec{v}ij}=1
$$
to be determined. Young's inequality yields
\begin{eqnarray*}
&&\Big|\Big|Q^+_{\vec{i}_1}P_{\vec{l}_1}\Pi_+(D)\psi_1\cdots Q^+_{\vec{i}_{d_0}}P_{\vec{l}_{d_0}}\Pi_+(D)\psi_{d_0}\cdot Q^-_{\vec{j}_1}P_{\vec{v}_1}\Pi_-(D)\varphi_1\cdots Q^-_{\vec{j}_{d_0}}P_{\vec{v}_{d_0}}\Pi_-(D)\varphi_{d_0}\Big|\Big|_{L_t^{p_{\vec{i}\vec{j}\vec{l}\vec{v}ij}}L_x^{\tilde{p}_{\vec{i}\vec{j}\vec{l}\vec{v}ij}}}\\
&\leq&\prod\limits_{k=1}^{d_0}\prod\limits_{r_k=1}^{m_k}\Big|\Big|Q^+_{i^{r_k}_k}P_{l^{r_k}_k}\Pi_+(D)\psi_k\Big|\Big|_{L_t^{p_k^{r_k}}L_x^{\tilde{p}_k^{r_k}}}\cdot\prod\limits_{k=1}^{d_0}\prod\limits_{s_k=1}^{n_k}\Big|\Big|Q^-_{j^{s_k}_k}P_{v^{s_k}_k}\Pi_-(D)\varphi_k\Big|\Big|_{L_t^{q_k^{s_k}}L_x^{\tilde{q}_k^{s_k}}}
\end{eqnarray*}
with
$$
\sum\limits_{k=1}^{d_0}\sum\limits_{r_k=1}^{m_k}1/p_k^{r_k}+\sum\limits_{k=1}^{d_0}\sum\limits_{s_k=1}^{n_k}1/q_k^{s_k}=1/p_{\vec{i}\vec{j}\vec{l}\vec{v}ij}
$$
and
$$
\sum\limits_{k=1}^{d_0}\sum\limits_{r_k=1}^{m_k}1/\tilde{p}_k^{r_k}+\sum\limits_{k=1}^{d_0}\sum\limits_{s_k=1}^{n_k}1/\tilde{q}_k^{s_k}=1/\tilde{p}_{\vec{i}\vec{j}\vec{l}\vec{v}ij}
$$
to be determined. For any $q_1,q_2\in\mathbb{R}$, we define
$$
\lceil q_1,q_2\rceil:=\min\{z\in\mathbb{Z}|z\geq 1/q_1-1/q_2\}d+2d.
$$
From Gagliardo-Nirenberg-Moser estimates it follows that
\begin{eqnarray*}
&&\Big|\Big|Q^+_{\vec{i}_1}P_{\vec{l}_1}\Pi_+(D)\psi_1\cdots Q^+_{\vec{i}_{d_0}}P_{\vec{l}_{d_0}}\Pi_+(D)\psi_{d_0}\cdot Q^-_{\vec{j}_1}P_{\vec{v}_1}\Pi_-(D)\varphi_1\cdots Q^-_{\vec{j}_{d_0}}P_{\vec{v}_{d_0}}\Pi_-(D)\varphi_{d_0}\Big|\Big|_{L_t^{p_{\vec{i}\vec{j}\vec{l}\vec{v}ij}}L_x^{\tilde{p}_{\vec{i}\vec{j}\vec{l}\vec{v}ij}}}\\
&\lesssim&\sum\limits_{\substack{|\chi_{r_k}|=\lceil \hat{p}_k^{r_k},\tilde{p}_k^{r_k}\rceil\\|\alpha_{s_k}|=\lceil\hat{q}_k^{s_k}, \tilde{q}_k^{s_k}\rceil}}\prod\limits_{k=1}^{d_0}\prod\limits_{r_k=1}^{m_k}\Big\{\Big|\Big|Q^+_{i^{r_k}_k}P_{l^{r_k}_k}\Pi_+(D)\psi_k\Big|\Big|^{1-\theta_{r_k}}_{L_t^{p_k^{r_k}}L_x^{\hat{p}_k^{r_k}}}\cdot\Big|\Big|D^{\chi_{r_k}}Q^+_{i^{r_k}_k}P_{l^{r_k}_k}\Pi_+(D)\psi_k\Big|\Big|^{\theta_{r_k}}_{L_t^{p_k^{r_k}}L_x^{\hat{p}_k^{r_k}}}\Big\}\\
&&\cdot\prod\limits_{k=1}^{d_0}\prod\limits_{s_k=1}^{n_k}\left\{\Big|\Big|Q^-_{j^{s_k}_k}P_{v^{s_k}_k}\Pi_-(D)\varphi_k\Big|\Big|^{1-\lambda_{s_k}}_{L_t^{q_k^{s_k}}L_x^{\hat{q}_k^{s_k}}}\cdot\Big|\Big|D^{\alpha_{s_k}}Q^-_{j^{s_k}_k}P_{v^{s_k}_k}\Pi_-(D)\varphi_k\Big|\Big|^{\lambda_{s_k}}_{L_t^{q_k^{s_k}}L_x^{\hat{q}_k^{s_k}}}\right\}\\
\end{eqnarray*}
with
$$
1/p_k^{r_k}+1/\hat{p}_k^{r_k}=1/2\s\s\mbox{and}\s\s 1/q_k^{s_k}+1/\hat{q}_k^{s_k}=1/2,
$$
where
$$
\theta_{r_k}:=d/(\lceil\hat{p}_k^{r_k},\tilde{p}_k^{r_k}\rceil\hat{p}_k^{r_k})-d/(\lceil\hat{p}_k^{r_k},\tilde{p}_k^{r_k}\rceil\tilde{p}_k^{r_k})
$$
and
$$
\lambda_{s_k}:=d/(\lceil\hat{q}_k^{s_k},\tilde{q}_k^{s_k}\rceil\hat{q}_k^{s_k})-d/(\lceil\hat{q}_k^{s_k},\tilde{q}_k^{s_k}\rceil\tilde{q}_k^{s_k}).
$$
Using the same method leads to
\begin{eqnarray*}
&&\Big|\Big|\left\langle C_{\vec{i}\vec{j}\vec{l}\vec{v}},\beta Q^{\mp}_i\Pi_{\pm}(D)P_j\phi\right\rangle_{\mathbb{C}^{d_0}}\Big|\Big|_{L_t^{q_{\vec{i}\vec{j}\vec{l}\vec{v}ij}}L_x^{\tilde{q}_{\vec{i}\vec{j}\vec{l}\vec{v}ij}}}\\
&\lesssim&\sum\limits_{|\vartheta|=\lceil\hat{q}_{\vec{i}\vec{j}\vec{l}\vec{v}ij},\tilde{q}_{\vec{i}\vec{j}\vec{l}\vec{v}ij}\rceil}\Big|\Big|\left\langle C_{\vec{i}\vec{j}\vec{l}\vec{v}},\beta Q^{\mp}_i\Pi_{\pm}(D)P_j\phi\right\rangle_{\mathbb{C}^{d_0}}\Big|\Big|^{1-\sigma_{\vec{i}\vec{j}\vec{l}\vec{v}ij}}_{L_t^{q_{\vec{i}\vec{j}\vec{l}\vec{v}ij}}L_x^{\hat{q}_{\vec{i}\vec{j}\vec{l}\vec{v}ij}}}\\
&&\cdot\Big|\Big|\left\langle C_{\vec{i}\vec{j}\vec{l}\vec{v}},\beta D^{\vartheta}Q^{\mp}_i\Pi_{\pm}(D)P_j\phi\right\rangle_{\mathbb{C}^{d_0}}\Big|\Big|^{\sigma_{\vec{i}\vec{j}\vec{l}\vec{v}ij}}_{L_t^{q_{\vec{i}\vec{j}\vec{l}\vec{v}ij}}L_x^{\hat{q}_{\vec{i}\vec{j}\vec{l}\vec{v}ij}}}
\end{eqnarray*}
with
$$
\sigma_{\vec{i}\vec{j}\vec{l}\vec{v}ij}:=d/\left(\hat{q}_{\vec{i}\vec{j}\vec{l}\vec{v}ij}\lceil\hat{q}_{\vec{i}\vec{j}\vec{l}\vec{v}ij},\tilde{q}_{\vec{i}\vec{j}\vec{l}\vec{v}ij}\rceil\right)-d/\left(\tilde{q}_{\vec{i}\vec{j}\vec{l}\vec{v}ij}\lceil\hat{q}_{\vec{i}\vec{j}\vec{l}\vec{v}ij},\tilde{q}_{\vec{i}\vec{j}\vec{l}\vec{v}ij}\rceil\right)
$$
and
$$
1/\hat{q}_{\vec{i}\vec{j}\vec{l}\vec{v}ij}+1/q_{\vec{i}\vec{j}\vec{l}\vec{v}ij}=1/2.
$$
From the definition and Lemma \ref{lem2} it follows that
\begin{eqnarray}\label{13}
&&\Big|\Big|Q^+_{i_k^{r_k}}P_{l^{r_k}_k}\Pi_+(D)\psi_k\Big|\Big|_{L_t^{p_k^{r_k}}L_x^{\hat{p}_k^{r_k}}}\leq\Big|\Big|Q^+_{i_k^{r_k}}P_{l^{r_k}_k}\Pi_+(D)\psi_k\Big|\Big|_{L_t^{p_k^{r_k}}L_x^{\hat{p}_k^{r_k}}[j,j]}\nonumber\\
&\lesssim&2^{2j/p_k^{r_k}}||P_{l^{r_k}_k}\Pi_+(D)\psi_k||_{S_{j}^+}.
\end{eqnarray}
Similarly, we can obtain
\begin{eqnarray*}
\Big|\Big|Q^-_{j_k^{s_k}}P_{v^{s_k}_k}\Pi_-(D)\varphi_k\Big|\Big|_{L_t^{q_k^{s_k}}L_x^{\hat{q}_k^{s_k}}}\lesssim 2^{2j/q_k^{s_k}}||P_{v^{s_k}_k}\Pi_-(D)\varphi_k||_{S_{j}^-}.
\end{eqnarray*}
In conclusion, the next inequality holds true
\begin{eqnarray*}
&&\Big|\Big|Q^+_{\vec{i}_1}P_{\vec{l}_1}\Pi_+(D)\psi_1\cdots Q^+_{\vec{i}_{d_0}}P_{\vec{l}_{d_0}}\Pi_+(D)\psi_{d_0}\cdot Q^-_{\vec{j}_1}P_{\vec{v}_1}\Pi_-(D)\varphi_1\cdots Q^-_{\vec{j}_{d_0}}P_{\vec{v}_{d_0}}\Pi_-(D)\varphi_{d_0}\Big|\Big|_{L_t^{p_{\vec{i}\vec{j}\vec{l}\vec{v}ij}}L_x^{\tilde{p}_{\vec{i}\vec{j}\vec{l}\vec{v}ij}}}\\
&\lesssim&\sum\limits_{\substack{|\chi_{r_k}|=\lceil \hat{p}_k^{r_k},\tilde{p}_k^{r_k}\rceil\\|\alpha_{s_k}|=\lceil\hat{q}_k^{s_k}, \tilde{q}_k^{s_k}\rceil}}\prod\limits_{k=1}^{d_0}\prod\limits_{r_k=1}^{m_k}\left\{2^{2j/p_k^{r_k}}\Big|\Big|P_{l^{r_k}_k}\Pi_+(D)\psi_k\Big|\Big|_{S_{j}^+}^{1-\theta_{r_k}}\cdot\Big|\Big|D^{\chi_{r_k}}P_{l^{r_k}_k}\Pi_+(D)\psi_k\Big|\Big|_{S_{j}^+}^{\theta_{r_k}}\right\}\\
&&\cdot\prod\limits_{k=1}^{d_0}\prod\limits_{s_k=1}^{n_k}\left\{2^{2jq_k^{s_k}}\Big|\Big|P_{v^{s_k}_k}\Pi_-(D)\varphi_k\Big|\Big|^{1-\lambda_{s_k}}_{S_{j}^-}\cdot\Big|\Big|D^{\alpha_{s_k}}P_{v^{s_k}_k}\Pi_-(D)\varphi_k\Big|\Big|^{\lambda_{s_k}}_{S_{j}^-}\right\}\\
&=&\sum\limits_{\substack{|\chi_{r_k}|=\lceil \hat{p}_k^{r_k},\tilde{p}_k^{r_k}\rceil\\|\alpha_{s_k}|=\lceil\hat{q}_k^{s_k}, \tilde{q}_k^{s_k}\rceil}}2^{2j/p_{\vec{i}\vec{j}\vec{l}\vec{v}ij}}\prod\limits_{k=1}^{d_0}\prod\limits_{s_k=1}^{n_k}\left\{\Big|\Big|P_{v^{s_k}_k}\Pi_-(D)\varphi_k\Big|\Big|^{1-\lambda_{s_k}}_{S_{j}^-}\cdot\Big|\Big|D^{\alpha_{s_k}}P_{v^{s_k}_k}\Pi_-(D)\varphi_k\Big|\Big|^{\lambda_{s_k}}_{S_{j}^-}\right\}\\
&&\cdot\prod\limits_{k=1}^{d_0}\prod\limits_{r_k=1}^{m_k}\left\{\Big|\Big|P_{l^{r_k}_k}\Pi_+(D)\psi_k\Big|\Big|_{S_{j}^+}^{1-\theta_{r_k}}\cdot\Big|\Big|D^{\chi_{r_k}}P_{l^{r_k}_k}\Pi_+(D)\psi_k\Big|\Big|_{S_{j}^+}^{\theta_{r_k}}\right\}.
\end{eqnarray*}
Moreover, the same idea also implies
\begin{eqnarray*}
\Big|\Big|\left\langle C_{\vec{i}\vec{j}\vec{l}\vec{v}},\beta Q^{\mp}_i\Pi_{\pm}(D)P_j\phi\right\rangle_{\mathbb{C}^{d_0}}\Big|\Big|_{L_t^{q_{\vec{i}\vec{j}\vec{l}\vec{v}ij}}L_x^{\hat{q}_{\vec{i}\vec{j}\vec{l}\vec{v}ij}}}\leq2^{2j/q_{\vec{i}\vec{j}\vec{l}\vec{v}ij}}\Big|\Big|\left\langle C_{\vec{i}\vec{j}\vec{l}\vec{v}},\beta\Pi_{\pm}(D)P_j\phi\right\rangle_{\mathbb{C}^{d_0}}\Big|\Big|_{S_{j}^{\mp}}
\end{eqnarray*}
and
\begin{eqnarray*}
&&\sum\limits_{|\vartheta|=\lceil\hat{q}_{\vec{i}\vec{j}\vec{l}\vec{v}ij},\tilde{q}_{\vec{i}\vec{j}\vec{l}\vec{v}ij}\rceil}\Big|\Big|\left\langle C_{\vec{i}\vec{j}\vec{l}\vec{v}},\beta D^{\vartheta}Q^{\mp}_i\Pi_{\pm}(D)P_j\phi\right\rangle_{\mathbb{C}^{d_0}}\Big|\Big|_{L_t^{q_{\vec{i}\vec{j}\vec{l}\vec{v}ij}}L_x^{\hat{q}_{\vec{i}\vec{j}\vec{l}\vec{v}ij}}}\\
&\leq&\sum\limits_{|\vartheta|=\lceil\hat{q}_{\vec{i}\vec{j}\vec{l}\vec{v}ij},\tilde{q}_{\vec{i}\vec{j}\vec{l}\vec{v}ij}\rceil}2^{2j/q_{\vec{i}\vec{j}\vec{l}\vec{v}ij}}\Big|\Big|\left\langle C_{\vec{i}\vec{j}\vec{l}\vec{v}},\beta D^{\vartheta}\Pi_{\pm}(D)P_j\phi\right\rangle_{\mathbb{C}^{d_0}}\Big|\Big|_{S_{j}^{\mp}}.
\end{eqnarray*}
Hence we have
\begin{eqnarray*}
&&\Big|\Big|\left\langle C_{\vec{i}\vec{j}\vec{l}\vec{v}},\beta Q^{\mp}_i\Pi_{\pm}(D)P_j\phi\right\rangle_{\mathbb{C}^{d_0}}\Big|\Big|_{L_t^{q_{\vec{i}\vec{j}\vec{l}\vec{v}ij}}L_x^{\tilde{q}_{\vec{i}\vec{j}\vec{l}\vec{v}ij}}}\\
&\leq&\sum\limits_{|\vartheta|=\lceil\hat{q}_{\vec{i}\vec{j}\vec{l}\vec{v}ij},\tilde{q}_{\vec{i}\vec{j}\vec{l}\vec{v}ij}\rceil}2^{2j/q_{\vec{i}\vec{j}\vec{l}\vec{v}ij}}\Big|\Big|\left\langle C_{\vec{i}\vec{j}\vec{l}\vec{v}},\beta\Pi_{\pm}(D)P_j\phi\right\rangle_{\mathbb{C}^{d_0}}\Big|\Big|^{1-\sigma_{\vec{i}\vec{j}\vec{l}\vec{v}ij}}_{S_{j}^{\mp}}\\
&&\Big|\Big|\left\langle C_{\vec{i}\vec{j}\vec{l}\vec{v}},\beta D^{\vartheta}\Pi_{\pm}(D)P_j\phi\right\rangle_{\mathbb{C}^{d_0}}\Big|\Big|^{\sigma_{\vec{i}\vec{j}\vec{l}\vec{v}ij}}_{S_{j}^{\mp}}.
\end{eqnarray*}
In a word, It is not difficult to get
\begin{eqnarray*}
&&\Bigg|\int_{\mathbb{R}}\int_{\mathbb{T}^d}Q^+_{\vec{i}_1}P_{\vec{l}_1}\Pi_+(D)\psi_1\cdots Q^+_{\vec{i}_{d_0}}P_{\vec{l}_{d_0}}\Pi_+(D)\psi_{d_0}\cdot Q^-_{\vec{j}_1}P_{\vec{v}_1}\Pi_-(D)\varphi_1\cdots Q^-_{\vec{j}_{d_0}}P_{\vec{v}_{d_0}}\Pi_-(D)\varphi_{d_0}\\
&&\left\langle C_{\vec{i}\vec{j}\vec{l}\vec{v}},\beta Q^{\mp}_i\Pi_{\pm}(D)P_j\phi\right\rangle_{\mathbb{C}^{d_0}}\Bigg|\\
&\lesssim&2^{2j}\sum\limits_{\substack{|\chi_{r_k}|=\lceil\hat{p}_k^{r_k},\tilde{p}_k^{r_k}\rceil\\|\alpha_{s_k}|=\lceil\hat{q}_k^{s_k}, \tilde{q}_k^{s_k}\rceil\\|\vartheta|=\lceil\hat{q}_{\vec{i}\vec{j}\vec{l}\vec{v}ij},\tilde{q}_{\vec{i}\vec{j}\vec{l}\vec{v}ij}\rceil}}\prod\limits_{k=1}^{d_0}\prod\limits_{s_k=1}^{n_k}\left\{\Big|\Big|P_{v^{s_k}_k}\Pi_-(D)\varphi_k\Big|\Big|^{1-\lambda_{s_k}}_{S_{j}^-}\cdot\Big|\Big|D^{\alpha_{s_k}}P_{v^{s_k}_k}\Pi_-(D)\varphi_k\Big|\Big|^{\lambda_{s_k}}_{S_{j}^-}\right\}\\
&&\cdot\prod\limits_{k=1}^{d_0}\prod\limits_{r_k=1}^{m_k}\left\{\Big|\Big|P_{l_k^{r_k}}\Pi_+(D)\psi_k\Big|\Big|_{S_{j}^+}^{1-\theta_{r_k}}\cdot\Big|\Big|D^{\chi_{r_k}}P_{l_k^{r_k}}\Pi_+(D)\psi_k\Big|\Big|_{S_{j}^+}^{\theta_{r_k}}\right\}\\
&&\Big|\Big|\left\langle C_{\vec{i}\vec{j}\vec{l}\vec{v}},\beta\Pi_{\pm}(D)P_j\phi\right\rangle_{\mathbb{C}^{d_0}}\Big|\Big|^{1-\sigma_{\vec{i}\vec{j}\vec{l}\vec{v}ij}}_{S_{j}^{\mp}}\Big|\Big|\left\langle C_{\vec{i}\vec{j}\vec{l}\vec{v}},\beta D^{\vartheta}\Pi_{\pm}(D)P_j\phi\right\rangle_{\mathbb{C}^{d_0}}\Big|\Big|^{\sigma_{\vec{i}\vec{j}\vec{l}\vec{v}ij}}_{S_{j}^{\mp}}.
\end{eqnarray*}

In the next, we shall estimate
$$\Big|\Big|\left\langle C_{\vec{i}\vec{j}\vec{l}\vec{v}},\beta D^{\theta}\Pi_{\pm}(D)P_j\phi\right\rangle_{\mathbb{C}^{d_0}}\Big|\Big|_{S_{j}^{\mp}}$$
term by term. From Bernstein-Type Lemma, i.e. Lemma 2.1 in chapter 2 of \cite{BCD}, it follows that
\begin{eqnarray*}
\Big|\Big|\left\langle C_{\vec{i}\vec{j}\vec{l}\vec{v}},\beta D^{\theta}\Pi_{\pm}(D)P_j\phi\right\rangle_{\mathbb{C}^{d_0}}\Big|\Big|_{L^{\infty}_tL^2_x}\leq C^{|\theta|+1}2^{|\theta|j}\cdot\Big|\Big|\left\langle C_{\vec{i}\vec{j}\vec{l}\vec{v}},\beta \Pi_{\pm}(D)P_j\phi\right\rangle_{\mathbb{C}^{d_0}}\Big|\Big|_{L^{\infty}_tL^2_x}
\end{eqnarray*}
for an absolute constant $C$, since
$$
\mbox{supp}\left\{\mathcal{F}_{x\mapsto\xi}\left(\left\langle C_{\vec{i}\vec{j}\vec{l}\vec{v}},\beta \Pi_{\pm}(D)P_j\phi\right\rangle_{\mathbb{C}^{d_0}}\right)\right\}\subset\{\xi\in\mathbb{Z}^d|2^j\leq|\xi|\leq2^{j+2}\}.
$$
By the same approach we will get
$$
\Big|\Big|Q_k^{\mp}\left\langle C_{\vec{i}\vec{j}\vec{l}\vec{v}},\beta D^{\theta}\Pi_{\pm}(D)P_j\phi\right\rangle_{\mathbb{C}^{d_0}}\Big|\Big|_{L_t^2L_x^2}\leq C^{|\theta|+1}2^{|\theta|j}\cdot\Big|\Big|Q_k^{\mp}\left\langle C_{\vec{i}\vec{j}\vec{l}\vec{v}},\beta \Pi_{\pm}(D)P_j\phi\right\rangle_{\mathbb{C}^{d_0}}\Big|\Big|_{L^2_tL^2_x}
$$
which implies
$$
\Big|\Big|\left\langle C_{\vec{i}\vec{j}\vec{l}\vec{v}},\beta D^{\theta}\Pi_{\pm}(D)P_j\phi\right\rangle_{\mathbb{C}^{d_0}}\Big|\Big|_{\dot{X}^{\mp,1/2,\infty}}\leq C^{|\theta|+1}2^{|\theta|j}\cdot\Big|\Big|\left\langle C_{\vec{i}\vec{j}\vec{l}\vec{v}},\beta \Pi_{\pm}(D)P_j\phi\right\rangle_{\mathbb{C}^{d_0}}\Big|\Big|_{\dot{X}^{\mp,1/2,\infty}},
$$
and
\begin{eqnarray*}
&&\Big|\Big|\left\langle C_{\vec{i}\vec{j}\vec{l}\vec{v}},\beta D^{\theta}\Pi_{\pm}(D)P_j\phi\right\rangle_{\mathbb{C}^{d_0}}\Big|\Big|_{L^a_tL^b_x[l,k']}\\
&=&\sum\limits_{\kappa\in\mathcal{K}_l}\sum\limits_{n\in\Xi_{k'}}\Big|\Big|\Gamma_{k',n}P_{\kappa}\left\langle C_{\vec{i}\vec{j}\vec{l}\vec{v}},\beta D^{\theta}\Pi_{\pm}(D)P_j\phi\right\rangle_{\mathbb{C}^{d_0}}\Big|\Big|_{L_t^aL_x^b}\\
&=&\sum\limits_{\kappa\in\mathcal{K}_l}\sum\limits_{n\in\Xi_{k'}}\Big|\Big|D^{\theta}\Gamma_{k',n}P_{\kappa}\left\langle C_{\vec{i}\vec{j}\vec{l}\vec{v}},\beta\Pi_{\pm}(D)P_j\phi\right\rangle_{\mathbb{C}^{d_0}}\Big|\Big|_{L_t^aL_x^b}\\
&\leq&\sum\limits_{\kappa\in\mathcal{K}_l}\sum\limits_{n\in\Xi_{k'}}C^{|\theta|+1}2^{|\theta|j}\cdot\Big|\Big|\Gamma_{k',n}P_{\kappa}\left\langle C_{\vec{i}\vec{j}\vec{l}\vec{v}},\beta \Pi_{\pm}(D)P_j\phi\right\rangle_{\mathbb{C}^{d_0}}\Big|\Big|_{L_t^aL_x^b}\\
&=&C^{|\theta|+1}2^{|\theta|j}\cdot\Big|\Big|\left\langle C_{\vec{i}\vec{j}\vec{l}\vec{v}},\beta\Pi_{\pm}(D)P_j\phi\right\rangle_{\mathbb{C}^{d_0}}\Big|\Big|_{L^a_tL^b_x[l,k']}.
\end{eqnarray*}
In conclusion, we have obtained
\begin{eqnarray}\label{14}
\Big|\Big|\left\langle C_{\vec{i}\vec{j}\vec{l}\vec{v}},\beta D^{\theta}\Pi_{\pm}(D)P_j\phi\right\rangle_{\mathbb{C}^{d_0}}\Big|\Big|_{S_{j}^{\mp}}\leq C^{|\theta|+1}2^{|\theta|j}\cdot\Big|\Big|\left\langle C_{\vec{i}\vec{j}\vec{l}\vec{v}},\beta\Pi_{\pm}(D)P_j\phi\right\rangle_{\mathbb{C}^{d_0}}\Big|\Big|_{S_{j}^{\mp}}.
\end{eqnarray}

The same idea leads to
$$
||D^{\alpha}P_{v^{s_k}_k}\Pi_-(D)\varphi_k||_{S^-_{j}}\leq C^{|\alpha|+1}2^{|\alpha|j}||P_{v^{s_k}_k}\Pi_-(D)\varphi_k||_{S^-_{j}}
$$
and
$$
||D^{\chi}P_{l_k^{r_k}}\Pi_+(D)\psi_k||_{S^+_{j}}\leq C^{|\chi|+1}2^{|\chi|j}||P_{l_k^{r_k}}\Pi_+(D)\psi_k||_{S^+_{j}}.
$$

So we get
\begin{eqnarray*}
&&\Bigg|\int_{\mathbb{R}}\int_{\mathbb{T}^d}Q^+_{\vec{i}_1}P_{\vec{l}_1}\Pi_+(D)\psi_1\cdots Q^+_{\vec{i}_{d_0}}P_{\vec{l}_{d_0}}\Pi_+(D)\psi_{d_0}\cdot Q^-_{\vec{j}_1}P_{\vec{v}_1}\Pi_-(D)\varphi_1\cdots Q^-_{\vec{j}_{d_0}}P_{\vec{v}_{d_0}}\Pi_-(D)\varphi_{d_0}\\
&&\left\langle C_{\vec{i}\vec{j}\vec{l}\vec{v}},\beta Q^{\mp}_i\Pi_{\pm}(D)P_j\phi\right\rangle_{\mathbb{C}^{d_0}}\Bigg|\\
&\lesssim&2^{\varrho_{\vec{i}\vec{j}\vec{l}\vec{v}ij}}C^{\nu_{\vec{i}\vec{j}\vec{l}\vec{v}ij}}\prod\limits_{k=1}^{d_0}\prod\limits_{s_k=1}^{n_k}\Big|\Big|P_{v^{s_k}_k}\Pi_-(D)\varphi_k\Big|\Big|_{S_{j}^-}\cdot\prod\limits_{k=1}^{d_0}\prod\limits_{r_k=1}^{m_k}\Big|\Big|P_{l_k^{r_k}}\Pi_+(D)\psi_k\Big|\Big|_{S_{j}^+}\\
&&\cdot\Big|\Big|\left\langle C_{\vec{i}\vec{j}\vec{l}\vec{v}},\beta\Pi_{\pm}(D)P_j\phi\right\rangle_{\mathbb{C}^{d_0}}\Big|\Big|_{S_{j}^{\mp}}
\end{eqnarray*}
where
\begin{eqnarray*}
\varrho_{\vec{i}\vec{j}\vec{l}\vec{v}ij}:&=&\sum\limits_{k=1}^{d_0}\sum\limits_{s_k=1}^{n_k}\lambda_{s_k}\lceil\hat{q}_k^{s_k},\tilde{q}_k^{s_k}\rceil j+\sum\limits_{k=1}^{d_0}\sum\limits_{r_k=1}^{m_k}\theta_{r_k}\lceil\hat{p}_k^{r_k},\tilde{p}_k^{r_k}\rceil j+\sigma_{\vec{i}\vec{j}\vec{l}\vec{v}ij}\lceil\hat{q}_{\vec{i}\vec{j}\vec{l}\vec{v}ij},\tilde{q}_{\vec{i}\vec{j}\vec{l}\vec{v}ij}\rceil j+2j\\
&=&2j+jd|\mathbf{p}|/2-3jd/2
\end{eqnarray*}
and
\begin{eqnarray*}
\nu_{\vec{i}\vec{j}\vec{l}\vec{v}ij}:&=&\sum\limits_{k=1}^{d_0}\sum\limits_{s_k=1}^{n_k}\lambda_{s_k}(\lceil\hat{q}_k^{s_k},\tilde{q}_k^{s_k}\rceil+1)+\sum\limits_{k=1}^{d_0}\sum\limits_{r_k=1}^{m_k}\theta_{r_k}(\lceil\hat{p}_k^{r_k},\tilde{p}_k^{r_k}\rceil+1)+\sigma_{\vec{i}\vec{j}\vec{l}\vec{v}ij}(\lceil\hat{q}_{\vec{i}\vec{j}\vec{l}\vec{v}ij},\tilde{q}_{\vec{i}\vec{j}\vec{l}\vec{v}ij}\rceil+1)\\
&=&d|\mathbf{p}|/2-3d/2+\sum\limits_{k=1}^{d_0}\sum\limits_{s_k=1}^{n_k}\lambda_{s_k}+\sum\limits_{k=1}^{d_0}\sum\limits_{r_k=1}^{m_k}\theta_{r_k}+\sigma_{\vec{i}\vec{j}\vec{l}\vec{v}ij}.
\end{eqnarray*}
Because of
\begin{eqnarray*}
||(\vec{l}_1,\vec{l}_2,\cdots,\vec{l}_{d_0},\vec{v}_1,\vec{v}_2,\cdots,\vec{v}_{d_0},j)||\leq1
\end{eqnarray*}
we have
\begin{eqnarray}\label{15}
\Big|\Big|P_{v^{s_k}_k}\Pi_-(D)\varphi_k\Big|\Big|_{S_{j}^-}\leq\max\left\{2^{(v_k^{s_k}+1-j)/a},2^{(v_k^{s_k}+1-j)/b}\right\}\Big|\Big|P_{v^{s_k}_k}\Pi_-(D)\varphi_k\Big|\Big|_{S_{v^{s_k}_k+1}^-}
\end{eqnarray}
and
\begin{eqnarray*}
\Big|\Big|P_{l_k^{r_k}}\Pi_+(D)\psi_k\Big|\Big|_{S_{j}^+}\leq\max\left\{2^{(l_k^{r_k}+1-j)/a},2^{(l_k^{r_k}+1-j)/b}\right\}\Big|\Big|P_{l_k^{r_k}}\Pi_+(D)\psi_k\Big|\Big|_{S_{l_k^{r_k}+1}^+}
\end{eqnarray*}
which implies
\begin{eqnarray*}
&&\Bigg|\int_{\mathbb{R}}\int_{\mathbb{T}^d}Q^+_{\vec{i}_1}P_{\vec{l}_1}\Pi_+(D)\psi_1\cdots Q^+_{\vec{i}_{d_0}}P_{\vec{l}_{d_0}}\Pi_+(D)\psi_{d_0}\cdot Q^-_{\vec{j}_1}P_{\vec{v}_1}\Pi_-(D)\varphi_1\cdots Q^-_{\vec{j}_{d_0}}P_{\vec{v}_{d_0}}\Pi_-(D)\varphi_{d_0}\\
&&\left\langle C_{\vec{i}\vec{j}\vec{l}\vec{v}},\beta Q^{\mp}_i\Pi_{\pm}(D)P_j\phi\right\rangle_{\mathbb{C}^{d_0}}\Bigg|\\
&\lesssim&2^{2j+jd|\mathbf{p}|/2-3jd/2}C^{\nu_{\vec{i}\vec{j}\vec{l}\vec{v}ij}}\Big|\Big|\left\langle C_{\nu_{\vec{i}\vec{j}\vec{l}\vec{v}ij}},\beta\Pi_{\pm}(D)P_j\phi\right\rangle_{\mathbb{C}^{d_0}}\Big|\Big|_{S_j^{\mp}}\\
&&\prod\limits_{k=1}^{d_0}\prod\limits_{s_k=1}^{n_k}2^{\max\{1/a,1/b\}(v_k^{s_k}+1-j)}\Big|\Big|P_{v^{s_k}_k}\Pi_-(D)\varphi_k\Big|\Big|_{S_{v^{s_k}_k+1}^-}\\
&&\cdot\prod\limits_{k=1}^{d_0}\prod\limits_{r_k=1}^{m_k}2^{\max\{1/a,1/b\}(l_k^{r_k}+1-j)}\Big|\Big|P_{l_k^{r_k}}\Pi_+(D)\psi_k\Big|\Big|_{S_{l_k^{r_k}+1}^+}.
\end{eqnarray*}
From
\begin{eqnarray*}
||f||_{S_j^{\pm}}\lesssim||f||_{S_{j+1}^{\pm}}\s\s\mbox{and}\s\s\max\{1/a,1/b\}\leq1/a+1/b=1/2\leq d/2
\end{eqnarray*}
it follows that
\begin{eqnarray*}
&&\Bigg|\int_{\mathbb{R}}\int_{\mathbb{T}^d}Q^+_{\vec{i}_1}P_{\vec{l}_1}\Pi_+(D)\psi_1\cdots Q^+_{\vec{i}_{d_0}}P_{\vec{l}_{d_0}}\Pi_+(D)\psi_{d_0}\cdot Q^-_{\vec{j}_1}P_{\vec{v}_1}\Pi_-(D)\varphi_1\cdots Q^-_{\vec{j}_{d_0}}P_{\vec{v}_{d_0}}\Pi_-(D)\varphi_{d_0}\\
&&\left\langle C_{\vec{i}\vec{j}\vec{l}\vec{v}},\beta Q^{\mp}_i\Pi_{\pm}(D)P_j\phi\right\rangle_{\mathbb{C}^{d_0}}\Bigg|\\
&\lesssim&2^{q}C^{\nu_{\vec{i}\vec{j}\vec{l}\vec{v}ij}}\Big|\Big|\left\langle C_{\vec{i}\vec{j}\vec{l}\vec{v}},\beta\Pi_{\pm}(D)P_j\phi\right\rangle_{\mathbb{C}^{d_0}}\Big|\Big|_{S_{j+1}^{\mp}}\prod\limits_{k=1}^{d_0}\prod\limits_{s_k=1}^{n_k}\Big|\Big|P_{v^{s_k}_k}\Pi_-(D)\varphi_k\Big|\Big|_{S_{v^{s_k}_k+1}^-}\\
&&\prod\limits_{k=1}^{d_0}\prod\limits_{r_k=1}^{m_k}\Big|\Big|P_{l_k^{r_k}}\Pi_+(D)\psi_k\Big|\Big|_{S_{l_k^{r_k}+1}^+}
\end{eqnarray*}
with
\begin{eqnarray*}
q:=2j-3jd/2+d\left(\sum\limits_{k=1}^{d_0}\sum\limits_{s_k=1}^{n_k}v_k^{s_k}+\sum\limits_{k=1}^{d_0}\sum\limits_{r_k=1}^{m_k}l_k^{r_k}\right)/2+|\mathbf{p}|d/2.
\end{eqnarray*}
So we get
\begin{eqnarray*}
&&\sum\limits_{i\in\mathbb{Z}}\Bigg|\int_{\mathbb{R}}\int_{\mathbb{T}^d}\left\langle F\left(\Pi_+(D)\psi+\Pi_-(D)\varphi\right),\beta Q^{\mp}_i\Pi_{\pm}(D)P_j\phi\right\rangle_{\mathbb{C}^{d_0}}\,dxdt\Bigg|\\
&\lesssim&\sum\limits_i\sum\limits_{|\mathbf{p}|=0}^{\infty}\sum\limits_{\mathbf{m}+\mathbf{n}=\mathbf{p}}\sum\limits_{\substack{\vec{i}:=(\vec{i}_1,\vec{i}_2,\cdots,\vec{i}_{d_0})\\ \vec{j}:=(\vec{j}_1,\vec{j}_2,\cdots,\vec{j}_{d_0})}}\sum\limits_{\substack{\vec{l}:=(\vec{l}_1,\vec{l}_2,\cdots,\vec{l}_{d_0})\\ \vec{v}:=(\vec{v}_1,\vec{v}_2,\cdots,\vec{v}_{d_0})}}2^qC^{\nu_{\vec{i}\vec{j}\vec{l}\vec{v}ij}}\\
&&\Big|\Big|\left\langle C_{\vec{i}\vec{j}\vec{l}\vec{v}},\beta\Pi_{\pm}(D)P_j\phi\right\rangle_{\mathbb{C}^{d_0}}\Big|\Big|_{S_{j+1}^{\mp}}\prod\limits_{k=1}^{d_0}\prod\limits_{s_k=1}^{n_k}\Big|\Big|P_{v^{s_k}_k}\Pi_-(D)\varphi_k\Big|\Big|_{S_{v^{s_k}_k+1}^-}\\
&&\prod\limits_{k=1}^{d_0}\prod\limits_{r_k=1}^{m_k}\Big|\Big|P_{l_k^{r_k}}\Pi_+(D)\psi_k\Big|\Big|_{S_{l_k^{r_k}+1}^+}.
\end{eqnarray*}
In the sequel, we shall estimate
$$
\Big|\Big|\left\langle C_{\vec{i}\vec{j}\vec{l}\vec{v}},\beta\Pi_{\pm}(D)P_j\phi\right\rangle_{\mathbb{C}^{d_0}}\Big|\Big|_{S_{j+1}^{\mp}}.
$$

From the definition it follows that
\begin{eqnarray*}
&&\Big|\Big|\left\langle C_{\vec{i}\vec{j}\vec{l}\vec{v}},\beta\Pi_{\pm}(D)P_j\phi\right\rangle_{\mathbb{C}^{d_0}}\Big|\Big|_{L_t^{\infty}L_x^2}\\
&=&\frac{1}{2}\Big|\Big|\varphi_j(\xi)\left\langle\gamma^0C_{\vec{i}\vec{j}\vec{l}\vec{v}},\mathcal{F}_{x\mapsto\xi}\phi\right\rangle_{\mathbb{C}^{d_0}}\mp\sum\limits_{\dot{k}=1}^d\frac{\xi_{\dot{k}}\varphi_j(\xi)}{\langle\xi\rangle}\left\langle\gamma^{\dot{k}}C_{\vec{i}\vec{j}\vec{l}\vec{v}},\mathcal{F}_{x\mapsto\xi}\phi\right\rangle_{\mathbb{C}^{d_0}}\\
&&\pm\frac{\varphi_j(\xi)}{\langle\xi\rangle}\left\langle C_{\vec{i}\vec{j}\vec{l}\vec{v}},\mathcal{F}_{x\mapsto\xi}\phi\right\rangle_{\mathbb{C}^{d_0}}\Big|\Big|_{L_t^{\infty}l_{\xi}^2}\\
&\leq&\frac{1}{2}\Big|\Big|\varphi_j(\xi)\left\langle\gamma^0C_{\vec{i}\vec{j}\vec{l}\vec{v}},\mathcal{F}_{x\mapsto\xi}\phi\right\rangle_{\mathbb{C}^{d_0}}\Big|\Big|_{L_t^{\infty}l_{\xi}^2}+\frac{1}{2}\sum\limits_{\dot{k}=1}^d\Big|\Big|\frac{\xi_{\dot{k}}\varphi_j(\xi)}{\langle\xi\rangle}\left\langle\gamma^{\dot{k}}C_{\vec{i}\vec{j}\vec{l}\vec{v}},\mathcal{F}_{x\mapsto\xi}\phi\right\rangle_{\mathbb{C}^{d_0}}\Big|\Big|_{L_t^{\infty}l_{\xi}^2}\\
&&+\frac{1}{2}\Big|\Big|\frac{\varphi_j(\xi)}{\langle\xi\rangle}\left\langle C_{\vec{i}\vec{j}\vec{l}\vec{v}},\mathcal{F}_{x\mapsto\xi}\phi\right\rangle_{\mathbb{C}^{d_0}}\Big|\Big|_{L_t^{\infty}l_{\xi}^2}\\
&\leq&\frac{1}{2}\Big|\Big|\varphi_j(\xi)\left\langle\gamma^0C_{\vec{i}\vec{j}\vec{l}\vec{v}},\mathcal{F}_{x\mapsto\xi}\phi\right\rangle_{\mathbb{C}^{d_0}}\Big|\Big|_{L_t^{\infty}l_{\xi}^2}+\frac{1}{2}\sum\limits_{\dot{k}=1}^d\Big|\Big|\xi_{\dot{k}}\varphi_j(\xi)\left\langle\gamma^{\dot{k}}C_{\vec{i}\vec{j}\vec{l}\vec{v}},\mathcal{F}_{x\mapsto\xi}\phi\right\rangle_{\mathbb{C}^{d_0}}\Big|\Big|_{L_t^{\infty}l_{\xi}^2}\\
&&+\frac{1}{2}\Big|\Big|\varphi_j(\xi)\left\langle C_{\vec{i}\vec{j}\vec{l}\vec{v}},\mathcal{F}_{x\mapsto\xi}\phi\right\rangle_{\mathbb{C}^{d_0}}\Big|\Big|_{L_t^{\infty}l_{\xi}^2}\\
&=&\frac{1}{2}\Big|\Big|\left\langle\gamma^0C_{\vec{i}\vec{j}\vec{l}\vec{v}},P_j\phi\right\rangle_{\mathbb{C}^{d_0}}\Big|\Big|_{L_t^{\infty}L_x^2}+\frac{1}{2}\sum\limits_{\dot{k}=1}^d\Big|\Big|D_{\dot{k}}\left\langle\gamma^{\dot{k}}C_{\vec{i}\vec{j}\vec{l}\vec{v}},P_j\phi\right\rangle_{\mathbb{C}^{d_0}}\Big|\Big|_{L_t^{\infty}L_x^2}\\
&&+\frac{1}{2}\Big|\Big|\left\langle C_{\vec{i}\vec{j}\vec{l}\vec{v}},P_j\phi\right\rangle_{\mathbb{C}^{d_0}}\Big|\Big|_{L_t^{\infty}L_x^2}.
\end{eqnarray*}
From Bernstein-Type Lemma, i.e. Lemma 2.1 in chapter 2 of \cite{BCD}, it follows that
$$
\sum\limits_{\dot{k}=1}^d\Big|\Big|D_{\dot{k}}\left\langle\gamma^{\dot{k}}C_{\vec{i}\vec{j}\vec{l}\vec{v}},P_j\phi\right\rangle_{\mathbb{C}^{d_0}}\Big|\Big|_{L_t^{\infty}L_x^2}\lesssim2^{j}\sum\limits_{\dot{k}=1}^d\Big|\Big|\left\langle\gamma^{\dot{k}}C_{\vec{i}\vec{j}\vec{l}\vec{v}},P_j\phi\right\rangle_{\mathbb{C}^{d_0}}\Big|\Big|_{L_t^{\infty}L_x^2}.
$$
Combining the above results yields
\begin{eqnarray*}
&&\Big|\Big|\left\langle C_{\vec{i}\vec{j}\vec{l}\vec{v}},\beta\Pi_{\pm}(D)P_j\phi\right\rangle_{\mathbb{C}^{d_0}}\Big|\Big|_{L_t^{\infty}L_x^2}\\
&\lesssim&2^j\Big|\Big|\left\langle\gamma^0C_{\vec{i}\vec{j}\vec{l}\vec{v}},P_j\phi\right\rangle_{\mathbb{C}^{d_0}}\Big|\Big|_{L_t^{\infty}L_x^2}+2^j\sum\limits_{\dot{k}=1}^d\Big|\Big|\left\langle\gamma^{\dot{k}}C_{\vec{i}\vec{j}\vec{l}\vec{v}},P_j\phi\right\rangle_{\mathbb{C}^{d_0}}\Big|\Big|_{L_t^{\infty}L_x^2}\\
&&+2^j\Big|\Big|\left\langle C_{\vec{i}\vec{j}\vec{l}\vec{v}},P_j\phi\right\rangle_{\mathbb{C}^{d_0}}\Big|\Big|_{L_t^{\infty}L_x^2}.
\end{eqnarray*}
The same idea implies
\begin{eqnarray*}
&&\Big|\Big|\left\langle C_{\vec{i}\vec{j}\vec{l}\vec{v}},\beta\Pi_{\pm}(D)P_j\phi\right\rangle_{\mathbb{C}^{d_0}}\Big|\Big|_{\dot{X}^{\mp,1/2,\infty}}\\
&\lesssim&2^j\Big|\Big|\left\langle\gamma^0C_{\vec{i}\vec{j}\vec{l}\vec{v}},P_j\phi\right\rangle_{\mathbb{C}^{d_0}}\Big|\Big|_{\dot{X}^{\mp,1/2,\infty}}+2^j\sum\limits_{\dot{k}=1}^d\Big|\Big|\left\langle\gamma^{\dot{k}}C_{\vec{i}\vec{j}\vec{l}\vec{v}},P_j\phi\right\rangle_{\mathbb{C}^{d_0}}\Big|\Big|_{\dot{X}^{\mp,1/2,\infty}}\\
&&+2^j\Big|\Big|\left\langle C_{\vec{i}\vec{j}\vec{l}\vec{v}},P_j\phi\right\rangle_{\mathbb{C}^{d_0}}\Big|\Big|_{\dot{X}^{\mp,1/2,\infty}}.
\end{eqnarray*}
As for
$$
\Big|\Big|\left\langle C_{\vec{i}\vec{j}\vec{l}\vec{v}},\beta\Pi_{\pm}(D)P_j\phi\right\rangle_{\mathbb{C}^{d_0}}\Big|\Big|_{L_t^aL_x^b[l,k']},
$$
from the definition it follows that we only need to control
$$
\sum\limits_{\kappa\in\mathcal{K}_l}\sum\limits_{n\in\Xi_{k'}}\Big|\Big|\Gamma_{k',n}P_{\kappa}\left\langle C_{\vec{i}\vec{j}\vec{l}\vec{v}},\beta\Pi_{\pm}(D)P_j\phi\right\rangle_{\mathbb{C}^{d_0}}\Big|\Big|_{L_t^aL_x^b}.
$$
By Gagliardo-Nirenberg-Moser estimates we can get
\begin{eqnarray*}
&&\Big|\Big|\Gamma_{k',n}P_{\kappa}\left\langle C_{\vec{i}\vec{j}\vec{l}\vec{v}},\beta\Pi_{\pm}(D)P_j\phi\right\rangle_{\mathbb{C}^{d_0}}\Big|\Big|_{L_t^aL_x^b}\\
&\lesssim&\sum\limits_{|\alpha|=\lfloor d/2-d/b\rfloor+1}\Big|\Big|\Gamma_{k',n}P_{\kappa}\left\langle C_{\vec{i}\vec{j}\vec{l}\vec{v}},\beta\Pi_{\pm}(D)P_j\phi\right\rangle_{\mathbb{C}^{d_0}}\Big|\Big|^{\theta}_{L_t^aL_x^2}\\
&&\Big|\Big|D^{\alpha}\Gamma_{k',n}P_{\kappa}\left\langle C_{\vec{i}\vec{j}\vec{l}\vec{v}},\beta\Pi_{\pm}(D)P_j\phi\right\rangle_{\mathbb{C}^{d_0}}\Big|\Big|^{1-\theta}_{L_t^aL_x^2}
\end{eqnarray*}
with $\theta:=d/\{2(\lfloor d/2-d/b\rfloor+1)\}-d/\{b(\lfloor d/2-d/b\rfloor+1)\}$, where $\lfloor x\rfloor$ means the integer part of a real number $x$. From Bernstein-Type Lemma, i.e. Lemma 2.1 in chapter 2 of \cite{BCD}, it follows that
\begin{eqnarray*}
&&\Big|\Big|D^{\alpha}\Gamma_{k',n}P_{\kappa}\left\langle C_{\vec{i}\vec{j}\vec{l}\vec{v}},\beta\Pi_{\pm}(D)P_j\phi\right\rangle_{\mathbb{C}^{d_0}}\Big|\Big|_{L_t^aL_x^2}\\
&\lesssim&2^{j(\lfloor d/2-d/b\rfloor+1)}\Big|\Big|\Gamma_{k',n}P_{\kappa}\left\langle C_{\vec{i}\vec{j}\vec{l}\vec{v}},\beta\Pi_{\pm}(D)P_j\phi\right\rangle_{\mathbb{C}^{d_0}}\Big|\Big|_{L_t^aL_x^2}
\end{eqnarray*}
which implies
\begin{eqnarray*}
&&\Big|\Big|\Gamma_{k',n}P_{\kappa}\left\langle C_{\vec{i}\vec{j}\vec{l}\vec{v}},\beta\Pi_{\pm}(D)P_j\phi\right\rangle_{\mathbb{C}^{d_0}}\Big|\Big|_{L_t^aL_x^b}\\
&\lesssim&2^{j(\lfloor d/2-d/b\rfloor+1)(1-\theta)}\Big|\Big|\Gamma_{k',n}P_{\kappa}\left\langle C_{\vec{i}\vec{j}\vec{l}\vec{v}},\beta\Pi_{\pm}(D)P_j\phi\right\rangle_{\mathbb{C}^{d_0}}\Big|\Big|_{L_t^aL_x^2}\\
&\lesssim&2^{j(\lfloor d/2-d/b\rfloor+1)(1-\theta)}\Big|\Big|\Gamma_{k',n}P_{\kappa}\left\langle\gamma^0C_{\vec{i}\vec{j}\vec{l}\vec{v}},P_j\phi\right\rangle_{\mathbb{C}^{d_0}}\Big|\Big|_{L_t^aL_x^2}\\
&&+2^{j(\lfloor d/2-d/b\rfloor+1)(1-\theta)}\sum\limits_{\dot{k}=1}^d\Big|\Big|\Gamma_{k',n}P_{\kappa}\left\langle\gamma^{\dot{k}}C_{\vec{i}\vec{j}\vec{l}\vec{v}},P_j\phi\right\rangle_{\mathbb{C}^{d_0}}\Big|\Big|_{L_t^aL_x^2}\\
&&+2^{j(\lfloor d/2-d/b\rfloor+1)(1-\theta)}\Big|\Big|\Gamma_{k',n}P_{\kappa}\left\langle C_{\vec{i}\vec{j}\vec{l}\vec{v}},P_j\phi\right\rangle_{\mathbb{C}^{d_0}}\Big|\Big|_{L_t^aL_x^2}\\
&\lesssim&2^{j(\lfloor d/2-d/b\rfloor+1)(1-\theta)}\Big|\Big|\Gamma_{k',n}P_{\kappa}\left\langle\gamma^0C_{\vec{i}\vec{j}\vec{l}\vec{v}},P_j\phi\right\rangle_{\mathbb{C}^{d_0}}\Big|\Big|_{L_t^aL_x^b}\\
&&+2^{j(\lfloor d/2-d/b\rfloor+1)(1-\theta)}\sum\limits_{\dot{k}=1}^d\Big|\Big|\Gamma_{k',n}P_{\kappa}\left\langle\gamma^{\dot{k}}C_{\vec{i}\vec{j}\vec{l}\vec{v}},P_j\phi\right\rangle_{\mathbb{C}^{d_0}}\Big|\Big|_{L_t^aL_x^b}\\
&&+2^{j(\lfloor d/2-d/b\rfloor+1)(1-\theta)}\Big|\Big|\Gamma_{k',n}P_{\kappa}\left\langle C_{\vec{i}\vec{j}\vec{l}\vec{v}},P_j\phi\right\rangle_{\mathbb{C}^{d_0}}\Big|\Big|_{L_t^aL_x^b},
\end{eqnarray*}
where we have used the fact that flat torus has finite measure. In summary, we have obtained
\begin{eqnarray*}
&&\Big|\Big|\left\langle C_{\vec{i}\vec{j}\vec{l}\vec{v}},\beta\Pi_{\pm}(D)P_j\phi\right\rangle_{\mathbb{C}^{d_0}}\Big|\Big|_{L_t^aL_x^b[l,k']}\\
&\lesssim&2^{j(\lfloor d/2-d/b\rfloor+1)(1-\theta)}\Big|\Big|\left\langle\gamma^0C_{\vec{i}\vec{j}\vec{l}\vec{v}},P_j\phi\right\rangle_{\mathbb{C}^{d_0}}\Big|\Big|_{L_t^aL_x^b[l,k']}\\
&&+2^{j(\lfloor d/2-d/b\rfloor+1)(1-\theta)}\sum\limits_{\dot{k}=1}^d\Big|\Big|\left\langle\gamma^{\dot{k}}C_{\vec{i}\vec{j}\vec{l}\vec{v}},P_j\phi\right\rangle_{\mathbb{C}^{d_0}}\Big|\Big|_{L_t^aL_x^b[l,k']}\\
&&+2^{j(\lfloor d/2-d/b\rfloor+1)(1-\theta)}\Big|\Big|\left\langle C_{\vec{i}\vec{j}\vec{l}\vec{v}},P_j\phi\right\rangle_{\mathbb{C}^{d_0}}\Big|\Big|_{L_t^aL_x^b[l,k']}.
\end{eqnarray*}

So it is not difficult to get
\begin{eqnarray}\label{16}
&&\Big|\Big|\left\langle C_{\vec{i}\vec{j}\vec{l}\vec{v}},\beta \Pi_{\pm}(D)P_j\phi\right\rangle_{\mathbb{C}^{d_0}}\Big|\Big|_{S_{j+1}^{\mp}}\\
&\lesssim&2^{j(\max\{\lfloor d/2-d/b\rfloor,\lfloor d/2-d/a\rfloor\}+1)}\Big|\Big|\left\langle\gamma^0C_{\vec{i}\vec{j}\vec{l}\vec{v}},P_j\phi\right\rangle_{\mathbb{C}^{d_0}}\Big|\Big|_{S_{j+1}^{\mp}}\nonumber\\
&&+2^{j(\max\{\lfloor d/2-d/b\rfloor,\lfloor d/2-d/a\rfloor\}+1)}\sum\limits_{\dot{k}=1}^d\Big|\Big|\left\langle\gamma^{\dot{k}}C_{\vec{i}\vec{j}\vec{l}\vec{v}},P_j\phi\right\rangle_{\mathbb{C}^{d_0}}\Big|\Big|_{S_{j+1}^{\mp}}\nonumber\\
&&+2^{j(\max\{\lfloor d/2-d/b\rfloor,\lfloor d/2-d/a\rfloor\}+1)}\Big|\Big|\left\langle C_{\vec{i}\vec{j}\vec{l}\vec{v}},P_j\phi\right\rangle_{\mathbb{C}^{d_0}}\Big|\Big|_{S_{j+1}^{\mp}},\nonumber
\end{eqnarray}
which implies
\begin{eqnarray}\label{4}
&&\sum\limits_{i\in\mathbb{Z}}\Bigg|\int_{\mathbb{R}}\int_{\mathbb{T}^d}\left\langle F\left(\Pi_+(D)\psi+\Pi_-(D)\varphi\right),\beta Q^{\mp}_i\Pi_{\pm}(D)P_j\phi\right\rangle_{\mathbb{C}^{d_0}}\,dxdt\Bigg|\nonumber\\
&\lesssim&\sum\limits_i\sum\limits_{|\mathbf{p}|=1}^{\infty}\sum\limits_{\mathbf{m}+\mathbf{n}=\mathbf{p}}\sum\limits_{\substack{\vec{i}:=(\vec{i}_1,\vec{i}_2,\cdots,\vec{i}_{d_0})\\ \vec{j}:=(\vec{j}_1,\vec{j}_2,\cdots,\vec{j}_{d_0})}}\sum\limits_{\substack{\vec{l}:=(\vec{l}_1,\vec{l}_2,\cdots,\vec{l}_{d_0})\\ \vec{v}:=(\vec{v}_1,\vec{v}_2,\cdots,\vec{v}_{d_0})}}2^{q+j(\max\{\lfloor d/2-d/b\rfloor,\lfloor d/2-d/a\rfloor\}+1)}\\
&&C^{\nu_{\vec{i}\vec{j}\vec{l}\vec{v}ij}}\prod\limits_{k=1}^{d_0}\prod\limits_{s_k=1}^{n_k}\Big|\Big|P_{v^{s_k}_k}\Pi_-(D)\varphi_k\Big|\Big|_{S_{v^{s_k}_k+1}^-}\prod\limits_{k=1}^{d_0}\prod\limits_{r_k=1}^{m_k}\Big|\Big|P_{l_k^{r_k}}\Pi_+(D)\psi_k\Big|\Big|_{S_{l_k^{r_k}+1}^+}\nonumber\\
&&\Bigg\{\Big|\Big|\left\langle\gamma^0C_{\vec{i}\vec{j}\vec{l}\vec{v}},P_j\phi\right\rangle_{\mathbb{C}^{d_0}}\Big|\Big|_{S_{j+1}^{\mp}}+\sum\limits_{\dot{k}=1}^d\Big|\Big|\left\langle\gamma^{\dot{k}}C_{\vec{i}\vec{j}\vec{l}\vec{v}},P_j\phi\right\rangle_{\mathbb{C}^{d_0}}\Big|\Big|_{S_{j+1}^{\mp}}+\Big|\Big|\left\langle C_{\vec{i}\vec{j}\vec{l}\vec{v}},P_j\phi\right\rangle_{\mathbb{C}^{d_0}}\Big|\Big|_{S_{j+1}^{\mp}}\Bigg\}\nonumber.
\end{eqnarray}
From direct computation it follows that
\begin{eqnarray*}
&&\sum\limits_{k=1}^{d_0}\sum\limits_{s_k=1}^{n_k}\lambda_{s_k}+\sum\limits_{k=1}^{d_0}\sum\limits_{r_k=1}^{m_k}\theta_{r_k}+\sigma_{\vec{i}\vec{j}\vec{l}\vec{v}ij}\\
&=&\sum\limits_{k=1}^{d_0}\sum\limits_{s_k=1}^{n_k}\frac{1/2-1/q^{s_k}_k-1/\tilde{q}^{s_k}_k}{2+\lceil1/2-1/\tilde{q}^{s_k}_k-1/q^{s_k}_k\rceil}+\sum\limits_{k=1}^{d_0}\sum\limits_{r_k=1}^{m_k}\frac{1/2-1/p^{r_k}_k-1/\tilde{p}^{r_k}_k}{2+\lceil1/2-1/\tilde{p}^{r_k}_k-1/p^{r_k}_k\rceil}\\
&&+\frac{1/2-1/q_{\vec{i}\vec{j}\vec{l}\vec{v}ij}-1/\tilde{q}_{\vec{i}\vec{j}\vec{l}\vec{v}ij}}{2+\lceil 1/2-1/\tilde{q}_{\vec{i}\vec{j}\vec{l}\vec{v}ij}-1/q_{\vec{i}\vec{j}\vec{l}\vec{v}ij}\rceil}
\end{eqnarray*}
with $\lceil x\rceil:=\lfloor x\rfloor+1$. Choose $p_k^{r_k}$, $\tilde{p}_k^{r_k}$, $q_k^{s_k}$, $\tilde{q}_k^{s_k}$, $q_{\vec{i}\vec{j}\vec{l}\vec{v}ij}$ and $\tilde{q}_{\vec{i}\vec{j}\vec{l}\vec{v}ij}$ such that
$$
1/p_k^{r_k}=1/\tilde{p}_k^{r_k}:=2^{-|i_k^{r_k}|-1}/(|\mathbf{p}|+1)
$$
and
$$
1/q_k^{s_k}=1/\tilde{q}_k^{s_k}:=2^{-|j_k^{s_k}|-1}/(|\mathbf{p}|+1)
$$
and
$$ 1/q_{\vec{i}\vec{j}\vec{l}\vec{v}ij}=1/\tilde{q}_{\vec{i}\vec{j}\vec{l}\vec{v}ij}:=1-\sum\limits_{k=1}^{d_0}\sum\limits_{s_k=1}^{n_k}2^{-|j_k^{s_k}|-1}/(|\mathbf{p}|+1)-\sum\limits_{k=1}^{d_0}\sum\limits_{r_k=1}^{m_k}2^{-|i_k^{r_k}|-1}/(|\mathbf{p}|+1),
$$
where $f$ and $g$ are to be determined. Then we have
\begin{eqnarray*}
&&\sum\limits_{k=1}^{d_0}\sum\limits_{s_k=1}^{n_k}\lambda_{s_k}+\sum\limits_{k=1}^{d_0}\sum\limits_{r_k=1}^{m_k}\theta_{r_k}+\sigma_{\vec{i}\vec{j}\vec{l}\vec{v}ij}\\
&=&\sum\limits_{k=1}^{d_0}\sum\limits_{s_k=1}^{n_k}\frac{1/2-2^{-|j_k^{s_k}|}/(|\mathbf{p}|+1)}{2+\lceil1/2-2^{-|j_k^{s_k}|}/(|\mathbf{p}|+1)\rceil}+\sum\limits_{k=1}^{d_0}\sum\limits_{r_k=1}^{m_k}\frac{1/2-2^{-|i_k^{r_k}|}/(|\mathbf{p}|+1)}{2+\lceil1/2-2^{-|i_k^{r_k}|}/(|\mathbf{p}|+1)\rceil}\\
&&+\frac{\sum\limits_{k=1}^{d_0}\sum\limits_{s_k=1}^{n_k}2^{-|j_k^{s_k}|}/(|\mathbf{p}|+1)+\sum\limits_{k=1}^{d_0}\sum\limits_{r_k=1}^{m_k}2^{-|i_k^{r_k}|}/(|\mathbf{p}|+1)-3/2}{2+\left\lceil\sum\limits_{k=1}^{d_0}\sum\limits_{s_k=1}^{n_k}2^{-|j_k^{s_k}|}/(|\mathbf{p}|+1)+\sum\limits_{k=1}^{d_0}\sum\limits_{r_k=1}^{m_k}2^{-|i_k^{r_k}|}/(|\mathbf{p}|+1)-3/2\right\rceil}\\
&\leq&-\sum\limits_{k=1}^{d_0}\sum\limits_{s_k=1}^{n_k}2^{-|j_k^{s_k}|-1}/(|\mathbf{p}|+1)-\sum\limits_{k=1}^{d_0}\sum\limits_{r_k=1}^{m_k}2^{-|i_k^{r_k}|-1}/(|\mathbf{p}|+1)+1+|\mathbf{p}|/4
\end{eqnarray*}
and
\begin{eqnarray*}
\nu_{\vec{i}\vec{j}\vec{l}\vec{v}ij}&\leq& -\sum\limits_{k=1}^{d_0}\sum\limits_{s_k=1}^{n_k}2^{-|j_k^{s_k}|-1}/(|\mathbf{p}|+1)-\sum\limits_{k=1}^{d_0}\sum\limits_{r_k=1}^{m_k}2^{-|i_k^{r_k}|-1}/(|\mathbf{p}|+1)\\
&&+d|\mathbf{p}|/2-3d/2+1+|\mathbf{p}|/4.
\end{eqnarray*}
Substituting the above inequality into (\ref{4}) leads to
\begin{eqnarray*}
&&\sum\limits_{i\in\mathbb{Z}}\Bigg|\int_{\mathbb{R}}\int_{\mathbb{T}^d}\left\langle F\left(\Pi_+(D)\psi+\Pi_-(D)\varphi\right),\beta Q^{\mp}_i\Pi_{\pm}(D)P_j\phi\right\rangle_{\mathbb{C}^{d_0}}\,dxdt\Bigg|\nonumber\\
&\lesssim&\sum\limits_i\sum\limits_{|\mathbf{p}|=1}^{\infty}\sum\limits_{\mathbf{m}+\mathbf{n}=\mathbf{p}}\sum\limits_{\substack{\vec{i}:=(\vec{i}_1,\vec{i}_2,\cdots,\vec{i}_{d_0})\\ \vec{j}:=(\vec{j}_1,\vec{j}_2,\cdots,\vec{j}_{d_0})}}\sum\limits_{\substack{\vec{l}:=(\vec{l}_1,\vec{l}_2,\cdots,\vec{l}_{d_0})\\ \vec{v}:=(\vec{v}_1,\vec{v}_2,\cdots,\vec{v}_{d_0})}}2^{q+j(\max\{\lfloor d/2-d/b\rfloor,\lfloor d/2-d/a\rfloor\}+1)}\\
&&C^{-\sum\limits_{k=1}^{d_0}\sum\limits_{s_k=1}^{n_k}2^{-|j_k^{s_k}|-1}/(|\mathbf{p}|+1)-\sum\limits_{k=1}^{d_0}\sum\limits_{r_k=1}^{m_k}2^{-|i_k^{r_k}|-1}/(|\mathbf{p}|+1)+d|\mathbf{p}|/2+|\mathbf{p}|/4}\\
&&\prod\limits_{k=1}^{d_0}\prod\limits_{s_k=1}^{n_k}\Big|\Big|P_{v^{s_k}_k}\Pi_-(D)\varphi_k\Big|\Big|_{S_{v^{s_k}_k+1}^-}\prod\limits_{k=1}^{d_0}\prod\limits_{r_k=1}^{m_k}\Big|\Big|P_{l_k^{r_k}}\Pi_+(D)\psi_k\Big|\Big|_{S_{l_k^{r_k}+1}^+}\\
&&\Bigg\{\Big|\Big|\left\langle\gamma^0C_{\vec{i}\vec{j}\vec{l}\vec{v}},P_j\phi\right\rangle_{\mathbb{C}^{d_0}}\Big|\Big|_{S_{j+1}^{\mp}}+\sum\limits_{\dot{k}=1}^d\Big|\Big|\left\langle\gamma^{\dot{k}}C_{\vec{i}\vec{j}\vec{l}\vec{v}},P_j\phi\right\rangle_{\mathbb{C}^{d_0}}\Big|\Big|_{S_{j+1}^{\mp}}+\Big|\Big|\left\langle C_{\vec{i}\vec{j}\vec{l}\vec{v}},P_j\phi\right\rangle_{\mathbb{C}^{d_0}}\Big|\Big|_{S_{j+1}^{\mp}}\Bigg\}\\
&\lesssim&\sum\limits_i\sum\limits_{|\mathbf{p}|=1}^{\infty}\sum\limits_{\mathbf{m}+\mathbf{n}=\mathbf{p}}\sum\limits_{\substack{\vec{i}:=(\vec{i}_1,\vec{i}_2,\cdots,\vec{i}_{d_0})\\ \vec{j}:=(\vec{j}_1,\vec{j}_2,\cdots,\vec{j}_{d_0})}}\sum\limits_{\substack{\vec{l}:=(\vec{l}_1,\vec{l}_2,\cdots,\vec{l}_{d_0})\\ \vec{v}:=(\vec{v}_1,\vec{v}_2,\cdots,\vec{v}_{d_0})}}2^{q+j(\max\{\lfloor d/2-d/b\rfloor,\lfloor d/2-d/a\rfloor\}+1)}\\
&&C^{-\sum\limits_{k=1}^{d_0}\sum\limits_{s_k=1}^{n_k}2^{-|j_k^{s_k}|-1}/(|\mathbf{p}|+1)-\sum\limits_{k=1}^{d_0}\sum\limits_{r_k=1}^{m_k}2^{-|i_k^{r_k}|-1}/(|\mathbf{p}|+1)+d|\mathbf{p}|/2+|\mathbf{p}|/4}\\
&&\prod\limits_{k=1}^{d_0}\prod\limits_{s_k=1}^{n_k}\Big|\Big|P_{v^{s_k}_k}\Pi_-(D)\varphi_k\Big|\Big|_{S_{v^{s_k}_k+1}^-}\prod\limits_{k=1}^{d_0}\prod\limits_{r_k=1}^{m_k}\Big|\Big|P_{l_k^{r_k}}\Pi_+(D)\psi_k\Big|\Big|_{S_{l_k^{r_k}+1}^+}\\
&&\Bigg\{|\gamma^0C_{\vec{i}\vec{j}\vec{l}\vec{v}}|\cdot\Big|\Big|P_j\phi\Big|\Big|_{S_{j+1}^{\mp}}+\sum\limits_{\dot{k}=1}^d|\gamma^{\dot{k}}C_{\vec{i}\vec{j}\vec{l}\vec{v}}|\cdot\Big|\Big|P_j\phi\Big|\Big|_{S_{j+1}^{\mp}}+|C_{\vec{i}\vec{j}\vec{l}\vec{v}}|\cdot\Big|\Big|P_j\phi\Big|\Big|_{S_{j+1}^{\mp}}\Bigg\}.
\end{eqnarray*}
Recalling
$$
||(\vec{i}_1,\vec{i}_2,\cdots,\vec{i}_{d_0},\vec{j}_1,\vec{j}_2\cdots,\vec{j}_{d_0},i)||\leq1
$$
we have
\begin{eqnarray*}
&&\sum\limits_{i\in\mathbb{Z}}\Bigg|\int_{\mathbb{R}}\int_{\mathbb{T}^d}\left\langle F\left(\Pi_+(D)\psi+\Pi_-(D)\varphi\right),\beta Q^{\mp}_i\Pi_{\pm}(D)P_j\phi\right\rangle_{\mathbb{C}^{d_0}}\,dxdt\Bigg|\nonumber\\
&\lesssim&\sum\limits_{i\in\mathbb{Z}}\sum\limits_{|\mathbf{p}|=1}^{\infty}\sum\limits_{\mathbf{m}+\mathbf{n}=\mathbf{p}}\sum\sum\limits_{\substack{\vec{l}:=(\vec{l}_1,\vec{l}_2,\cdots,\vec{l}_{d_0})\\ \vec{v}:=(\vec{v}_1,\vec{v}_2,\cdots,\vec{v}_{d_0})}}2^{q+j(\max\{\lfloor d/2-d/b\rfloor,\lfloor d/2-d/a\rfloor\}+1)}3^{|\mathbf{p}|}\\
&&C^{d|\mathbf{p}|/2+|\mathbf{p}|/4-2^{-|i|}|\mathbf{p}|/(|\mathbf{p}|+1)}\prod\limits_{k=1}^{d_0}\prod\limits_{s_k=1}^{n_k}\Big|\Big|P_{v^{s_k}_k}\Pi_-(D)\varphi_k\Big|\Big|_{S_{v^{s_k}_k+1}^-}\prod\limits_{k=1}^{d_0}\prod\limits_{r_k=1}^{m_k}\Big|\Big|P_{l_k^{r_k}}\Pi_+(D)\psi_k\Big|\Big|_{S_{l_k^{r_k}+1}^+}\\
&&\Bigg\{|\gamma^0C_{\vec{i}\vec{j}\vec{l}\vec{v}}|\cdot\Big|\Big|P_j\phi\Big|\Big|_{S_{j+1}^{\mp}}+\sum\limits_{\dot{k}=1}^d|\gamma^{\dot{k}} C_{\vec{i}\vec{j}\vec{l}\vec{v}}|\cdot\Big|\Big|P_j\phi\Big|\Big|_{S_{j+1}^{\mp}}+|C_{\vec{i}\vec{j}\vec{l}\vec{v}}|\cdot\Big|\Big|P_j\phi\Big|\Big|_{S_{j+1}^{\mp}}\Bigg\},\nonumber
\end{eqnarray*}
where
\begin{eqnarray*}
\sum:=\sum\limits_{\vec{i}_1=i-1}^{i+1}\sum\limits_{\vec{i}_2=i-1}^{i+1}\cdots\sum\limits_{\vec{i}_{d_0}=i-1}^{i+1}\sum\limits_{\vec{j}_1=i-1}^{i+1}\sum\limits_{\vec{j}_2=i-1}^{i+1}\cdots\sum\limits_{\vec{j}_{d_0}=i-1}^{i+1}
\end{eqnarray*}
and
\begin{eqnarray*}
\sum\limits_{\vec{i}_k=i-1}^{i+1}:=\sum\limits_{i^1_k=i-1}^{i+1}\sum\limits_{i^2_k=i-1}^{i+1}\cdots\sum\limits_{i^{m_k}_k=i-1}^{i+1}\s\mbox{and}\s\sum\limits_{\vec{j}_k=i-1}^{i+1}:=\sum\limits_{j^1_k=i-1}^{i+1}\sum\limits_{j^2_k=i-1}^{i+1}\cdots\sum\limits_{j^{n_k}_k=i-1}^{i+1}.
\end{eqnarray*}
Define
$$
B^1_{\vec{l}\vec{v}}:=\sup\limits_{i\in\mathbb{Z}}\sum|C_{\vec{i}\vec{j}\vec{l}\vec{v}}|,\s B^2_{\vec{l}\vec{v}}:=\sup\limits_{i\in\mathbb{Z}}\sum|\gamma^0C_{\vec{i}\vec{j}\vec{l}\vec{v}}|\s\mbox{and}\s B^3_{\vec{l}\vec{v}}:=\sup\limits_{i\in\mathbb{Z}}\sum\sum\limits_{\dot{k}=1}^d|\gamma^{\dot{k}} C_{\vec{i}\vec{j}\vec{l}\vec{v}}|.
$$
Therefore, we can obtain
\begin{eqnarray*}
&&\sum\limits_{i\in\mathbb{Z}}\Bigg|\int_{\mathbb{R}}\int_{\mathbb{T}^d}\left\langle F\left(\Pi_+(D)\psi+\Pi_-(D)\varphi\right),\beta Q^{\mp}_i\Pi_{\pm}(D)P_j\phi\right\rangle_{\mathbb{C}^{d_0}}\,dxdt\Bigg|\nonumber\\
&\lesssim&\sum\limits_{|\mathbf{p}|=1}^{\infty}C^{d|\mathbf{p}|/2+|\mathbf{p}|/4}3^{|\mathbf{p}|}\sum\limits_{\mathbf{m}+\mathbf{n}=\mathbf{p}}\sum\limits_{i\in\mathbb{Z}}C^{-2^{-|i|}|\mathbf{p}|/(|\mathbf{p}|+1)}\sum\limits_{\substack{\vec{l}:=(\vec{l}_1,\vec{l}_2,\cdots,\vec{l}_{d_0})\\ \vec{v}:=(\vec{v}_1,\vec{v}_2,\cdots,\vec{v}_{d_0})}}\\
&&2^{q+j(\max\{\lfloor d/2-d/b\rfloor,\lfloor d/2-d/a\rfloor\}+1)}\prod\limits_{k=1}^{d_0}\prod\limits_{s_k=1}^{n_k}\Big|\Big|P_{v^{s_k}_k}\Pi_-(D)\varphi_k\Big|\Big|_{S_{v^{s_k}_k+1}^-}\\
&&\Bigg\{B^2_{\vec{l}\vec{v}}\cdot\Big|\Big|P_j\phi\Big|\Big|_{S_{j+1}^{\mp}}+B^3_{\vec{l}\vec{v}}\cdot\Big|\Big|P_j\phi\Big|\Big|_{S_{j+1}^{\mp}}+B^1_{\vec{l}\vec{v}}\cdot\Big|\Big|P_j\phi\Big|\Big|_{S_{j+1}^{\mp}}\Bigg\}\\
&&\prod\limits_{k=1}^{d_0}\prod\limits_{r_k=1}^{m_k}\Big|\Big|P_{l_k^{r_k}}\Pi_+(D)\psi_k\Big|\Big|_{S_{l_k^{r_k}+1}^+}.
\end{eqnarray*}
If letting
\begin{eqnarray}\label{17}
r(|\mathbf{p}|):=\sum\limits_{i\in\mathbb{Z}}C^{-2^{-|i|}|\mathbf{p}|/(|\mathbf{p}|+1)}
\end{eqnarray}
we shall get $r(|\mathbf{p}|)<\infty$. Indeed, the following inequality holds true
$$
\sum\limits_{i\in\mathbb{Z}}C^{-2^{-|i|}|\mathbf{p}|/(|\mathbf{p}|+1)}\leq\sum\limits_{j\in\mathbb{Z}}C^{-|\mathbf{p}|/\{(|j|+1)(|\mathbf{p}|+1)\}}.
$$
Since for large positive real number $x$, it is not difficult to get
$$
\exp(-|\mathbf{p}|\log C/\{(x+1)(|\mathbf{p}|+1)\})\leq |\mathbf{p}|^2(\log C)^2/\{(x+1)^2(|\mathbf{p}|+1)^2\},
$$
one are able to have
$$
\int_0^{\infty}C^{-|\mathbf{p}|/\{(x+1)(|\mathbf{p}|+1)\}}\,dx<\infty
$$
which implies the convergence of $r(|\mathbf{p}|)$. The above result leads to
\begin{eqnarray*}
&&\sum\limits_{i\in\mathbb{Z}}\Bigg|\int_{\mathbb{R}}\int_{\mathbb{T}^d}\left\langle F\left(\Pi_+(D)\psi+\Pi_-(D)\varphi\right),\beta Q^{\mp}_i\Pi_{\pm}(D)P_j\phi\right\rangle_{\mathbb{C}^{d_0}}\,dxdt\Bigg|\nonumber\\
&\lesssim&2^{j(\max\{\lfloor d/2-d/b\rfloor,\lfloor d/2-d/a\rfloor\}+3)-3jd/2}\sum\limits_{|\mathbf{p}|=1}^{\infty}C^{d|\mathbf{p}|/2+|\mathbf{p}|/4}r(|\mathbf{p}|)2^{d|\mathbf{p}|/2}3^{|\mathbf{p}|}\\
&&\sum\limits_{\mathbf{m}+\mathbf{n}=\mathbf{p}}\sum\limits_{\vec{v}:=(\vec{v}_1,\vec{v}_2,\cdots,\vec{v}_{d_0})}\prod\limits_{k=1}^{d_0}\prod\limits_{s_k=1}^{n_k}\left\{2^{dv_k^{s_k}/2}\Big|\Big|P_{v^{s_k}_k}\Pi_-(D)\varphi_k\Big|\Big|_{S_{v^{s_k}_k+1}^-}\right\}\\
&&\sum\limits_{\vec{l}:=(\vec{l}_1,\vec{l}_2,\cdots,\vec{l}_{d_0})}\prod\limits_{k=1}^{d_0}\prod\limits_{r_k=1}^{m_k}\left\{2^{dl_k^{r_k}/2}\Big|\Big|P_{l_k^{r_k}}\Pi_+(D)\psi_k\Big|\Big|_{S_{l_k^{r_k}+1}^+}\right\}\\
&&\Bigg\{B^2_{\vec{l}\vec{v}}\cdot\Big|\Big|P_j\phi\Big|\Big|_{S_{j+1}^{\mp}}+B^3_{\vec{l}\vec{v}}\cdot\Big|\Big|P_j\phi\Big|\Big|_{S_{j+1}^{\mp}}+B^1_{\vec{l}\vec{v}}\cdot\Big|\Big|P_j\phi\Big|\Big|_{S_{j+1}^{\mp}}\Bigg\}.
\end{eqnarray*}
For $q=1,2,3$, we set
$$
B^q_{\mathbf{mn}}:=\sup\limits_{\vec{l}_1}\sup\limits_{\vec{l}_2}\cdots\sup\limits_{\vec{l}_{d_0}}\sup\limits_{\vec{v}_1}\sup\limits_{\vec{v}_2}\cdots\sup\limits_{\vec{v}_{d_0}}B^q_{\vec{l}\vec{v}}
$$
with
$$
\sup\limits_{\vec{l}_k}:=\sup\limits_{l^1_k}\cdots\sup\limits_{l^{m_k}_k}\s\mbox{and}\s\sup\limits_{\vec{v}_k}:=\sup\limits_{v^1_k}\cdots\sup\limits_{v^{n_k}_k}.
$$
Moreover, we also define
$$
B_{r}:=\sup\limits_{q=1,2,3}\sup\limits_{|\mathbf{m}|+|\mathbf{n}|=r}B^q_{\mathbf{mn}}.
$$
Adopting the above notation one can get
\begin{eqnarray*}
&&\sum\limits_{i\in\mathbb{Z}}\Bigg|\int_{\mathbb{R}}\int_{\mathbb{T}^d}\left\langle F\left(\Pi_+(D)\psi+\Pi_-(D)\varphi\right),\beta Q^{\mp}_i\Pi_{\pm}(D)P_j\phi\right\rangle_{\mathbb{C}^{d_0}}\,dxdt\Bigg|\nonumber\\
&\lesssim&2^{j(\max\{\lfloor d/2-d/b\rfloor,\lfloor d/2-d/a\rfloor\}+3)-3jd/2}\Big|\Big|P_j\phi\Big|\Big|_{S_{j+1}^{\mp}}\sum\limits_{|\mathbf{p}|=1}^{\infty}B_{|\mathbf{p}|}C^{d|\mathbf{p}|/2+|\mathbf{p}|/4}r(|\mathbf{p}|)2^{d|\mathbf{p}|/2}3^{|\mathbf{p}|}\\
&&\sum\limits_{\mathbf{m}+\mathbf{n}=\mathbf{p}}\sum\limits_{\vec{v}:=(\vec{v}_1,\vec{v}_2,\cdots,\vec{v}_{d_0})}\prod\limits_{k=1}^{d_0}\prod\limits_{s_k=1}^{n_k}\left\{2^{dv_k^{s_k}/2}\Big|\Big|P_{v^{s_k}_k}\Pi_-(D)\varphi_k\Big|\Big|_{S_{v^{s_k}_k+1}^-}\right\}\\
&&\sum\limits_{\vec{l}:=(\vec{l}_1,\vec{l}_2,\cdots,\vec{l}_{d_0})}\prod\limits_{k=1}^{d_0}\prod\limits_{r_k=1}^{m_k}\left\{2^{dl_k^{r_k}/2}\Big|\Big|P_{l_k^{r_k}}\Pi_+(D)\psi_k\Big|\Big|_{S_{l_k^{r_k}+1}^+}\right\}\\
&\leq&2^{j(\max\{\lfloor d/2-d/b\rfloor,\lfloor d/2-d/a\rfloor\}+3)-3jd/2}\Big|\Big|P_j\phi\Big|\Big|_{S_{j+1}^{\mp}}\sum\limits_{|\mathbf{p}|=1}^{\infty}B_{|\mathbf{p}|}C^{d|\mathbf{p}|/2+|\mathbf{p}|/4}r(|\mathbf{p}|)2^{d|\mathbf{p}|/2}3^{|\mathbf{p}|}\\
&&\sum\limits_{|\mathbf{m}|+|\mathbf{n}|=|\mathbf{p}|}\Big|\Big|\Pi_-(D)\varphi\Big|\Big|_{S^{-,d/2}}^{|\mathbf{n}|}\Big|\Big|\Pi_+(D)\psi\Big|\Big|_{S^{+,d/2}}^{|\mathbf{m}|}.
\end{eqnarray*}
Define
$$
h(t,\tau):=\sum\limits_{p=1}^{\infty}B_pC^{dp/2+p/4}r(p)2^{dp/2}3^p\sum\limits_{m+n=p}t^{n}\tau^{m}\s\s\mbox{for}\s\s t\in[0,\epsilon],\s\tau\in[0,\epsilon]
$$
and we will check that $h$ is bounded. Since
\begin{eqnarray*}
\sum\limits_{m+n=p}t^{n}\tau^{m}\leq\left(\sum\limits_{n=0}^pt^{n}\right)\left(\sum\limits_{n=0}^p\tau^{n}\right)\leq (1+p)^2\epsilon^2
\end{eqnarray*}
we have
\begin{eqnarray*}
h(t,\tau)\leq\epsilon^2\sum\limits_{p=1}^{\infty}B_pC^{dp/2+p/4}r(p)2^{dp/2}3^p(p+1)^2.
\end{eqnarray*}
Moreover, from
$$\limsup\limits_{p\rightarrow\infty}B_p^{1/p}<C^{-d/2-1/4}2^{-d/2}/3$$
it follows that
$$
\sum\limits_{p=1}^{\infty}B_pC^{dp/2+p/4}r(p)2^{dp/2}3^p(p+1)^2<\infty,
$$
where we have employ the fact
$$
\lim\limits_{p\rightarrow\infty}r(p)=\sum\limits_{i\in\mathbb{Z}}C^{-2^{-|i|}}<\infty.
$$
In summary, we have obtained
\begin{eqnarray*}
&&\sum\limits_{i\in\mathbb{Z}}\Bigg|\int_{\mathbb{R}}\int_{\mathbb{T}^d}\left\langle F\left(\Pi_+(D)\psi+\Pi_-(D)\varphi\right),\beta Q^{\mp}_i\Pi_{\pm}(D)P_j\phi\right\rangle_{\mathbb{C}^{d_0}}\,dxdt\Bigg|\nonumber\\
&\lesssim&2^{j(\max\{\lfloor d/2-d/b\rfloor,\lfloor d/2-d/a\rfloor\}+3)-3jd/2}\Big|\Big|P_j\phi\Big|\Big|_{S_{j+1}^{\mp}}h(||\Pi_+(D)\psi||_{S^{+,d/2}},||\Pi_-(D)\varphi||_{S^{-,d/2}}).
\end{eqnarray*}
At last, we shall prove
\begin{eqnarray}\label{5}
\Big|\Big|P_j\phi\Big|\Big|_{S_{j+1}^{\mp}}\lesssim2^{j(\max\{\lfloor d/2-d/a,\lfloor d/2-d/b\rfloor\}+1)}\Big|\Big|\tilde{P}_{j+1}\phi\Big|\Big|_{S_{j+1}^{\mp}}.
\end{eqnarray}
Indeed, from the definition it follows that
\begin{eqnarray*}
&&\Big|\Big|P_j\phi\Big|\Big|_{S_{j+1}^{\mp}}=||P_j\phi||_{L_t^{\infty}L_x^2}+||P_j\phi||_{\dot{X}^{\mp,1/2,\infty}}\\
&&+\sup\limits_{\substack{0\leq k'\leq j+1\\(d+2)(j+1)/(2d-2)\leq l\leq j+1}}\left\{2^{-(k'+j+1)/a}||P_j\phi||_{L_t^aL_x^b[l,k']}+2^{-(k'+j+1)/b}||P_j\phi||_{L_t^bL_x^a[l,k']}\right\}.
\end{eqnarray*}
For simplicity we only control
$$
||P_j\phi||_{L_t^aL_x^b[l,k']}.
$$
From the definition and Gagliardo-Nirenberg-Moser estimates it follows that
\begin{eqnarray*}
||P_j\phi||_{L_t^aL_x^b[l,k']}&=&\sum\limits_{\kappa\in\mathcal{K}_l}\sum\limits_{n\in\Xi_{k'}}||\Gamma_{k',n}P_{\kappa}P_j\phi||_{L_t^aL_x^b}\\
&\lesssim&\sum\limits_{\kappa\in\mathcal{K}_l}\sum\limits_{n\in\Xi_{k'}}\sum\limits_{|\alpha|=\lfloor d/2-d/b\rfloor+1}\Big|\Big|\Gamma_{k',n}P_{\kappa}P_j\phi\Big|\Big|^{\theta}_{L_t^aL_x^2}\Big|\Big|D^{\alpha}\Gamma_{k',n}P_{\kappa}P_j\phi\Big|\Big|^{1-\theta}_{L_t^aL_x^2}.
\end{eqnarray*}
Employing Bernstein-Type Lemma, i.e. Lemma 2.1 in chapter 2 of \cite{BCD} leads to
\begin{eqnarray*}
\Big|\Big|D^{\alpha}\Gamma_{k',n}P_{\kappa}P_j\phi\Big|\Big|_{L_t^aL_x^2}\leq C^{|\alpha|+1}2^{|\alpha|j}\cdot\Big|\Big|\Gamma_{k',n}P_{\kappa}P_j\phi\Big|\Big|_{L^a_tL^2_x}
\end{eqnarray*}
which implies
\begin{eqnarray*}
||P_j\phi||_{L_t^aL_x^b[l,k']}&\lesssim&2^{j(\lfloor d/2-d/b\rfloor+1)}\sum\limits_{\kappa\in\mathcal{K}_l}\sum\limits_{n\in\Xi_{k'}}\Big|\Big|\Gamma_{k',n}P_{\kappa}P_j\phi\Big|\Big|_{L_t^aL_x^2}\\
&\leq&2^{j(\lfloor d/2-d/b\rfloor+1)}\sum\limits_{\kappa\in\mathcal{K}_l}\sum\limits_{n\in\Xi_{k'}}\Big|\Big|\Gamma_{k',n}P_{\kappa}\tilde{P}_{j+1}\phi\Big|\Big|_{L_t^aL_x^2}\\
&\lesssim&2^{j(\lfloor d/2-d/b\rfloor+1)}\sum\limits_{\kappa\in\mathcal{K}_l}\sum\limits_{n\in\Xi_{k'}}\Big|\Big|\Gamma_{k',n}P_{\kappa}\tilde{P}_{j+1}\phi\Big|\Big|_{L_t^aL_x^b},
\end{eqnarray*}
where we have used the fact that flat torus has finite measure. The same idea also yields
\begin{eqnarray*}
||P_j\phi||_{L_t^bL_x^a[l,k']}\lesssim2^{j(\lfloor d/2-d/a\rfloor+1)}\sum\limits_{\kappa\in\mathcal{K}_l}\sum\limits_{n\in\Xi_{k'}}\Big|\Big|\Gamma_{k',n}P_{\kappa}\tilde{P}_{j+1}\phi\Big|\Big|_{L_t^bL_x^a}
\end{eqnarray*}
telling us that (\ref{5}) is true.

In conclusion, we have obtained
\begin{eqnarray*}
&&\sum\limits_{i\in\mathbb{Z}}\Bigg|\int_{\mathbb{R}}\int_{\mathbb{T}^d}\left\langle F\left(\Pi_+(D)\psi+\Pi_-(D)\varphi\right),\beta Q^{\mp}_i\Pi_{\pm}(D)P_j\phi\right\rangle_{\mathbb{C}^{d_0}}\,dxdt\Bigg|\nonumber\\
&\lesssim&2^{j(2\max\{\lfloor d/2-d/b\rfloor,\lfloor d/2-d/a\rfloor\}+4)-3jd/2}\Big|\Big|\tilde{P}_{j+1}\phi\Big|\Big|_{S_{j+1}^{\mp}}h(||\Pi_+(D)\psi||_{S^{+,d/2}},||\Pi_-(D)\varphi||_{S^{-,d/2}})\\
&\lesssim&\epsilon^2 2^{j(2\max\{\lfloor d/2-d/b\rfloor,\lfloor d/2-d/a\rfloor\}+4)-3jd/2}\Big|\Big|\tilde{P}_{j+1}\phi\Big|\Big|_{S_{j+1}^{\mp}}.
\end{eqnarray*}

Choose $a=b=4$ and
$$
y:=3d/2-2\max\{\lfloor d/2-d/b\rfloor,\lfloor d/2-d/a\rfloor\}-4
$$

Since $d\geq 9$, it is not difficult to get $-y<-d/2$, which means that $\mathcal{T}$ maps $Y^{d/2}$ into itself.

\subsection{Step 2: $\mathcal{T}$ is a contraction mapping}
\begin{eqnarray}\label{8}
\left\{ \begin{aligned}
         \mathcal{T}\psi^1(t)-\mathcal{T}\psi^2(t):&=\i\int_0^t e^{\i(t-s)\langle D\rangle}\Pi_+(D)\{\beta F(\Pi_+(D)\psi^1+\Pi_-(D)\varphi^1)\}\\
                                                   &-\i\int_0^t e^{\i(t-s)\langle D\rangle}\Pi_+(D)\{\beta F(\Pi_+(D)\psi^2+\Pi_-(D)\varphi^2)\}\\
         \mathcal{T}\varphi^1(t)-\mathcal{T}\varphi^2(t):&=\i\int_0^t e^{-\i(t-s)\langle D\rangle}\Pi_-(D)\{\beta F(\Pi_+(D)\psi^1+\Pi_-(D)\varphi^1)\}\\
                                                         &-\i\int_0^t e^{-\i(t-s)\langle D\rangle}\Pi_-(D)\{\beta F(\Pi_+(D)\psi^2+\Pi_-(D)\varphi^2)\}
\end{aligned} \right.
\end{eqnarray}
Because of Newton-Leibnitz formula we can get
$$
F(u^1)-F(u^2)=\int_0^1\frac{\partial F}{\partial u}(s(u^1-u^2)+u^2)(u^1-u^2)\,ds.
$$
Recalling (\ref{7}) yields
\begin{eqnarray*}
\frac{\partial F}{\partial u}(s(u^1-u^2)+u^2)(u^1-u^2)=\sum\limits_{|\mathbf{p}|=1}^{\infty}\sum\limits_{i=1}^{d_0}\mathbf{c_p}p_i\{s(u^1-u^2)+u^2\}^{(\mathbf{p}-e_i)^+}(u^1_i-u^2_i)
\end{eqnarray*}
where $\mathbf{p}=(p_1,p_2,\cdots,p_{d_0})$ and $(\mathbf{p}-e_i)^+$ means a $d_0$-dimensional vector whose the $i$-th component is $\max\{p_i-1,0\}$ and the others are $p_k$ for $k\not=i$. Expand the above expression as
$$
(u^1+u^2)^{(\mathbf{p}-e_i)^+}:=\sum\limits_{\mathbf{m}_i+\mathbf{n}_i=(\mathbf{p}-e_i)^+}c_{\mathbf{m}_i\mathbf{n}_i}(u^1)^{\mathbf{m}_i}(u^2)^{\mathbf{n}_i}\s\s\mbox{with}\s c_{\mathbf{m}_i\mathbf{n}_i}\in\mathbb{N}
$$
which implies
\begin{eqnarray}\label{9}
&&F(u^1)-F(u^2)=\int_0^1\frac{\partial F}{\partial u}(s(u^1-u^2)+u^2)(u^1-u^2)\,ds\\
&=&\sum\limits_{|\mathbf{p}|=1}^{\infty}\sum\limits_{i=1}^{d_0}\sum\limits_{\mathbf{m}_i+\mathbf{n}_i=(\mathbf{p}-e_i)^+}\mathbf{c_p}p_i\frac{c_{\mathbf{m}_i\mathbf{n}_i}}{|\mathbf{m}_i|+1}(u^1-u^2)^{\mathbf{m}_i}(u^2)^{\mathbf{n}_i}(u^1_i-u^2_i).\nonumber
\end{eqnarray}

From the definition it follows that
\begin{eqnarray*}
&&||\mathcal{T}\psi^1-\mathcal{T}\psi^2||_{S^{+,d/2}}=\sum\limits_{j=0}^{\infty}2^{dj/2}||P_j\mathcal{T}\psi^1-P_j\mathcal{T}\psi^2||_{S^+_{j+1}}\\
&\leq&\sum\limits_{j=0}^{\infty}2^{dj/2}\sup\limits_{g\in G^+_{j+1}}\Big|\int_{\mathbb{R}}\int_{\mathbb{T}^d}\langle P_j\Pi_+(D)\beta\{F(\Pi_+(D)\psi^1+\Pi_-(D)\varphi^1)\\
&&-F(\Pi_+(D)\psi^2+\Pi_-(D)\varphi^2)\},g \rangle_{\mathbb{C}^{d_0}}\,dxdt\Big|\\
&\leq&\sum\limits_{j=0}^{\infty}2^{dj/2}\sup\limits_{\phi\not=0}\Big|\int_{\mathbb{R}}\int_{\mathbb{T}^d}\langle P_j\Pi_+(D)\beta\{F(\Pi_+(D)\psi^1+\Pi_-(D)\varphi^1)\\
&&-F(\Pi_+(D)\psi^2+\Pi_-(D)\varphi^2)\},\tilde{P}_{j+1}\phi \rangle_{\mathbb{C}^{d_0}}\,dxdt\Big|/||\tilde{P}_{j+1}\phi||_{S^+_{j+1}}
\end{eqnarray*}
where we have used Lemma \ref{lem1}. Then it suffices to prove
\begin{eqnarray*}
&&\Big|\int_{\mathbb{R}}\int_{\mathbb{T}^d}\langle F(\Pi_+(D)\psi^1+\Pi_-(D)\varphi^1)-F(\Pi_+(D)\psi^2+\Pi_-(D)\varphi^2),\beta\Pi_{\mp}(D)P_j\phi \rangle_{\mathbb{C}^{d_0}}\,dxdt\Big|\\
&\lesssim&2^{-zj}||\tilde{P}_{j+1}\phi||_{S^{\pm}_{j+1}}\left\{||\psi^1-\psi^2||_{S^{+,d/2}}+||\varphi^1-\varphi^2||_{S^{-,d/2}}\right\}\\
&&g(||\Pi_+(D)(\psi^1-\psi^2)||_{S^{+,d/2}},||\Pi_-(D)(\varphi^1-\varphi^2)||_{S^{-,d/2}},||\Pi_+(D)\psi^2||_{S^{+,d/2}},||\Pi_-(D)\varphi^2||_{S^{-,d/2}})
\end{eqnarray*}
with $z>d/2$ and $0\leq g<1/2$ to be determined. Indeed, we have
\begin{eqnarray*}
&&\int_{\mathbb{R}}\int_{\mathbb{T}^d}\langle F(\Pi_+(D)\psi^1+\Pi_-(D)\varphi^1)-F(\Pi_+(D)\psi^2+\Pi_-(D)\varphi^2),\beta\Pi_{\mp}(D)P_j\phi \rangle_{\mathbb{C}^{d_0}}\,dxdt\\
&=&\sum\limits_{|\mathbf{p}|=1}^{\infty}\sum\limits_{i=1}^{d_0}\sum\limits_{\mathbf{m}_i+\mathbf{n}_i=(\mathbf{p}-e_i)^+}p_i\frac{c_{\mathbf{m}_i\mathbf{n}_i}}{|\mathbf{m}_i|+1}\int_{\mathbb{R}}\int_{\mathbb{T}^d}\langle\beta\Pi_{\mp}(D)P_j\phi,\mathbf{c_p}\rangle_{\mathbb{C}^{d_0}}\\ &&(\Pi_+(D)\psi^1-\Pi_+(D)\psi^2+\Pi_-(D)\varphi^1-\Pi_-(D)\varphi^2)^{\mathbf{m}_i}(\Pi_+(D)\psi^2+\Pi_-(D)\varphi^2)^{\mathbf{n}_i}\\
&&(\Pi_+(D)\psi_i^1-\Pi_+(D)\psi_i^2+\Pi_-(D)\varphi_i^1-\Pi_-(D)\varphi_i^2)\,dxdt.
\end{eqnarray*}
Employing the following formula
$$
(u^1+u^2)^{\mathbf{m}_i}:=\sum\limits_{\mathbf{k}_i+\mathbf{l}_i=\mathbf{m}_i}a_{\mathbf{k}_i\mathbf{l}_i}(u^1)^{\mathbf{k}_i}(u^2)^{\mathbf{l}_i}\s\s\mbox{with}\s a_{\mathbf{k}_i\mathbf{l}_i}\in\mathbb{N}
$$
and
$$
(u^1+u^2)^{\mathbf{n}_i}:=\sum\limits_{\mathbf{r}_i+\mathbf{s}_i=\mathbf{n}_i}b_{\mathbf{r}_i\mathbf{s}_i}(u^1)^{\mathbf{r}_i}(u^2)^{\mathbf{s}_i}\s\s\mbox{with}\s b_{\mathbf{r}_i\mathbf{s}_i}\in\mathbb{N}
$$
yields
\begin{eqnarray*}
&&\int_{\mathbb{R}}\int_{\mathbb{T}^d}\langle F(\Pi_+(D)\psi^1+\Pi_-(D)\varphi^1)-F(\Pi_+(D)\psi^2+\Pi_-(D)\varphi^2),\beta\Pi_{\mp}(D)P_j\phi \rangle_{\mathbb{C}^{d_0}}\,dxdt\\
&=&\sum\limits_{|\mathbf{p}|=1}^{\infty}\sum\limits_{i=1}^{d_0}\sum\limits_{\mathbf{m}_i+\mathbf{n}_i=(\mathbf{p}-e_i)^+}\sum\limits_{\mathbf{k}_i+\mathbf{l}_i=\mathbf{m}_i}\sum\limits_{\mathbf{r}_i+\mathbf{s}_i=\mathbf{n}_i}a_{\mathbf{k}_i\mathbf{l}_i}b_{\mathbf{r}_i\mathbf{s}_i}p_i\frac{c_{\mathbf{m}_i\mathbf{n}_i}}{|\mathbf{m}_i|+1}\int_{\mathbb{R}}\int_{\mathbb{T}^d}\langle\mathbf{c_p},\beta\Pi_{\mp}(D)P_j\phi\rangle_{\mathbb{C}^{d_0}}\\ &&(\Pi_+(D)\psi^1-\Pi_+(D)\psi^2)^{\mathbf{k}_i}(\Pi_-(D)\varphi^1-\Pi_-(D)\varphi^2)^{\mathbf{l}_i}(\Pi_+(D)\psi^2)^{\mathbf{r}_i}(\Pi_-(D)\varphi^2)^{\mathbf{s}_i}\\
&&(\Pi_+(D)\psi_i^1-\Pi_+(D)\psi_i^2+\Pi_-(D)\varphi_i^1-\Pi_-(D)\varphi_i^2)\,dxdt\\
&=&\sum\limits_{|\mathbf{p}|=1}^{\infty}\sum\limits_{i=1}^{d_0}\sum\limits_{\mathbf{m}_i+\mathbf{n}_i=(\mathbf{p}-e_i)^+}\sum\limits_{\mathbf{k}_i+\mathbf{l}_i=\mathbf{m}_i}\sum\limits_{\mathbf{r}_i+\mathbf{s}_i=\mathbf{n}_i}\\
&&\int_{\mathbb{R}}\int_{\mathbb{T}^d}\left\langle C^i_{\varsigma\overrightarrow{wvfb\tau\iota\sigma\varrho}\zeta\chi},\sum\limits_{\varsigma\in\mathbb{Z}}\beta Q^{\pm}_{\varsigma}\Pi_{\mp}(D)P_j\phi\right\rangle_{\mathbb{C}^{d_0}}\left\{\sum\limits_{w,v\in\mathbb{Z}}Q^+_wP_v\Pi_+(D)(\psi^1-\psi^2)\right\}^{\mathbf{k}_i}\\
&&\left\{\sum\limits_{f,b\in\mathbb{Z}}Q^{-}_fP_b\Pi_-(D)(\varphi^1-\varphi^2)\right\}^{\mathbf{l}_i}\left(\sum\limits_{\tau,\iota\in\mathbb{Z}}Q^+_\tau P_\iota\Pi_+(D)\psi^2\right)^{\mathbf{r}_i}\left(\sum\limits_{\sigma,\varrho\in\mathbb{Z}}Q^-_\sigma P_\varrho\Pi_-(D)\varphi^2\right)^{\mathbf{s}_i}\\
&&\left\{\sum\limits_{\zeta,\chi\in\mathbb{Z}}Q^+_\zeta P_\chi\Pi_+(D)(\psi_i^1-\psi_i^2)\right\}\,dxdt\\
&&+\sum\limits_{|\mathbf{p}|=1}^{\infty}\sum\limits_{i=1}^{d_0}\sum\limits_{\mathbf{m}_i+\mathbf{n}_i=(\mathbf{p}-e_i)^+}\sum\limits_{\mathbf{k}_i+\mathbf{l}_i=\mathbf{m}_i}\sum\limits_{\mathbf{r}_i+\mathbf{s}_i=\mathbf{n}_i}\\
&&\int_{\mathbb{R}}\int_{\mathbb{T}^d}\left\langle C^i_{\varsigma\overrightarrow{wvfb\tau\iota\sigma\varrho}\vartheta\varepsilon},\sum\limits_{\varsigma\in\mathbb{Z}}\beta Q^{\pm}_{\varsigma}\Pi_{\mp}(D)P_j\phi\right\rangle_{\mathbb{C}^{d_0}}\left\{\sum\limits_{w,v\in\mathbb{Z}}Q^+_wP_v\Pi_+(D)(\psi^1-\psi^2)\right\}^{\mathbf{k}_i}\\
&&\left\{\sum\limits_{f,b\in\mathbb{Z}}Q^{-}_fP_b\Pi_-(D)(\varphi^1-\varphi^2)\right\}^{\mathbf{l}_i}\left(\sum\limits_{\tau,\iota\in\mathbb{Z}}Q^+_\tau P_\iota\Pi_+(D)\psi^2\right)^{\mathbf{r}_i}\left(\sum\limits_{\sigma,\varrho\in\mathbb{Z}}Q^-_\sigma P_\varrho\Pi_-(D)\varphi^2\right)^{\mathbf{s}_i}\\
&&\left\{\sum\limits_{\vartheta,\varepsilon\in\mathbb{Z}}Q^-_\vartheta P_\varepsilon\Pi_-(D)(\varphi_i^1-\varphi_i^2)\right\}\,dxdt.
\end{eqnarray*}
Moreover, from the definition it follows that
\begin{eqnarray*}
&&\left\{\sum\limits_{w,v\in\mathbb{Z}}Q^+_wP_v\Pi_+(D)(\psi^1-\psi^2)\right\}^{\mathbf{k}_i}\\
&=&\sum\limits_{\substack{\vec{w}:=(\vec{w}_1,\vec{w}_2,\cdots,\vec{w}_{d_0})\\ \vec{v}:=(\vec{v}_1,\vec{v}_2,\cdots,\vec{v}_{d_0})}}Q^+_{\vec{w}_1}P_{\vec{v}_1}\Pi_+(D)(\psi^1_1-\psi^2_1)\cdots Q^+_{\vec{w}_{d_0}}P_{\vec{v}_{d_0}}\Pi_+(D)(\psi^1_{d_0}-\psi^2_{d_0}),
\end{eqnarray*}
\begin{eqnarray*}
&&\left\{\sum\limits_{f,b\in\mathbb{Z}}Q^-_fP_b\Pi_-(D)(\varphi^1-\varphi^2)\right\}^{\mathbf{l}_i}\\
&=&\sum\limits_{\substack{\vec{f}:=(\vec{f}_1,\vec{f}_2,\cdots,\vec{f}_{d_0})\\ \vec{b}:=(\vec{b}_1,\vec{b}_2,\cdots,\vec{b}_{d_0})}}Q^-_{\vec{f}_1}P_{\vec{b}_1}\Pi_-(D)(\varphi^1_1-\varphi^2_1)\cdots Q^-_{\vec{f}_{d_0}}P_{\vec{b}_{d_0}}\Pi_-(D)(\varphi^1_{d_0}-\varphi^2_{d_0}),
\end{eqnarray*}
\begin{eqnarray*}
\left(\sum\limits_{\tau,\iota\in\mathbb{Z}}Q^+_\tau P_\iota\Pi_+(D)\psi^2\right)^{\mathbf{r}_i}=\sum\limits_{\substack{\vec{\tau}:=(\vec{\tau}_1,\vec{\tau}_2,\cdots,\vec{\tau}_{d_0})\\ \vec{\iota}:=(\vec{\iota}_1,\vec{\iota}_2,\cdots,\vec{\iota}_{d_0})}}Q^+_{\vec{\tau}_1}P_{\vec{\iota}_1}\Pi_+(D)\psi^2_1\cdots Q^+_{\vec{\tau}_{d_0}}P_{\vec{\iota}_{d_0}}\Pi_+(D)\psi^2_{d_0},
\end{eqnarray*}
\begin{eqnarray*}
\left(\sum\limits_{\sigma,\varrho\in\mathbb{Z}}Q^-_\sigma P_\varrho\Pi_-(D)\varphi^2\right)^{\mathbf{s}_i}=\sum\limits_{\substack{\vec{\sigma}:=(\vec{\sigma}_1,\vec{\sigma}_2,\cdots,\vec{\sigma}_{d_0})\\ \vec{\varrho}:=(\vec{\varrho}_1,\vec{\varrho}_2,\cdots,\vec{\varrho}_{d_0})}}Q^-_{\vec{\sigma}_1}P_{\vec{\varrho}_1}\Pi_-(D)\varphi^2_1\cdots Q^-_{\vec{\sigma}_{d_0}}P_{\vec{\varrho}_{d_0}}\Pi_-(D)\varphi^2_{d_0}
\end{eqnarray*}
with
$$
\mathbf{k}_i:=(k_i^1,k_i^2,\cdots,k_i^{d_0}),\s\mathbf{l}_i:=(l_i^1,l_i^2,\cdots,l_i^{d_0}),\s\mathbf{r}_i:=(r_i^1,r_i^2,\cdots,r_i^{d_0}),\s\mathbf{s}_i:=(s_i^1,s_i^2,\cdots,s_i^{d_0})
$$
and
$$
Q^+_{\vec{w}_t}P_{\vec{v}_t}\Pi_+(D)(\psi^1_t-\psi^2_t):=Q^+_{w^1_t}P_{v^1_t}\Pi_+(D)(\psi^1_t-\psi^2_t)\cdots Q^+_{w^{k^t_i}_t}P_{v^{k^t_i}_t}\Pi_+(D)(\psi^1_t-\psi^2_t),
$$
$$
Q^-_{\vec{f}_t}P_{\vec{b}_t}\Pi_+(D)(\varphi^1_t-\varphi^2_t):=Q^-_{f^1_t}P_{b^1_t}\Pi_-(D)(\varphi^1_t-\varphi^2_t)\cdots Q^-_{f^{l^t_i}_t}P_{v^{l^t_i}_t}\Pi_+(D)(\psi^1_t-\psi^2_t),
$$
$$
Q^+_{\vec{\tau}_t}P_{\vec{\iota}_t}\Pi_+(D)\psi^2_t:=Q^+_{\tau^1_t}P_{\iota^1_t}\Pi_+(D)\psi^2_t\cdots Q^+_{\tau^{r^t_i}_t}P_{\iota^{r^t_i}_t}\Pi_+(D)\psi^2_t,
$$
$$
Q^-_{\vec{\sigma}_t}P_{\vec{\varrho}_t}\Pi_-(D)\varphi^2_t:=Q^+_{\sigma^1_t}P_{\varrho^1_t}\Pi_+(D)\varphi^2_t\cdots Q^+_{\sigma^{s^t_i}_t}P_{\varrho^{s^t_i}_t}\Pi_+(D)\varphi^2_t.
$$
Applying the same approach of (\ref{10}) we know that if
$$
\max\{||(\varsigma,\vec{w},\vec{f},\vec{\tau},\vec{\sigma},\zeta)||,||(\varsigma,\vec{w},\vec{f},\vec{\tau},\vec{\sigma},\vartheta)||,||(j,\vec{v},\vec{b},\vec{\iota},\vec{\varrho},\chi)||,||(j,\vec{v},\vec{b},\vec{\iota},\vec{\varrho},\varepsilon)||\}\geq 2,
$$
then
\begin{eqnarray*}
&&\int_{\mathbb{R}}\int_{\mathbb{T}^d}\left\langle C^i_{\varsigma\overrightarrow{wvfb\tau\iota\sigma\varrho}\zeta\chi},\beta Q^{\pm}_{\varsigma}\Pi_{\mp}(D)P_j\phi\right\rangle_{\mathbb{C}^{d_0}}Q^+_{\vec{w}_1}P_{\vec{v}_1}\Pi_+(D)(\psi^1_1-\psi^2_1)\\
&&\cdots Q^+_{\vec{w}_{d_0}}P_{\vec{v}_{d_0}}\Pi_+(D)(\psi^1_{d_0}-\psi^2_{d_0})Q^+_{\vec{\tau}_1}P_{\vec{\iota}_1}\Pi_+(D)\psi^2_1\cdots Q^+_{\vec{\tau}_{d_0}}P_{\vec{\iota}_{d_0}}\Pi_+(D)\psi^2_{d_0}\\
&&\cdot Q^-_{\vec{\sigma}_1}P_{\vec{\varrho}_1}\Pi_-(D)\varphi^2_1\cdots Q^-_{\vec{\sigma}_{d_0}}P_{\vec{\varrho}_{d_0}}\Pi_-(D)\varphi^2_{d_0}Q^-_{\vec{f}_1}P_{\vec{b}_1}\Pi_-(D)(\varphi^1_1-\varphi^2_1)\\
&&\cdots Q^-_{\vec{f}_{d_0}}P_{\vec{b}_{d_0}}\Pi_-(D)(\varphi^1_{d_0}-\varphi^2_{d_0})\cdot Q^+_\zeta P_\chi\Pi_+(D)(\psi_i^1-\psi_i^2)\,dxdt=0.
\end{eqnarray*}
and
\begin{eqnarray*}
&&\int_{\mathbb{R}}\int_{\mathbb{T}^d}\left\langle C^i_{\varsigma\overrightarrow{wvfb\tau\iota\sigma\varrho}\vartheta\varepsilon},\beta Q^{\pm}_{\varsigma}\Pi_{\mp}(D)P_j\phi\right\rangle_{\mathbb{C}^{d_0}}Q^+_{\vec{w}_1}P_{\vec{v}_1}\Pi_+(D)(\psi^1_1-\psi^2_1)\\
&&\cdots Q^+_{\vec{w}_{d_0}}P_{\vec{v}_{d_0}}\Pi_+(D)(\psi^1_{d_0}-\psi^2_{d_0})Q^+_{\vec{\tau}_1}P_{\vec{\iota}_1}\Pi_+(D)\psi^2_1\cdots Q^+_{\vec{\tau}_{d_0}}P_{\vec{\iota}_{d_0}}\Pi_+(D)\psi^2_{d_0}\\
&&\cdot Q^-_{\vec{\sigma}_1}P_{\vec{\varrho}_1}\Pi_-(D)\varphi^2_1\cdots Q^-_{\vec{\sigma}_{d_0}}P_{\vec{\varrho}_{d_0}}\Pi_-(D)\varphi^2_{d_0}Q^-_{\vec{f}_1}P_{\vec{b}_1}\Pi_-(D)(\varphi^1_1-\varphi^2_1)\\
&&\cdots Q^-_{\vec{f}_{d_0}}P_{\vec{b}_{d_0}}\Pi_-(D)(\varphi^1_{d_0}-\varphi^2_{d_0})\cdot Q^-_\vartheta P_\varepsilon\Pi_-(D)(\varphi_i^1-\varphi_i^2)\,dxdt=0,
\end{eqnarray*}
where
\begin{eqnarray*}
||(\varsigma,\vec{w},\vec{f},\vec{\tau},\vec{\sigma},\zeta)||:&=&\max\{|\varsigma-\zeta|,|\varsigma-w^a_t|,|\varsigma-f^b_t|,|\varsigma-\tau^c_t|,|\varsigma-\sigma^d_t|,|\zeta-w^a_t|,|\zeta-f^b_t|,|\zeta-\tau^c_t|,\\
&&|\zeta-\sigma^d_t|,|w_t^a-w_s^{\dot{a}}|,|f_t^b-f_s^{\dot{b}}|,|\tau_t^c-\tau_s^{\dot{c}}|,|\sigma_t^d-\sigma_s^{\dot{d}}|,|w_t^a-f_s^b|,|w_t^a-\tau_s^c|,\\
&&|w_t^a-\sigma_s^d|,|f_t^a-\tau_s^b|,|f_t^a-\sigma_s^b|,|\tau_t^a-\sigma_s^b|\,\,|1\leq t,s\leq d_0,1\leq a,\dot{a}\leq k_i^t,\\
&&1\leq b,\dot{b}\leq l_i^t,1\leq c,\dot{c}\leq r_i^t,1\leq d,\dot{d}\leq s_i^t\},
\end{eqnarray*}
\begin{eqnarray*}
||(\varsigma,\vec{w},\vec{f},\vec{\tau},\vec{\sigma},\vartheta)||:&=&\max\{|\varsigma-\vartheta|,|\varsigma-w^a_t|,|\varsigma-f^b_t|,|\varsigma-\tau^c_t|,|\varsigma-\sigma^d_t|,|\vartheta-w^a_t|,|\vartheta-f^b_t|,|\vartheta-\tau^c_t|,\\
&&|\vartheta-\sigma^d_t|,|w_t^a-w_s^{\dot{a}}|,|f_t^b-f_s^{\dot{b}}|,|\tau_t^c-\tau_s^{\dot{c}}|,|\sigma_t^d-\sigma_s^{\dot{d}}|,|w_t^a-f_s^b|,|w_t^a-\tau_s^c|,\\
&&|w_t^a-\sigma_s^d|,|f_t^a-\tau_s^b|,|f_t^a-\sigma_s^b|,|\tau_t^a-\sigma_s^b|\,\,|1\leq t,s\leq d_0,1\leq a,\dot{a}\leq k_i^t,\\
&&1\leq b,\dot{b}\leq l_i^t,1\leq c,\dot{c}\leq r_i^t,1\leq d,\dot{d}\leq s_i^t\},
\end{eqnarray*}
\begin{eqnarray*}
||(j,\vec{v},\vec{b},\vec{\iota},\vec{\varrho},\chi)||:&=&\max\{|j-\varrho|,|j-v^a_t|,|j-b^{\bar{b}}_t|,|j-\iota^c_t|,|j-\varrho^d_t|,|\chi-v^a_t|,|\chi-b^{\bar{b}}_t|,|\chi-\iota^c_t|,\\
&&|\chi-\varrho^d_t|,|v_t^a-v_s^{\dot{a}}|,|b_t^{\bar{b}}-b_s^{\dot{b}}|,|\iota_t^c-\iota_s^{\dot{c}}|,|\varrho_t^d-\varrho_s^{\dot{d}}|,|v_t^a-b_s^{\bar{b}}|,|v_t^a-\iota_s^c|,\\
&&|v_t^a-\varrho_s^d|,|b_t^a-\iota_s^{\bar{b}}|,|b_t^a-\varrho_s^{\bar{b}}|,|\iota_t^a-\varrho_s^{\bar{b}}|\,\,|1\leq t,s\leq d_0,1\leq a,\dot{a}\leq k_i^t,\\
&&1\leq \bar{b},\dot{b}\leq l_i^t,1\leq c,\dot{c}\leq r_i^t,1\leq d,\dot{d}\leq s_i^t\},
\end{eqnarray*}
\begin{eqnarray*}
||(j,\vec{v},\vec{b},\vec{\iota},\vec{\varrho},\varepsilon)||:&=&\max\{|j-\varrho|,|j-v^a_t|,|j-b^{\bar{b}}_t|,|j-\iota^c_t|,|j-\varrho^d_t|,|\varepsilon-v^a_t|,|\varepsilon-b^{\bar{b}}_t|,|\varepsilon-\iota^c_t|,\\
&&|\varepsilon-\varrho^d_t|,|v_t^a-v_s^{\dot{a}}|,|b_t^{\bar{b}}-b_s^{\dot{b}}|,|\iota_t^c-\iota_s^{\dot{c}}|,|\varrho_t^d-\varrho_s^{\dot{d}}|,|v_t^a-b_s^{\bar{b}}|,|v_t^a-\iota_s^c|,\\
&&|v_t^a-\varrho_s^d|,|b_t^a-\iota_s^{\bar{b}}|,|b_t^a-\varrho_s^{\bar{b}}|,|\iota_t^a-\varrho_s^{\bar{b}}|\,\,|1\leq t,s\leq d_0,1\leq a,\dot{a}\leq k_i^t,\\
&&1\leq \bar{b},\dot{b}\leq l_i^t,1\leq c,\dot{c}\leq r_i^t,1\leq d,\dot{d}\leq s_i^t\}.
\end{eqnarray*}
Using Young's Inequality leads to
\begin{eqnarray}\label{11}
&&\Big|\int_{\mathbb{R}}\int_{\mathbb{T}^d}\left\langle C^i_{\varsigma\overrightarrow{wvfb\tau\iota\sigma\varrho}\zeta\chi},\beta Q^{\pm}_{\varsigma}\Pi_{\mp}(D)P_j\phi\right\rangle_{\mathbb{C}^{d_0}}Q^+_{\vec{w}_1}P_{\vec{v}_1}\Pi_+(D)(\psi^1_1-\psi^2_1)\nonumber\\
&&\cdots Q^+_{\vec{w}_{d_0}}P_{\vec{v}_{d_0}}\Pi_+(D)(\psi^1_{d_0}-\psi^2_{d_0})Q^+_{\vec{\tau}_1}P_{\vec{\iota}_1}\Pi_+(D)\psi^2_1\cdots Q^+_{\vec{\tau}_{d_0}}P_{\vec{\iota}_{d_0}}\Pi_+(D)\psi^2_{d_0}\nonumber\\
&&\cdot Q^-_{\vec{\sigma}_1}P_{\vec{\varrho}_1}\Pi_-(D)\varphi^2_1\cdots Q^-_{\vec{\sigma}_{d_0}}P_{\vec{\varrho}_{d_0}}\Pi_-(D)\varphi^2_{d_0}Q^-_{\vec{f}_1}P_{\vec{b}_1}\Pi_-(D)(\varphi^1_1-\varphi^2_1)\nonumber\\
&&\cdots Q^-_{\vec{f}_{d_0}}P_{\vec{b}_{d_0}}\Pi_-(D)(\varphi^1_{d_0}-\varphi^2_{d_0})\cdot Q^+_\zeta P_\chi\Pi_+(D)(\psi_i^1-\psi_i^2)\Big|\\
&\leq&\Big|\Big|\left\langle C^i_{\varsigma\overrightarrow{wvfb\tau\iota\sigma\varrho}\zeta\chi},\beta Q^{\pm}_{\varsigma}\Pi_{\mp}(D)P_j\phi\right\rangle_{\mathbb{C}^{d_0}}\Big|\Big|_{L^{p_{\varsigma j\overrightarrow{wv\tau\iota\sigma\varrho fb}\zeta\chi}}_tL^{\hat{p}_{\varsigma j\overrightarrow{wv\tau\iota\sigma\varrho fb}\zeta\chi}}_x}\nonumber\\
&&\prod\limits_{t=1}^{d_0}\prod\limits_{r_t=1}^{k_i^t}\Big|\Big|Q^+_{w_t^{r_t}}P_{v_t^{r_t}}\Pi_+(D)(\psi^1_t-\psi_t^2)\Big|\Big|_{L^{p_t^{r_t}}_tL_x^{\hat{p}_t^{r_t}}}\prod\limits_{t=1}^{d_0}\prod\limits_{s_t=1}^{r_i^t}\Big|\Big|Q^+_{\tau_t^{s_t}}P_{\iota_t^{s_t}}\Pi_+(D)\psi_t^2\Big|\Big|_{L^{q_t^{s_t}}_tL_x^{\hat{q}_t^{s_t}}}\nonumber\\
&&\prod\limits_{t=1}^{d_0}\prod\limits_{a_t=1}^{s_i^t}\Big|\Big|Q^-_{\sigma_t^{a_t}}P_{\varrho_t^{a_t}}\Pi_-(D)\psi_t^2\Big|\Big|_{L^{m_t^{a_t}}_tL_x^{\hat{m}_t^{a_t}}}\prod\limits_{t=1}^{d_0}\prod\limits_{c_t=1}^{l_i^t}\Big|\Big|Q^-_{f_t^{c_t}}P_{b_t^{c_t}}\Pi_-(D)(\varphi^1_t-\varphi_t^2)\Big|\Big|_{L^{n_t^{c_t}}_tL_x^{\hat{n}_t^{c_t}}}\nonumber\\
&&\Big|\Big|Q^+_{\zeta}P_{\chi}\Pi_+(D)(\psi_i^1-\psi_i^2)\Big|\Big|_{L^{p_i}_tL^{\hat{p}_i}_x}\nonumber
\end{eqnarray}
with
\begin{eqnarray*}
1/p_{\varsigma j\overrightarrow{wv\tau\iota\sigma\varrho fb}\zeta\chi}+\sum\limits_{t=1}^{d_0}\sum\limits_{r_t=1}^{k^t_i}1/p_t^{r_t}+\sum\limits_{t=1}^{d_0}\sum\limits_{s_t=1}^{r^t_i}1/q_t^{s_t}+\sum\limits_{t=1}^{d_0}\sum\limits_{a_t=1}^{s^t_i}1/m_t^{a_t}+\sum\limits_{t=1}^{d_0}\sum\limits_{c_t=1}^{l^t_i}1/n_t^{c_t}+1/p_i=1
\end{eqnarray*}
and
\begin{eqnarray*}
1/\hat{p}_{\varsigma j\overrightarrow{wv\tau\iota\sigma\varrho fb}\zeta\chi}+\sum\limits_{t=1}^{d_0}\sum\limits_{r_t=1}^{k^t_i}1/\hat{p}_t^{r_t}+\sum\limits_{t=1}^{d_0}\sum\limits_{s_t=1}^{r^t_i}1/\hat{q}_t^{s_t}+\sum\limits_{t=1}^{d_0}\sum\limits_{a_t=1}^{s^t_i}1/\hat{m}_t^{a_t}+\sum\limits_{t=1}^{d_0}\sum\limits_{c_t=1}^{l^t_i}1/\hat{n}_t^{c_t}+1/\hat{p}_i=1
\end{eqnarray*}
to be determined and
\begin{eqnarray}\label{12}
&&\Big|\int_{\mathbb{R}}\int_{\mathbb{T}^d}\left\langle C^i_{\varsigma\overrightarrow{wvfb\tau\iota\sigma\varrho}\vartheta\varepsilon},\beta Q^{\mp}_{\varsigma}\Pi_{\pm}(D)P_j\phi\right\rangle_{\mathbb{C}^{d_0}}Q^+_{\vec{w}_1}P_{\vec{v}_1}\Pi_+(D)(\psi^1_1-\psi^2_1)\nonumber\\
&&\cdots Q^+_{\vec{w}_{d_0}}P_{\vec{v}_{d_0}}\Pi_+(D)(\psi^1_{d_0}-\psi^2_{d_0})Q^+_{\vec{\tau}_1}P_{\vec{\iota}_1}\Pi_+(D)\psi^2_1\cdots Q^+_{\vec{\tau}_{d_0}}P_{\vec{\iota}_{d_0}}\Pi_+(D)\psi^2_{d_0}\nonumber\\
&&\cdot Q^-_{\vec{\sigma}_1}P_{\vec{\varrho}_1}\Pi_-(D)\varphi^2_1\cdots Q^-_{\vec{\sigma}_{d_0}}P_{\vec{\varrho}_{d_0}}\Pi_-(D)\varphi^2_{d_0}Q^-_{\vec{f}_1}P_{\vec{b}_1}\Pi_-(D)(\varphi^1_1-\varphi^2_1)\nonumber\\
&&\cdots Q^-_{\vec{f}_{d_0}}P_{\vec{b}_{d_0}}\Pi_-(D)(\varphi^1_{d_0}-\varphi^2_{d_0})\cdot Q^-_\vartheta P_\varepsilon\Pi_-(D)(\varphi_i^1-\varphi_i^2)\Big|\\
&\leq&\Big|\Big|\left\langle C^i_{\varsigma\overrightarrow{wvfb\tau\iota\sigma\varrho}\vartheta\varepsilon},\beta Q^{\mp}_{\varsigma}\Pi_{\pm}(D)P_j\phi\right\rangle_{\mathbb{C}^{d_0}}\Big|\Big|_{L^{p_{\varsigma j\overrightarrow{wv\tau\iota\sigma\varrho fb}\vartheta\varepsilon}}_tL^{\hat{p}_{\varsigma j\overrightarrow{wv\tau\iota\sigma\varrho fb}\vartheta\varepsilon}}_x}\nonumber\\
&&\prod\limits_{t=1}^{d_0}\prod\limits_{r_t=1}^{k_i^t}\Big|\Big|Q^+_{w_t^{r_t}}P_{v_t^{r_t}}\Pi_+(D)(\psi^1_t-\psi_t^2)\Big|\Big|_{L^{g_t^{r_t}}_tL_x^{\hat{g}_t^{r_t}}}\prod\limits_{t=1}^{d_0}\prod\limits_{s_t=1}^{r_i^t}\Big|\Big|Q^+_{\tau_t^{s_t}}P_{\iota_t^{s_t}}\Pi_+(D)\psi_t^2\Big|\Big|_{L^{h_t^{s_t}}_tL_x^{\hat{h}_t^{s_t}}}\nonumber\\
&&\prod\limits_{t=1}^{d_0}\prod\limits_{a_t=1}^{s_i^t}\Big|\Big|Q^-_{\sigma_t^{a_t}}P_{\varrho_t^{a_t}}\Pi_-(D)\psi_t^2\Big|\Big|_{L^{d_t^{a_t}}_tL_x^{\hat{d}_t^{a_t}}}\prod\limits_{t=1}^{d_0}\prod\limits_{c_t=1}^{l_i^t}\Big|\Big|Q^-_{f_t^{c_t}}P_{b_t^{c_t}}\Pi_-(D)(\varphi^1_t-\varphi_t^2)\Big|\Big|_{L^{o_t^{c_t}}_tL_x^{\hat{o}_t^{c_t}}}\nonumber\\
&&\Big|\Big|Q^-_\vartheta P_\varepsilon\Pi_-(D)(\varphi_i^1-\varphi_i^2)\Big|\Big|_{L^{q_i}_tL^{\hat{q}_i}_x}\nonumber
\end{eqnarray}
with
\begin{eqnarray*}
1/p_{\varsigma j\overrightarrow{wv\tau\iota\sigma\varrho fb}\vartheta\varepsilon}+\sum\limits_{t=1}^{d_0}\sum\limits_{r_t=1}^{k^t_i}1/g_t^{r_t}+\sum\limits_{t=1}^{d_0}\sum\limits_{s_t=1}^{r^t_i}1/h_t^{s_t}+\sum\limits_{t=1}^{d_0}\sum\limits_{a_t=1}^{s^t_i}1/d_t^{a_t}+\sum\limits_{t=1}^{d_0}\sum\limits_{c_t=1}^{l^t_i}1/o_t^{c_t}+1/q_i=1
\end{eqnarray*}
and
\begin{eqnarray*}
1/\hat{p}_{\varsigma j\overrightarrow{wv\tau\iota\sigma\varrho fb}\vartheta\varepsilon}+\sum\limits_{t=1}^{d_0}\sum\limits_{r_t=1}^{k^t_i}1/\hat{g}_t^{r_t}+\sum\limits_{t=1}^{d_0}\sum\limits_{s_t=1}^{r^t_i}1/\hat{h}_t^{s_t}+\sum\limits_{t=1}^{d_0}\sum\limits_{a_t=1}^{s^t_i}1/\hat{d}_t^{a_t}+\sum\limits_{t=1}^{d_0}\sum\limits_{c_t=1}^{l^t_i}1/\hat{o}_t^{c_t}+1/\hat{q}_i=1
\end{eqnarray*}
to be determined. For simplicity we just compute (\ref{12}), since the same method still works for (\ref{11}). By Gagliardo-Nirenberg-Moser estimates we can get
\begin{eqnarray*}
&&\Big|\Big|\left\langle C^i_{\varsigma\overrightarrow{wvfb\tau\iota\sigma\varrho}\vartheta\varepsilon},\beta Q^{\mp}_{\varsigma}\Pi_{\pm}(D)P_j\phi\right\rangle_{\mathbb{C}^{d_0}}\Big|\Big|_{L^{p_{\varsigma j\overrightarrow{wv\tau\iota\sigma\varrho fb}\vartheta\varepsilon}}_tL^{\hat{p}_{\varsigma j\overrightarrow{wv\tau\iota\sigma\varrho fb}\vartheta\varepsilon}}_x}\\
&\lesssim&\sum\limits_{|\alpha|=\lceil\tilde{p}_{\varsigma j\overrightarrow{wv\tau\iota\sigma\varrho fb}\vartheta\varepsilon},\hat{p}_{\varsigma j\overrightarrow{wv\tau\iota\sigma\varrho fb}\vartheta\varepsilon}\rceil}\\
&&\Big|\Big|\left\langle C^i_{\varsigma\overrightarrow{wvfb\tau\iota\sigma\varrho}\vartheta\varepsilon},\beta Q^{\mp}_{\varsigma}\Pi_{\pm}(D)P_j\phi\right\rangle_{\mathbb{C}^{d_0}}\Big|\Big|^{1-\theta}_{L^{p_{\varsigma j\overrightarrow{wv\tau\iota\sigma\varrho fb}\vartheta\varepsilon}}_tL^{\tilde{p}_{\varsigma j\overrightarrow{wv\tau\iota\sigma\varrho fb}\vartheta\varepsilon}}_x}\\
&&\Big|\Big|\left\langle C^i_{\varsigma\overrightarrow{wvfb\tau\iota\sigma\varrho}\vartheta\varepsilon},\beta D^{\alpha}Q^{\mp}_{\varsigma}\Pi_{\pm}(D)P_j\phi\right\rangle_{\mathbb{C}^{d_0}}\Big|\Big|^{\theta}_{L^{p_{\varsigma j\overrightarrow{wv\tau\iota\sigma\varrho fb}\vartheta\varepsilon}}_tL^{\tilde{p}_{\varsigma j\overrightarrow{wv\tau\iota\sigma\varrho fb}\vartheta\varepsilon}}_x}
\end{eqnarray*}
with
$$
1/p_{\varsigma j\overrightarrow{wv\tau\iota\sigma\varrho fb}\vartheta\varepsilon}+1/\tilde{p}_{\varsigma j\overrightarrow{wv\tau\iota\sigma\varrho fb}\vartheta\varepsilon}=1/2
$$
and
\begin{eqnarray*}
\theta:&=&d/\left(\tilde{p}_{\varsigma j\overrightarrow{wv\tau\iota\sigma\varrho fb}\vartheta\varepsilon}\left\lceil\tilde{p}_{\varsigma j\overrightarrow{wv\tau\iota\sigma\varrho fb}\vartheta\varepsilon},\hat{p}_{\varsigma j\overrightarrow{wv\tau\iota\sigma\varrho fb}\vartheta\varepsilon}\right\rceil\right)\\
&&-d/\left(\hat{p}_{\varsigma j\overrightarrow{wv\tau\iota\sigma\varrho fb}\vartheta\varepsilon}\left\lceil\tilde{p}_{\varsigma j\overrightarrow{wv\tau\iota\sigma\varrho fb}\vartheta\varepsilon},\hat{p}_{\varsigma j\overrightarrow{wv\tau\iota\sigma\varrho fb}\vartheta\varepsilon}\right\rceil\right).
\end{eqnarray*}
Stepping modeled on (\ref{13}) and (\ref{14}) leads to
\begin{eqnarray*}
&&\Big|\Big|\left\langle C^i_{\varsigma\overrightarrow{wvfb\tau\iota\sigma\varrho}\vartheta\varepsilon},\beta D^{\alpha}Q^{\mp}_{\varsigma}\Pi_{\pm}(D)P_j\phi\right\rangle_{\mathbb{C}^{d_0}}\Big|\Big|_{L^{p_{\varsigma j\overrightarrow{wv\tau\iota\sigma\varrho fb}\vartheta\varepsilon}}_tL^{\tilde{p}_{\varsigma j\overrightarrow{wv\tau\iota\sigma\varrho fb}\vartheta\varepsilon}}_x}\\
&\lesssim&2^{2j/p_{\varsigma j\overrightarrow{wv\tau\iota\sigma\varrho fb}\vartheta\varepsilon}}\Big|\Big|\left\langle C^i_{\varsigma\overrightarrow{wvfb\tau\iota\sigma\varrho}\vartheta\varepsilon},\beta D^{\alpha}\Pi_{\mp}(D)P_j\phi\right\rangle_{\mathbb{C}^{d_0}}\Big|\Big|_{S_j^{\pm}}\\
&\leq&2^{2j/p_{\varsigma j\overrightarrow{wv\tau\iota\sigma\varrho fb}\vartheta\varepsilon}+|\alpha|j}C^{|\alpha|+1}\Big|\Big|\left\langle C^i_{\varsigma\overrightarrow{wvfb\tau\iota\sigma\varrho}\vartheta\varepsilon},\beta\Pi_{\mp}(D)P_j\phi\right\rangle_{\mathbb{C}^{d_0}}\Big|\Big|_{S_j^{\pm}},
\end{eqnarray*}
which implies
\begin{eqnarray*}
&&\Big|\Big|\left\langle C^i_{\varsigma\overrightarrow{wvfb\tau\iota\sigma\varrho}\vartheta\varepsilon},\beta Q^{\mp}_{\varsigma}\Pi_{\pm}(D)P_j\phi\right\rangle_{\mathbb{C}^{d_0}}\Big|\Big|_{L^{p_{\varsigma j\overrightarrow{wv\tau\iota\sigma\varrho fb}\vartheta\varepsilon}}_tL^{\hat{p}_{\varsigma j\overrightarrow{wv\tau\iota\sigma\varrho fb}\vartheta\varepsilon}}_x}\\
&\lesssim&2^{2j/p_{\varsigma j\overrightarrow{wv\tau\iota\sigma\varrho fb}\vartheta\varepsilon}+\left(d/\tilde{p}_{\varsigma j\overrightarrow{wv\tau\iota\sigma\varrho fb}\vartheta\varepsilon}-d/\hat{p}_{\varsigma j\overrightarrow{wv\tau\iota\sigma\varrho fb}\vartheta\varepsilon}\right)j}C^{d/\tilde{p}_{\varsigma j\overrightarrow{wv\tau\iota\sigma\varrho fb}\vartheta\varepsilon}-d/\hat{p}_{\varsigma j\overrightarrow{wv\tau\iota\sigma\varrho fb}\vartheta\varepsilon}+\theta}\\
&&\Big|\Big|\left\langle C^i_{\varsigma\overrightarrow{wvfb\tau\iota\sigma\varrho}\vartheta\varepsilon},\beta\Pi_{\pm}(D)P_j\phi\right\rangle_{\mathbb{C}^{d_0}}\Big|\Big|_{S_j^{\mp}}\\
&\lesssim&2^{2j/p_{\varsigma j\overrightarrow{wv\tau\iota\sigma\varrho fb}\vartheta\varepsilon}+\left(d/\tilde{p}_{\varsigma j\overrightarrow{wv\tau\iota\sigma\varrho fb}\vartheta\varepsilon}-d/\hat{p}_{\varsigma j\overrightarrow{wv\tau\iota\sigma\varrho fb}\vartheta\varepsilon}\right)j}C^{d/\tilde{p}_{\varsigma j\overrightarrow{wv\tau\iota\sigma\varrho fb}\vartheta\varepsilon}-d/\hat{p}_{\varsigma j\overrightarrow{wv\tau\iota\sigma\varrho fb}\vartheta\varepsilon}+\theta}\\
&&\Big|\Big|\left\langle C^i_{\varsigma\overrightarrow{wvfb\tau\iota\sigma\varrho}\vartheta\varepsilon},\beta\Pi_{\pm}(D)P_j\phi\right\rangle_{\mathbb{C}^{d_0}}\Big|\Big|_{S_{j+1}^{\mp}}.
\end{eqnarray*}
Employing the same approach of (\ref{15}) gives
\begin{eqnarray*}
&&\Big|\Big|Q^+_{w_t^{r_t}}P_{v_t^{r_t}}\Pi_+(D)(\psi^1_t-\psi_t^2)\Big|\Big|_{L^{g_t^{r_t}}_tL_x^{\hat{g}_t^{r_t}}}\\
&\lesssim& 2^{2j/g_t^{r_t}+\left(d/\tilde{g}_t^{r_t}-d/\hat{g}_t^{r_t}\right)j+(v_t^{r_t}+1-j)\max\{1/a,1/b\}}C^{d/\tilde{g}_t^{r_t}-d/\hat{g}_t^{r_t}+\theta_{r_t}}\\
&&\Big|\Big|P_{v_t^{r_t}}\Pi_+(D)(\psi^1_t-\psi_t^2)\Big|\Big|_{S^+_{v_t^{r_t}+1}}
\end{eqnarray*}
with
$$
1/g_t^{r_t}+1/\tilde{g}_t^{r_t}=1/2\s\mbox{and}\s\theta_{r_t}:=d/\left(\tilde{g}_t^{r_t}\left\lceil\tilde{g}_t^{r_t},\hat{g}_t^{r_t}\right\rceil\right)-d/\left(\hat{g}_t^{r_t}\left\lceil\tilde{g}_t^{r_t},\hat{g}_t^{r_t}\right\rceil\right)
$$
and
\begin{eqnarray*}
&&\Big|\Big|Q^+_{\tau_t^{s_t}}P_{\iota_t^{s_t}}\Pi_+(D)\psi_t^2\Big|\Big|_{L^{h_t^{s_t}}_tL_x^{\hat{h}_t^{s_t}}}\\
&\lesssim& 2^{2j/h_t^{s_t}+\left(d/\tilde{h}_t^{s_t}-d/\hat{h}_t^{s_t}\right)j+(\iota_t^{s_t}+1-j)\max\{1/a,1/b\}}C^{d/\tilde{h}_t^{s_t}-d/\hat{h}_t^{s_t}+\bar{\theta}_{s_t}}\\
&&\Big|\Big|P_{\iota_t^{s_t}}\Pi_+(D)\psi_t^2\Big|\Big|_{S^+_{\iota_t^{s_t}+1}}
\end{eqnarray*}
with
$$
1/h_t^{s_t}+1/\tilde{h}_t^{s_t}=1/2\s\mbox{and}\s\bar{\theta}_{s_t}:=d/\left(\tilde{h}_t^{s_t}\left\lceil\tilde{h}_t^{s_t},\hat{h}_t^{s_t}\right\rceil\right)-d/\left(\hat{h}_t^{s_t}\left\lceil\tilde{h}_t^{s_t},\hat{h}_t^{s_t}\right\rceil\right)
$$
and
\begin{eqnarray*}
&&\Big|\Big|Q^-_{\sigma_t^{a_t}}P_{\varrho_t^{a_t}}\Pi_-(D)\psi_t^2\Big|\Big|_{L^{d_t^{a_t}}_tL_x^{\hat{d}_t^{a_t}}}\\
&\lesssim&2^{2j/d_t^{a_t}+\left(d/\tilde{d}_t^{a_t}-d/\hat{d}_t^{a_t}\right)j+(\varrho_t^{a_t}+1-j)\max\{1/a,1/b\}}C^{d/\tilde{d}_t^{a_t}-d/\hat{d}_t^{a_t}+\pi_{a_t}}\\
&&\Big|\Big|P_{\varrho_t^{a_t}}\Pi_-(D)\psi_t^2\Big|\Big|_{S^-_{\varrho_t^{a_t}+1}}
\end{eqnarray*}
with
$$
1/d_t^{a_t}+1/\tilde{d}_t^{a_t}=1/2\s\mbox{and}\s\pi_{a_t}:=d/\left(\tilde{d}_t^{a_t}\left\lceil\tilde{d}_t^{a_t},\hat{d}_t^{a_t}\right\rceil\right)-d/\left(\hat{d}_t^{a_t}\left\lceil\tilde{d}_t^{a_t},\hat{d}_t^{a_t}\right\rceil\right)
$$
and
\begin{eqnarray*}
&&\Big|\Big|Q^-_{f_t^{c_t}}P_{b_t^{c_t}}\Pi_-(D)(\varphi^1_t-\varphi_t^2)\Big|\Big|_{L^{o_t^{c_t}}_tL_x^{\hat{o}_t^{c_t}}}\\
&\lesssim&2^{2j/o_t^{c_t}+\left(d/\tilde{o}_t^{c_t}-d/\hat{o}_t^{c_t}\right)j+(b_t^{c_t}+1-j)\max\{1/a,1/b\}}C^{d/\tilde{o}_t^{c_t}-d/\hat{o}_t^{c_t}+\bar{\pi}_{c_t}}\\
&&\Big|\Big|P_{b_t^{c_t}}\Pi_-(D)(\varphi^1_t-\varphi_t^2)\Big|\Big|_{S^-_{b_t^{c_t}+1}}
\end{eqnarray*}
with
$$
1/o_t^{c_t}+1/\tilde{o}_t^{c_t}=1/2\s\mbox{and}\s\bar{\pi}_{c_t}:=d/\left(\tilde{o}_t^{c_t}\left\lceil\tilde{o}_t^{c_t},\hat{o}_t^{c_t}\right\rceil\right)-d/\left(\hat{o}_t^{c_t}\left\lceil\tilde{o}_t^{c_t},\hat{o}_t^{c_t}\right\rceil\right)
$$
and
\begin{eqnarray*}
&&\Big|\Big|Q^-_{\vartheta}P_{\varepsilon}\Pi_-(D)(\varphi_i^1-\varphi_i^2)\Big|\Big|_{L^{q_i}_tL_x^{\hat{q}_i}}\\
&\lesssim&2^{2j/q_i+\left(d/\tilde{q}_i-d/\hat{q}_i\right)j+(\varepsilon+1-j)\max\{1/a,1/b\}}C^{d/\tilde{q}_i-d/\hat{q}_i+\bar{\pi}_i}\\
&&\Big|\Big|P_{\varepsilon}\Pi_-(D)(\varphi_i^1-\varphi_i^2)\Big|\Big|_{S^+_{\varepsilon+1}}
\end{eqnarray*}
with
$$
1/q_i+1/\tilde{q}_i=1/2\s\mbox{and}\s\bar{\pi}_i:=d/\left(\tilde{q}_i\left\lceil\tilde{q}_i,\hat{q}_i\right\rceil\right)-d/\left(\hat{q}_t^i\left\lceil\tilde{q}_i,\hat{q}_i\right\rceil\right).
$$
Therefore, we get
\begin{eqnarray*}
&&\Big|\int_{\mathbb{R}}\int_{\mathbb{T}^d}\left\langle C^i_{\varsigma\overrightarrow{wvfb\tau\iota\sigma\varrho}\vartheta\varepsilon},\beta Q^{\mp}_{\varsigma}\Pi_{\pm}(D)P_j\phi\right\rangle_{\mathbb{C}^{d_0}}Q^+_{\vec{w}_1}P_{\vec{v}_1}\Pi_+(D)(\psi^1_1-\psi^2_1)\nonumber\\
&&\cdots Q^+_{\vec{w}_{d_0}}P_{\vec{v}_{d_0}}\Pi_+(D)(\psi^1_{d_0}-\psi^2_{d_0})Q^+_{\vec{\tau}_1}P_{\vec{\iota}_1}\Pi_+(D)\psi^2_1\cdots Q^+_{\vec{\tau}_{d_0}}P_{\vec{\iota}_{d_0}}\Pi_+(D)\psi^2_{d_0}\nonumber\\
&&\cdot Q^-_{\vec{\sigma}_1}P_{\vec{\varrho}_1}\Pi_-(D)\varphi^2_1\cdots Q^-_{\vec{\sigma}_{d_0}}P_{\vec{\varrho}_{d_0}}\Pi_-(D)\varphi^2_{d_0}Q^-_{\vec{f}_1}P_{\vec{b}_1}\Pi_-(D)(\varphi^1_1-\varphi^2_1)\nonumber\\
&&\cdots Q^-_{\vec{f}_{d_0}}P_{\vec{b}_{d_0}}\Pi_-(D)(\varphi^1_{d_0}-\varphi^2_{d_0})\cdot Q^-_\vartheta P_\varepsilon\Pi_-(D)(\varphi_i^1-\varphi_i^2)\Big|\\
&\leq&2^{\Psi}C^{\Upsilon}\Big|\Big|\left\langle C^i_{\varsigma\overrightarrow{wvfb\tau\iota\sigma\varrho}\vartheta\varepsilon},\beta\Pi_{\pm}(D)P_j\phi\right\rangle_{\mathbb{C}^{d_0}}\Big|\Big|_{S_{j+1}^{\mp}}\\
&&\prod\limits_{t=1}^{d_0}\prod\limits_{r_t=1}^{k_i^t}\Big|\Big|P_{v_t^{r_t}}\Pi_+(D)(\psi^1_t-\psi_t^2)\Big|\Big|_{S^+_{v_t^{r_t}+1}}\prod\limits_{t=1}^{d_0}\prod\limits_{c_t=1}^{l_i^t}\Big|\Big|P_{b_t^{c_t}}\Pi_-(D)(\varphi^1_t-\varphi_t^2)\Big|\Big|_{S^-_{b_t^{c_t}+1}}\\
&&\prod\limits_{t=1}^{d_0}\prod\limits_{s_t=1}^{r_i^t}\Big|\Big|P_{\iota_t^{s_t}}\Pi_+(D)\psi_t^2\Big|\Big|_{S^+_{\iota_t^{s_t}+1}}\prod\limits_{t=1}^{d_0}\prod\limits_{a_t=1}^{s_i^t}\Big|\Big|P_{\varrho_t^{a_t}}\Pi_-(D)\psi_t^2\Big|\Big|_{S^-_{\varrho_t^{a_t}+1}}\\
&&\Big|\Big|P_{\varepsilon}\Pi_-(D)(\varphi_i^1-\varphi_i^2)\Big|\Big|_{S^+_{\varepsilon+1}}
\end{eqnarray*}
with
\begin{eqnarray*}
\Psi:&=&\frac{d}{2}\left(\sum\limits_{t=1}^{d_0}\sum\limits_{r_t=1}^{k_i^t}v_t^{r_t}+\sum\limits_{t=1}^{d_0}\sum\limits_{s_t=1}^{r_i^t}\iota_t^{s_t}+\sum\limits_{t=1}^{d_0}\sum\limits_{a_t=1}^{s_i^t}\varrho_t^{a_t}+\sum\limits_{t=1}^{d_0}\sum\limits_{c_t=1}^{l_i^t}b_t^{c_t}+\varepsilon\right)\\
&&+2j-3dj/2+d|(\mathbf{p}-e_i)^+|/2
\end{eqnarray*}
and
\begin{eqnarray*}
\Upsilon:&=&-3d/2+\frac{d}{2}|(\mathbf{p}-e_i)^+|+d/\left(\tilde{p}_{\varsigma j\overrightarrow{wv\tau\iota\sigma\varrho fb}\vartheta\varepsilon}\left\lceil\tilde{p}_{\varsigma j\overrightarrow{wv\tau\iota\sigma\varrho fb}\vartheta\varepsilon},\hat{p}_{\varsigma j\overrightarrow{wv\tau\iota\sigma\varrho fb}\vartheta\varepsilon}\right\rceil\right)\\
&&-d/\left(\hat{p}_{\varsigma j\overrightarrow{wv\tau\iota\sigma\varrho fb}\vartheta\varepsilon}\left\lceil\tilde{p}_{\varsigma j\overrightarrow{wv\tau\iota\sigma\varrho fb}\vartheta\varepsilon},\hat{p}_{\varsigma j\overrightarrow{wv\tau\iota\sigma\varrho fb}\vartheta\varepsilon}\right\rceil\right)\\
&&+\sum\limits_{t=1}^{d_0}\sum\limits_{r_t=1}^{k^t_i}\left\{d/\left(\tilde{g}_t^{r_t}\left\lceil\tilde{g}_t^{r_t},\hat{g}_t^{r_t}\right\rceil\right)-d/\left(\hat{g}_t^{r_t}\left\lceil\tilde{g}_t^{r_t},\hat{g}_t^{r_t}\right\rceil\right)\right\}\\
&&+\sum\limits_{t=1}^{d_0}\sum\limits_{s_t=1}^{r^t_i}\left\{d/\left(\tilde{h}_t^{s_t}\left\lceil\tilde{h}_t^{s_t},\hat{h}_t^{s_t}\right\rceil\right)-d/\left(\hat{h}_t^{s_t}\left\lceil\tilde{h}_t^{s_t},\hat{h}_t^{s_t}\right\rceil\right)\right\}\\
&&+\sum\limits_{t=1}^{d_0}\sum\limits_{a_t=1}^{s_i^t}\left\{d/\left(\tilde{d}_t^{a_t}\left\lceil\tilde{d}_t^{a_t},\hat{d}_t^{a_t}\right\rceil\right)-d/\left(\hat{d}_t^{a_t}\left\lceil\tilde{d}_t^{a_t},\hat{d}_t^{a_t}\right\rceil\right)\right\}\\
&&+\sum\limits_{t=1}^{d_0}\sum\limits_{c_t=1}^{l_i^t}\left\{d/\left(\tilde{o}_t^{c_t}\left\lceil\tilde{o}_t^{c_t},\hat{o}_t^{c_t}\right\rceil\right)-d/\left(\hat{o}_t^{c_t}\left\lceil\tilde{o}_t^{c_t},\hat{o}_t^{c_t}\right\rceil\right)\right\}\\
&&+d/\left(\tilde{q}_i\left\lceil\tilde{q}_i,\hat{q}_i\right\rceil\right)-d/\left(\hat{q}_i\left\lceil\tilde{q}_i,\hat{q}_i\right\rceil\right).
\end{eqnarray*}
Employing the same method of (\ref{16}) yields
\begin{eqnarray*}
&&\Big|\Big|\left\langle C^i_{\varsigma\overrightarrow{wvfb\tau\iota\sigma\varrho}\vartheta\varepsilon},\beta \Pi_{\pm}(D)P_j\phi\right\rangle_{\mathbb{C}^{d_0}}\Big|\Big|_{S_{j+1}^{\mp}}\\
&\lesssim&2^{j(\max\{\lfloor d/2-d/b\rfloor,\lfloor d/2-d/a\rfloor\}+1)}\Big|\Big|\left\langle\gamma^0 C^i_{\varsigma\overrightarrow{wvfb\tau\iota\sigma\varrho}\vartheta\varepsilon},P_j\phi\right\rangle_{\mathbb{C}^{d_0}}\Big|\Big|_{S_{j+1}^{\mp}}\nonumber\\
&&+2^{j(\max\{\lfloor d/2-d/b\rfloor,\lfloor d/2-d/a\rfloor\}+1)}\sum\limits_{\dot{k}=1}^d\Big|\Big|\left\langle\gamma^{\dot{k}}C^i_{\varsigma\overrightarrow{wvfb\tau\iota\sigma\varrho}\vartheta\varepsilon},P_j\phi\right\rangle_{\mathbb{C}^{d_0}}\Big|\Big|_{S_{j+1}^{\mp}}\nonumber\\
&&+2^{j(\max\{\lfloor d/2-d/b\rfloor,\lfloor d/2-d/a\rfloor\}+1)}\Big|\Big|\left\langle C^i_{\varsigma\overrightarrow{wvfb\tau\iota\sigma\varrho}\vartheta\varepsilon},P_j\phi\right\rangle_{\mathbb{C}^{d_0}}\Big|\Big|_{S_{j+1}^{\mp}}\nonumber\\
&\leq&2^{j(\max\{\lfloor d/2-d/b\rfloor,\lfloor d/2-d/a\rfloor\}+1)}\Big|\Big|P_j\phi\Big|\Big|_{S_{j+1}^{\mp}}\\
&&\left\{\Big|\gamma^0C^i_{\varsigma\overrightarrow{wvfb\tau\iota\sigma\varrho}\vartheta\varepsilon}\Big|+\sum\limits_{\dot{k}=1}^d\Big|\gamma^{\dot{k}}C^i_{\varsigma\overrightarrow{wvfb\tau\iota\sigma\varrho}\vartheta\varepsilon}\Big|+\Big|C^i_{\varsigma\overrightarrow{wvfb\tau\iota\sigma\varrho}\vartheta\varepsilon}\Big|\right\}\\
&\lesssim&2^{2j(\max\{\lfloor d/2-d/a,\lfloor d/2-d/b\rfloor\}+1)}\Big|\Big|\tilde{P}_{j+1}\phi\Big|\Big|_{S_{j+1}^{\mp}}\\
&&\left\{\Big|\gamma^0C^i_{\varsigma\overrightarrow{wvfb\tau\iota\sigma\varrho}\vartheta\varepsilon}\Big|+\sum\limits_{\dot{k}=1}^d\Big|\gamma^{\dot{k}}C^i_{\varsigma\overrightarrow{wvfb\tau\iota\sigma\varrho}\vartheta\varepsilon}\Big|+\Big|C^i_{\varsigma\overrightarrow{wvfb\tau\iota\sigma\varrho}\vartheta\varepsilon}\Big|\right\}.
\end{eqnarray*}
Now choose
\begin{eqnarray*}
1/g_t^{r_t}=1/\hat{g}_t^{r_t}:=\frac{2^{-|w_t^{r_t}|-1}}{|(\mathbf{p}-e_i)^+|+1}\s\mbox{and}\s 1/h_t^{s_t}=1/\hat{h}_t^{s_t}:=\frac{2^{-|\tau_t^{s_t}|-1}}{|(\mathbf{p}-e_i)^+|+1}
\end{eqnarray*}
\begin{eqnarray*}
1/d_t^{a_t}=1/\hat{d}_t^{a_t}:=\frac{2^{-|\sigma_t^{a_t}|-1}}{|(\mathbf{p}-e_i)^+|+1}\s\mbox{and}\s 1/o_t^{c_t}=1/\hat{o}_t^{c_t}:=\frac{2^{-|f_t^{c_t}|-1}}{|(\mathbf{p}-e_i)^+|+1}
\end{eqnarray*}
and
$$
1/q_i=1/\hat{q}_i:=\frac{2^{-|\vartheta|-1}}{|(\mathbf{p}-e_i)^+|+1}
$$
and
\begin{eqnarray*}
1/p_{\varsigma j\overrightarrow{wv\tau\iota\sigma\varrho fb}\vartheta\varepsilon}=1/\hat{p}_{\varsigma j\overrightarrow{wv\tau\iota\sigma\varrho fb}\vartheta\varepsilon}&:=&1-\frac{\sum\limits_{t=1}^{d_0}\sum\limits_{r_t=1}^{k^t_i}2^{-|w_t^{r_t}|-1}+\sum\limits_{t=1}^{d_0}\sum\limits_{s_t=1}^{r^t_i}2^{-|\tau_t^{s_t}|-1}}{|(\mathbf{p}-e_i)^+|+1}\\
&&-\frac{\sum\limits_{t=1}^{d_0}\sum\limits_{a_t=1}^{s_i^t}2^{-|\sigma_t^{a_t}|-1}+\sum\limits_{t=1}^{d_0}\sum\limits_{c_t=1}^{l_i^t}2^{-|f_t^{c_t}|-1}+2^{-|\vartheta|-1}}{|(\mathbf{p}-e_i)^+|+1}.
\end{eqnarray*}
Then
\begin{eqnarray*}
\Upsilon&=&\frac{-3/2+\frac{\sum\limits_{t=1}^{d_0}\sum\limits_{r_t=1}^{k^t_i}2^{-|w_t^{r_t}|}+\sum\limits_{t=1}^{d_0}\sum\limits_{s_t=1}^{r^t_i}2^{-|\tau_t^{s_t}|}+\sum\limits_{t=1}^{d_0}\sum\limits_{a_t=1}^{s_i^t}2^{-|\sigma_t^{a_t}|}+\sum\limits_{t=1}^{d_0}\sum\limits_{c_t=1}^{l_i^t}2^{-|f_t^{c_t}|}+2^{-|\vartheta|}}{|(\mathbf{p}-e_i)^+|+1}}{2+\left\lceil-3/2+\frac{\sum\limits_{t=1}^{d_0}\sum\limits_{r_t=1}^{k^t_i}2^{-|w_t^{r_t}|}+\sum\limits_{t=1}^{d_0}\sum\limits_{s_t=1}^{r^t_i}2^{-|\tau_t^{s_t}|}+\sum\limits_{t=1}^{d_0}\sum\limits_{a_t=1}^{s_i^t}2^{-|\sigma_t^{a_t}|}+\sum\limits_{t=1}^{d_0}\sum\limits_{c_t=1}^{l_i^t}2^{-|f_t^{c_t}|}+2^{-|\vartheta|}}{|(\mathbf{p}-e_i)^+|+1}\right\rceil}\\
&&+\sum\limits_{t=1}^{d_0}\sum\limits_{r_t=1}^{k^t_i}\frac{1/2-\frac{2^{-|w_t^{r_t}|}}{1+|(\mathbf{p}-e_i)^+|}}{2+\left\lceil1/2-\frac{2^{-|w_t^{r_t}|}}{1+|(\mathbf{p}-e_i)^+|}\right\rceil}+\sum\limits_{t=1}^{d_0}\sum\limits_{s_t=1}^{r^t_i}\frac{1/2-\frac{2^{-|\tau_t^{s_t}|}}{1+|(\mathbf{p}-e_i)^+|}}{2+\left\lceil1/2-\frac{2^{-|\tau_t^{s_t}|}}{1+|(\mathbf{p}-e_i)^+|}\right\rceil}\\
&&+\sum\limits_{t=1}^{d_0}\sum\limits_{a_t=1}^{s_i^t}\frac{1/2-\frac{2^{-|\sigma_t^{a_t}|}}{|(\mathbf{p}-e_i)^+|+1}}{2+\left\lceil1/2-\frac{2^{-|\sigma_t^{a_t}|}}{|(\mathbf{p}-e_i)^+|+1}\right\rceil}+\sum\limits_{t=1}^{d_0}\sum\limits_{c_t=1}^{l_i^t}\frac{1/2-\frac{2^{-|f_t^{c_t}|}}{|(\mathbf{p}-e_i)^+|+1}}{2+\left\lceil 1/2-\frac{2^{-|f_t^{c_t}|}}{1+|(\mathbf{p}-e_i)^+|}\right\rceil}\\
&&-3d/2+\frac{d}{2}|(\mathbf{p}-e_i)^+|+\frac{1/2-\frac{2^{-|\vartheta|}}{|(\mathbf{p}-e_i)^+|+1}}{2+\left\lceil 1/2-\frac{2^{-|\vartheta|}}{1+|(\mathbf{p}-e_i)^+|}\right\rceil}\\
&\leq&-\sum\limits_{t=1}^{d_0}\sum\limits_{r_t=1}^{k^t_i}\frac{2^{-|w_t^{r_t}|-1}}{1+|(\mathbf{p}-e_i)^+|}-\sum\limits_{t=1}^{d_0}\sum\limits_{s_t=1}^{r^t_i}\frac{2^{-|\tau_t^{s_t}|-1}}{1+|(\mathbf{p}-e_i)^+|}-\frac{2^{-|\vartheta|-1}}{1+|(\mathbf{p}-e_i)^+|}\\
&&-\sum\limits_{t=1}^{d_0}\sum\limits_{a_t=1}^{s_i^t}\frac{2^{-|\sigma_t^{a_t}|-1}}{|(\mathbf{p}-e_i)^+|+1}-\sum\limits_{t=1}^{d_0}\sum\limits_{c_t=1}^{l_i^t}\frac{2^{-|f_t^{c_t}|-1}}{1+|(\mathbf{p}-e_i)^+|}\\
&&+|(\mathbf{p}-e_i)^+|/4-3d/2+\frac{d}{2}|(\mathbf{p}-e_i)^+|+5/4.
\end{eqnarray*}
Hence we have
\begin{eqnarray*}
&&\Big|\int_{\mathbb{R}}\int_{\mathbb{T}^d}\left\langle C^i_{\varsigma\overrightarrow{wvfb\tau\iota\sigma\varrho}\vartheta\varepsilon},\beta Q^{\mp}_{\varsigma}\Pi_{\pm}(D)P_j\phi\right\rangle_{\mathbb{C}^{d_0}}Q^+_{\vec{w}_1}P_{\vec{v}_1}\Pi_+(D)(\psi^1_1-\psi^2_1)\nonumber\\
&&\cdots Q^+_{\vec{w}_{d_0}}P_{\vec{v}_{d_0}}\Pi_+(D)(\psi^1_{d_0}-\psi^2_{d_0})Q^+_{\vec{\tau}_1}P_{\vec{\iota}_1}\Pi_+(D)\psi^2_1\cdots Q^+_{\vec{\tau}_{d_0}}P_{\vec{\iota}_{d_0}}\Pi_+(D)\psi^2_{d_0}\nonumber\\
&&\cdot Q^-_{\vec{\sigma}_1}P_{\vec{\varrho}_1}\Pi_-(D)\varphi^2_1\cdots Q^-_{\vec{\sigma}_{d_0}}P_{\vec{\varrho}_{d_0}}\Pi_-(D)\varphi^2_{d_0}Q^-_{\vec{f}_1}P_{\vec{b}_1}\Pi_-(D)(\varphi^1_1-\varphi^2_1)\nonumber\\
&&\cdots Q^-_{\vec{f}_{d_0}}P_{\vec{b}_{d_0}}\Pi_-(D)(\varphi^1_{d_0}-\varphi^2_{d_0})\cdot Q^-_\vartheta P_\varepsilon\Pi_-(D)(\varphi_i^1-\varphi_i^2)\Big|\\
&\leq&2^{2j-3dj/2+d|(\mathbf{p}-e_i)^+|/2}C^{\Upsilon}\Big|\Big|\left\langle C^i_{\varsigma\overrightarrow{wvfb\tau\iota\sigma\varrho}\vartheta\varepsilon},\beta\Pi_{\pm}(D)P_j\phi\right\rangle_{\mathbb{C}^{d_0}}\Big|\Big|_{S_{j+1}^{\mp}}\\
&&\prod\limits_{t=1}^{d_0}\prod\limits_{r_t=1}^{k_i^t}\left\{2^{v_t^{r_t}d/2}\Big|\Big|P_{v_t^{r_t}}\Pi_+(D)(\psi^1_t-\psi_t^2)\Big|\Big|_{S^+_{v_t^{r_t}+1}}\right\}\\
&&\prod\limits_{t=1}^{d_0}\prod\limits_{c_t=1}^{l_i^t}\left\{2^{b_t^{c_t}d/2}\Big|\Big|P_{b_t^{c_t}}\Pi_-(D)(\varphi^1_t-\varphi_t^2)\Big|\Big|_{S^-_{b_t^{c_t}+1}}\right\}\\
&&\prod\limits_{t=1}^{d_0}\prod\limits_{s_t=1}^{r_i^t}\left\{2^{\iota_t^{s_t}d/2}\Big|\Big|P_{\iota_t^{s_t}}\Pi_+(D)\psi_t^2\Big|\Big|_{S^+_{\iota_t^{s_t}+1}}\right\}\prod\limits_{t=1}^{d_0}\prod\limits_{a_t=1}^{s_i^t}\left\{2^{\varrho_t^{a_t}d/2}\Big|\Big|P_{\varrho_t^{a_t}}\Pi_-(D)\psi_t^2\Big|\Big|_{S^-_{\varrho_t^{a_t}+1}}\right\}\\
&&\left\{2^{\varepsilon d/2}\Big|\Big|P_{\varepsilon}\Pi_-(D)(\varphi_i^1-\varphi_i^2)\Big|\Big|_{S^+_{\varepsilon+1}}\right\}.
\end{eqnarray*}
Stepping modeled on (\ref{16}) and (\ref{5}) leads to
\begin{eqnarray*}
&&\Big|\Big|\left\langle C^i_{\varsigma\overrightarrow{wvfb\tau\iota\sigma\varrho}\vartheta\varepsilon},\beta\Pi_{\pm}(D)P_j\phi\right\rangle_{\mathbb{C}^{d_0}}\Big|\Big|_{S_{j+1}^{\mp}}\\
&\lesssim&2^{j(\max\{\lfloor d/2-d/b\rfloor,\lfloor d/2-d/a\rfloor\}+1)}\Big|\Big|\left\langle\gamma^0C^i_{\varsigma\overrightarrow{wvfb\tau\iota\sigma\varrho}\vartheta\varepsilon},P_j\phi\right\rangle_{\mathbb{C}^{d_0}}\Big|\Big|_{S_{j+1}^{\mp}}\nonumber\\
&&+2^{j(\max\{\lfloor d/2-d/b\rfloor,\lfloor d/2-d/a\rfloor\}+1)}\sum\limits_{\dot{k}=1}^d\Big|\Big|\left\langle\gamma^{\dot{k}}C^i_{\varsigma\overrightarrow{wvfb\tau\iota\sigma\varrho}\vartheta\varepsilon},P_j\phi\right\rangle_{\mathbb{C}^{d_0}}\Big|\Big|_{S_{j+1}^{\mp}}\nonumber\\
&&+2^{j(\max\{\lfloor d/2-d/b\rfloor,\lfloor d/2-d/a\rfloor\}+1)}\Big|\Big|\left\langle C^i_{\varsigma\overrightarrow{wvfb\tau\iota\sigma\varrho}\vartheta\varepsilon},P_j\phi\right\rangle_{\mathbb{C}^{d_0}}\Big|\Big|_{S_{j+1}^{\mp}}\\
&\leq&2^{j(\max\{\lfloor d/2-d/b\rfloor,\lfloor d/2-d/a\rfloor\}+1)}\Big|\Big|P_j\phi\Big|\Big|_{S_{j+1}^{\mp}}\\
&&\left\{\Big|\gamma^0C^i_{\varsigma\overrightarrow{wvfb\tau\iota\sigma\varrho}\vartheta\varepsilon}\Big|+\sum\limits_{\dot{k}=1}^d\Big|\gamma^{\dot{k}}C^i_{\varsigma\overrightarrow{wvfb\tau\iota\sigma\varrho}\vartheta\varepsilon}\Big|+\Big|C^i_{\varsigma\overrightarrow{wvfb\tau\iota\sigma\varrho}\vartheta\varepsilon}\Big|\right\}\\
&\lesssim&2^{j(2\max\{\lfloor d/2-d/b\rfloor,\lfloor d/2-d/a\rfloor\}+2)}\Big|\Big|\tilde{P}_{j+1}\phi\Big|\Big|_{S_{j+1}^{\mp}}\\
&&\left\{\Big|\gamma^0C^i_{\varsigma\overrightarrow{wvfb\tau\iota\sigma\varrho}\vartheta\varepsilon}\Big|+\sum\limits_{\dot{k}=1}^d\Big|\gamma^{\dot{k}}C^i_{\varsigma\overrightarrow{wvfb\tau\iota\sigma\varrho}\vartheta\varepsilon}\Big|+\Big|C^i_{\varsigma\overrightarrow{wvfb\tau\iota\sigma\varrho}\vartheta\varepsilon}\Big|\right\}.
\end{eqnarray*}
Define
\begin{eqnarray*}
\sum:=\sum\limits_{\vec{w}=\varsigma-1}^{\varsigma+1}\sum\limits_{\vec{f}=\varsigma-1}^{\varsigma+1}\sum\limits_{\vec{\tau}=\varsigma-1}^{\varsigma+1}\sum\limits_{\vec{\sigma}=\varsigma-1}^{\varsigma+1}\sum\limits_{\vartheta=\varsigma-1}^{\varsigma+1},
\end{eqnarray*}
where
\begin{eqnarray*}
\sum\limits_{\vec{w}=\varsigma-1}^{\varsigma+1}:=\sum\limits_{\vec{w}_1=\varsigma-1}^{\varsigma+1}\sum\limits_{\vec{w}_2=\varsigma-1}^{\varsigma+1}\cdots\sum\limits_{\vec{w}_{d_0}=\varsigma-1}^{\varsigma+1}\s\mbox{with}\s\sum\limits_{\vec{w}_t=\varsigma-1}^{\varsigma+1}:=\sum\limits_{w^1_t=\varsigma-1}^{\varsigma+1}\sum\limits_{w^2_t=\varsigma-1}^{\varsigma+1}\cdots\sum\limits_{w^{k_i^t}_t=\varsigma-1}^{\varsigma+1}
\end{eqnarray*}
and
\begin{eqnarray*}
\sum\limits_{\vec{f}=\varsigma-1}^{\varsigma+1}:=\sum\limits_{\vec{f}_1=\varsigma-1}^{\varsigma+1}\sum\limits_{\vec{f}_2=\varsigma-1}^{\varsigma+1}\cdots\sum\limits_{\vec{f}_{d_0}=\varsigma-1}^{\varsigma+1}\s\mbox{with}\s\sum\limits_{\vec{f}_t=\varsigma-1}^{\varsigma+1}:=\sum\limits_{f^1_t=\varsigma-1}^{\varsigma+1}\sum\limits_{f^2_t=\varsigma-1}^{\varsigma+1}\cdots\sum\limits_{f^{l_i^t}_t=\varsigma-1}^{\varsigma+1}
\end{eqnarray*}
and
\begin{eqnarray*}
\sum\limits_{\vec{\tau}=\varsigma-1}^{\varsigma+1}:=\sum\limits_{\vec{\tau}_1=\varsigma-1}^{\varsigma+1}\sum\limits_{\vec{\tau}_2=\varsigma-1}^{\varsigma+1}\cdots\sum\limits_{\vec{\tau}_{d_0}=\varsigma-1}^{\varsigma+1}\s\mbox{with}\s\sum\limits_{\vec{\tau}_t=\varsigma-1}^{\varsigma+1}:=\sum\limits_{\tau^1_t=\varsigma-1}^{\varsigma+1}\sum\limits_{\tau^2_t=\varsigma-1}^{\varsigma+1}\cdots\sum\limits_{\tau^{r_i^t}_t=\varsigma-1}^{\varsigma+1}
\end{eqnarray*}
and
\begin{eqnarray*}
\sum\limits_{\vec{\sigma}=\varsigma-1}^{\varsigma+1}:=\sum\limits_{\vec{\sigma}_1=\varsigma-1}^{\varsigma+1}\sum\limits_{\vec{\sigma}_2=\varsigma-1}^{\varsigma+1}\cdots\sum\limits_{\vec{\sigma}_{d_0}=\varsigma-1}^{\varsigma+1}\s\mbox{with}\s\sum\limits_{\vec{\sigma}_t=\varsigma-1}^{\varsigma+1}:=\sum\limits_{\sigma^1_t=\varsigma-1}^{\varsigma+1}\sum\limits_{\sigma^2_t=\varsigma-1}^{\varsigma+1}\cdots\sum\limits_{\sigma^{s_i^t}_t=\varsigma-1}^{\varsigma+1}.
\end{eqnarray*}
Moreover, we define
$$
B^i_{\overrightarrow{vb\iota\varrho}\varepsilon}:=\sup\limits_{\varsigma\in\mathbb{Z}}\sum\left\{\Big|C^i_{\varsigma\overrightarrow{wvfb\tau\iota\sigma\varrho}\vartheta\varepsilon}\Big|+\Big|\gamma^0C^i_{\varsigma\overrightarrow{wvfb\tau\iota\sigma\varrho}\vartheta\varepsilon}\Big|+\sum\limits_{\dot{k}=1}^d\Big|\gamma^{\dot{k}} C^i_{\varsigma\overrightarrow{wvfb\tau\iota\sigma\varrho}\vartheta\varepsilon}\Big|\right\}
$$
and
$$
B^i_{\mathbf{klrs}}:=\sup\limits_{\vec{v}}\sup\limits_{\vec{b}}\sup\limits_{\vec{\iota}}\sup\limits_{\vec{\varrho}}\sup\limits_{\varepsilon}B^i_{\overrightarrow{vb\iota\varrho}\varepsilon}
$$
with
$$
\sup\limits_{\vec{v}}:=\sup\limits_{\vec{v}_1}\sup\limits_{\vec{v}_2}\cdots\sup\limits_{\vec{v}_{d_0}}\s\mbox{and}\s \sup\limits_{\vec{b}}:=\sup\limits_{\vec{b}_1}\sup\limits_{\vec{b}_2}\cdots\sup\limits_{\vec{b}_{d_0}}
$$
and
$$
\sup\limits_{\vec{\varrho}}:=\sup\limits_{\vec{\varrho}_1}\sup\limits_{\vec{\varrho}_2}\cdots\sup\limits_{\vec{\varrho}_{d_0}}\s\mbox{and}\s \sup\limits_{\vec{\iota}}:=\sup\limits_{\vec{\iota}_1}\sup\limits_{\vec{\iota}_2}\cdots\sup\limits_{\vec{\iota}_{d_0}},
$$
where
$$
\sup\limits_{\vec{v}_t}:=\sup\limits_{v^1_t}\sup\limits_{v^2_t}\cdots\sup\limits_{v^{k_i^t}_t}\s\mbox{and}\s \sup\limits_{\vec{b}_t}:=\sup\limits_{b^1_t}\sup\limits_{b^2_t}\cdots\sup\limits_{b^{l_i^t}_t}
$$
and
$$
\sup\limits_{\vec{\iota}_t}:=\sup\limits_{\iota^1_t}\sup\limits_{\iota^2_t}\cdots\sup\limits_{\iota^{r_i^t}_t}\s\mbox{and}\s \sup\limits_{\vec{\varrho}_t}:=\sup\limits_{\varrho^1_t}\sup\limits_{\varrho^2_t}\cdots\sup\limits_{\varrho^{s_i^t}_t}.
$$
At last, set
$$
B^i_r:=\sup\limits_{|\mathbf{k}|+|\mathbf{l}|+|\mathbf{r}|+|\mathbf{s}|=r}B^i_{\mathbf{klrs}}.
$$
Adopting the above notation we give
\begin{eqnarray*}
&&\sum\limits_{|\mathbf{p}|=1}^{\infty}\sum\limits_{i=1}^{d_0}\sum\limits_{\mathbf{m}_i+\mathbf{n}_i=(\mathbf{p}-e_i)^+}\sum\limits_{\mathbf{k}_i+\mathbf{l}_i=\mathbf{m}_i}\sum\limits_{\mathbf{r}_i+\mathbf{s}_i=\mathbf{n}_i}\sum\limits_{\substack{\vec{w}},\vec{v},\vec{f},\vec{b}}\sum\limits_{\vec{\tau},\vec{\iota},\vec{\sigma},\vec{\varrho}}\sum\limits_{\vartheta,\varepsilon}\sum\limits_{\varsigma\in\mathbb{Z}}\\
&&\Big|\int_{\mathbb{R}}\int_{\mathbb{T}^d}\left\langle C^i_{\varsigma\overrightarrow{wvfb\tau\iota\sigma\varrho}\vartheta\varepsilon},\beta Q^{\mp}_{\varsigma}\Pi_{\pm}(D)P_j\phi\right\rangle_{\mathbb{C}^{d_0}}Q^+_{\vec{w}_1}P_{\vec{v}_1}\Pi_+(D)(\psi^1_1-\psi^2_1)\nonumber\\
&&\cdots Q^+_{\vec{w}_{d_0}}P_{\vec{v}_{d_0}}\Pi_+(D)(\psi^1_{d_0}-\psi^2_{d_0})Q^+_{\vec{\tau}_1}P_{\vec{\iota}_1}\Pi_+(D)\psi^2_1\cdots Q^+_{\vec{\tau}_{d_0}}P_{\vec{\iota}_{d_0}}\Pi_+(D)\psi^2_{d_0}\nonumber\\
&&\cdot Q^-_{\vec{\sigma}_1}P_{\vec{\varrho}_1}\Pi_-(D)\varphi^2_1\cdots Q^-_{\vec{\sigma}_{d_0}}P_{\vec{\varrho}_{d_0}}\Pi_-(D)\varphi^2_{d_0}Q^-_{\vec{f}_1}P_{\vec{b}_1}\Pi_-(D)(\varphi^1_1-\varphi^2_1)\nonumber\\
&&\cdots Q^-_{\vec{f}_{d_0}}P_{\vec{b}_{d_0}}\Pi_-(D)(\varphi^1_{d_0}-\varphi^2_{d_0})\cdot Q^-_\vartheta P_\varepsilon\Pi_-(D)(\varphi_i^1-\varphi_i^2)\Big|\\
&\lesssim&\Big|\Big|\tilde{P}_{j+1}\phi\Big|\Big|_{S_{j+1}^{\mp}}\sum\limits_{|\mathbf{p}|=1}^{\infty}\sum\limits_{i=1}^{d_0}\sum\limits_{\mathbf{m}_i+\mathbf{n}_i=(\mathbf{p}-e_i)^+}2^{2j-3dj/2+d|(\mathbf{p}-e_i)^+|/2+j(2\max\{\lfloor d/2-d/b\rfloor,\lfloor d/2-d/a\rfloor\}+2)}\\
&&B^i_{|(\mathbf{p}-e_i)^+|}C^{|(\mathbf{p}-e_i)^+|/4-3d/2+d|(\mathbf{p}-e_i)^+|/2}\sum\limits_{\mathbf{k}_i+\mathbf{l}_i=\mathbf{m}_i}\sum\limits_{\mathbf{r}_i+\mathbf{s}_i=\mathbf{n}_i}\sum\limits_{\varsigma\in\mathbb{Z}}\sum\limits_{\vartheta}C^{-2^{-|\vartheta|-1}/(1+|(\mathbf{p}-e_i)^+|)}\\
&&\sum\limits_{\vec{w}:=(\vec{w}_1,\vec{w}_2,\cdots,\vec{w}_{d_0})}C^{-\sum\limits_{t=1}^{d_0}\sum\limits_{r_t=1}^{k_i^t}2^{-|w_t^{r_t}|-1}/(1+|(\mathbf{p}-e_i)^+|)}\sum\limits_{\vec{f}:=(\vec{f}_1,\vec{f}_2,\cdots,\vec{f}_{d_0})}C^{-\sum\limits_{t=1}^{d_0}\sum\limits_{c_t=1}^{l_i^t}2^{-|f_t^{c_t}|-1}/(1+|(\mathbf{p}-e_i)^+|)}\\
&&\sum\limits_{\vec{\tau}:=(\vec{\tau}_1,\vec{\tau}_2,\cdots,\vec{\tau}_{d_0})}C^{-\sum\limits_{t=1}^{d_0}\sum\limits_{s_t=1}^{r_i^t}2^{-|\tau_t^{s_t}|-1}/(1+|(\mathbf{p}-e_i)^+|)}\sum\limits_{\vec{\sigma}:=(\vec{\sigma}_1,\vec{\sigma}_2,\cdots,\vec{\sigma}_{d_0})}C^{-\sum\limits_{t=1}^{d_0}\sum\limits_{a_t=1}^{s_i^t}2^{-|\sigma_t^{a_t}|-1}/(1+|(\mathbf{p}-e_i)^+|)}\\
&&\sum\limits_{\vec{v}:=(\vec{v}_1,\vec{v}_2,\cdots,\vec{v}_{d_0})}\prod\limits_{t=1}^{d_0}\prod\limits_{r_t=1}^{k_i^t}\left\{2^{v_t^{r_t}d/2}\Big|\Big|P_{v_t^{r_t}}\Pi_+(D)(\psi^1_t-\psi_t^2)\Big|\Big|_{S^+_{v_t^{r_t}+1}}\right\}\\
&&\sum\limits_{\vec{b}:=(\vec{b}_1,\vec{b}_2,\cdots,\vec{b}_{d_0})}\prod\limits_{t=1}^{d_0}\prod\limits_{c_t=1}^{l_i^t}\left\{2^{b_t^{c_t}d/2}\Big|\Big|P_{b_t^{c_t}}\Pi_-(D)(\varphi^1_t-\varphi_t^2)\Big|\Big|_{S^-_{b_t^{c_t}+1}}\right\}\\
&&\sum\limits_{\vec{\iota}:=(\vec{\iota}_1,\vec{\iota}_2,\cdots,\vec{\iota}_{d_0})}\prod\limits_{t=1}^{d_0}\prod\limits_{s_t=1}^{r_i^t}\left\{2^{\iota_t^{s_t}d/2}\Big|\Big|P_{\iota_t^{s_t}}\Pi_+(D)\psi_t^2\Big|\Big|_{S^+_{\iota_t^{s_t}+1}}\right\}\\
&&\sum\limits_{\vec{\varrho}:=(\vec{\varrho}_1,\vec{\varrho}_2,\cdots,\vec{\varrho}_{d_0})}\prod\limits_{t=1}^{d_0}\prod\limits_{a_t=1}^{s_i^t}\left\{2^{\varrho_t^{a_t}d/2}\Big|\Big|P_{\varrho_t^{a_t}}\Pi_-(D)\psi_t^2\Big|\Big|_{S^-_{\varrho_t^{a_t}+1}}\right\}\\
&&\sum\limits_{\varepsilon}\left\{2^{\varepsilon d/2}\Big|\Big|P_{\varepsilon}\Pi_-(D)(\varphi_i^1-\varphi_i^2)\Big|\Big|_{S^+_{\varepsilon+1}}\right\}.
\end{eqnarray*}
Because
$$
\max\{||(\varsigma,\vec{w},\vec{f},\vec{\tau},\vec{\sigma},\zeta)||,||(\varsigma,\vec{w},\vec{f},\vec{\tau},\vec{\sigma},\vartheta)||,||(j,\vec{v},\vec{b},\vec{\iota},\vec{\varrho},\chi)||,||(j,\vec{v},\vec{b},\vec{\iota},\vec{\varrho},\varepsilon)||\}\leq 1,
$$
we have
\begin{eqnarray*}
&&\sum\limits_{|\mathbf{p}|=1}^{\infty}\sum\limits_{i=1}^{d_0}\sum\limits_{\mathbf{m}_i+\mathbf{n}_i=(\mathbf{p}-e_i)^+}\sum\limits_{\mathbf{k}_i+\mathbf{l}_i=\mathbf{m}_i}\sum\limits_{\mathbf{r}_i+\mathbf{s}_i=\mathbf{n}_i}\sum\limits_{\substack{\vec{w}},\vec{v},\vec{f},\vec{b}}\sum\limits_{\vec{\tau},\vec{\iota},\vec{\sigma},\vec{\varrho}}\sum\limits_{\vartheta,\varepsilon}\sum\limits_{\varsigma\in\mathbb{Z}}\\
&&\Big|\int_{\mathbb{R}}\int_{\mathbb{T}^d}\left\langle C^i_{\varsigma\overrightarrow{wvfb\tau\iota\sigma\varrho}\vartheta\varepsilon},\beta Q^{\mp}_{\varsigma}\Pi_{\pm}(D)P_j\phi\right\rangle_{\mathbb{C}^{d_0}}Q^+_{\vec{w}_1}P_{\vec{v}_1}\Pi_+(D)(\psi^1_1-\psi^2_1)\nonumber\\
&&\cdots Q^+_{\vec{w}_{d_0}}P_{\vec{v}_{d_0}}\Pi_+(D)(\psi^1_{d_0}-\psi^2_{d_0})Q^+_{\vec{\tau}_1}P_{\vec{\iota}_1}\Pi_+(D)\psi^2_1\cdots Q^+_{\vec{\tau}_{d_0}}P_{\vec{\iota}_{d_0}}\Pi_+(D)\psi^2_{d_0}\nonumber\\
&&\cdot Q^-_{\vec{\sigma}_1}P_{\vec{\varrho}_1}\Pi_-(D)\varphi^2_1\cdots Q^-_{\vec{\sigma}_{d_0}}P_{\vec{\varrho}_{d_0}}\Pi_-(D)\varphi^2_{d_0}Q^-_{\vec{f}_1}P_{\vec{b}_1}\Pi_-(D)(\varphi^1_1-\varphi^2_1)\nonumber\\
&&\cdots Q^-_{\vec{f}_{d_0}}P_{\vec{b}_{d_0}}\Pi_-(D)(\varphi^1_{d_0}-\varphi^2_{d_0})\cdot Q^-_\vartheta P_\varepsilon\Pi_-(D)(\varphi_i^1-\varphi_i^2)\Big|\\
&\lesssim&\Big|\Big|\tilde{P}_{j+1}\phi\Big|\Big|_{S_{j+1}^{\mp}}\sum\limits_{|\mathbf{p}|=1}^{\infty}\sum\limits_{i=1}^{d_0}\sum\limits_{\mathbf{m}_i+\mathbf{n}_i=(\mathbf{p}-e_i)^+}3^{1+|(\mathbf{p}-e_i)^+|}\\
&&2^{2j-3dj/2+d|(\mathbf{p}-e_i)^+|/2+j(2\max\{\lfloor d/2-d/b\rfloor,\lfloor d/2-d/a\rfloor\}+2)}\\
&&B^i_{|(\mathbf{p}-e_i)^+|}C^{|(\mathbf{p}-e_i)^+|/4-3d/2+d|(\mathbf{p}-e_i)^+|/2}\sum\limits_{\mathbf{k}_i+\mathbf{l}_i=\mathbf{m}_i}\sum\limits_{\mathbf{r}_i+\mathbf{s}_i=\mathbf{n}_i}\sum\limits_{\varsigma\in\mathbb{Z}}C^{-2^{-|\varsigma|}}\\
&&\sum\limits_{\vec{v}:=(\vec{v}_1,\vec{v}_2,\cdots,\vec{v}_{d_0})}\prod\limits_{t=1}^{d_0}\prod\limits_{r_t=1}^{k_i^t}\left\{2^{v_t^{r_t}d/2}\Big|\Big|P_{v_t^{r_t}}\Pi_+(D)(\psi^1_t-\psi_t^2)\Big|\Big|_{S^+_{v_t^{r_t}+1}}\right\}\\
&&\sum\limits_{\vec{b}:=(\vec{b}_1,\vec{b}_2,\cdots,\vec{b}_{d_0})}\prod\limits_{t=1}^{d_0}\prod\limits_{c_t=1}^{l_i^t}\left\{2^{b_t^{c_t}d/2}\Big|\Big|P_{b_t^{c_t}}\Pi_-(D)(\varphi^1_t-\varphi_t^2)\Big|\Big|_{S^-_{b_t^{c_t}+1}}\right\}\\
&&\sum\limits_{\vec{\iota}:=(\vec{\iota}_1,\vec{\iota}_2,\cdots,\vec{\iota}_{d_0})}\prod\limits_{t=1}^{d_0}\prod\limits_{s_t=1}^{r_i^t}\left\{2^{\iota_t^{s_t}d/2}\Big|\Big|P_{\iota_t^{s_t}}\Pi_+(D)\psi_t^2\Big|\Big|_{S^+_{\iota_t^{s_t}+1}}\right\}\\
&&\sum\limits_{\vec{\varrho}:=(\vec{\varrho}_1,\vec{\varrho}_2,\cdots,\vec{\varrho}_{d_0})}\prod\limits_{t=1}^{d_0}\prod\limits_{a_t=1}^{s_i^t}\left\{2^{\varrho_t^{a_t}d/2}\Big|\Big|P_{\varrho_t^{a_t}}\Pi_-(D)\psi_t^2\Big|\Big|_{S^-_{\varrho_t^{a_t}+1}}\right\}\\
&&\sum\limits_{\varepsilon}\left\{2^{\varepsilon d/2}\Big|\Big|P_{\varepsilon}\Pi_-(D)(\varphi_i^1-\varphi_i^2)\Big|\Big|_{S^+_{\varepsilon+1}}\right\}.
\end{eqnarray*}
Stepping modeled on (\ref{17}) leads to
$$
\sum\limits_{\varsigma\in\mathbb{Z}}C^{-2^{-|\varsigma|}}<\infty
$$
which implies
\begin{eqnarray*}
&&\sum\limits_{|\mathbf{p}|=1}^{\infty}\sum\limits_{i=1}^{d_0}\sum\limits_{\mathbf{m}_i+\mathbf{n}_i=(\mathbf{p}-e_i)^+}\sum\limits_{\mathbf{k}_i+\mathbf{l}_i=\mathbf{m}_i}\sum\limits_{\mathbf{r}_i+\mathbf{s}_i=\mathbf{n}_i}\sum\limits_{\substack{\vec{w}},\vec{v},\vec{f},\vec{b}}\sum\limits_{\vec{\tau},\vec{\iota},\vec{\sigma},\vec{\varrho}}\sum\limits_{\vartheta,\varepsilon}\sum\limits_{\varsigma\in\mathbb{Z}}\\
&&\Big|\int_{\mathbb{R}}\int_{\mathbb{T}^d}\left\langle C^i_{\varsigma\overrightarrow{wvfb\tau\iota\sigma\varrho}\vartheta\varepsilon},\beta Q^{\mp}_{\varsigma}\Pi_{\pm}(D)P_j\phi\right\rangle_{\mathbb{C}^{d_0}}Q^+_{\vec{w}_1}P_{\vec{v}_1}\Pi_+(D)(\psi^1_1-\psi^2_1)\nonumber\\
&&\cdots Q^+_{\vec{w}_{d_0}}P_{\vec{v}_{d_0}}\Pi_+(D)(\psi^1_{d_0}-\psi^2_{d_0})Q^+_{\vec{\tau}_1}P_{\vec{\iota}_1}\Pi_+(D)\psi^2_1\cdots Q^+_{\vec{\tau}_{d_0}}P_{\vec{\iota}_{d_0}}\Pi_+(D)\psi^2_{d_0}\nonumber\\
&&\cdot Q^-_{\vec{\sigma}_1}P_{\vec{\varrho}_1}\Pi_-(D)\varphi^2_1\cdots Q^-_{\vec{\sigma}_{d_0}}P_{\vec{\varrho}_{d_0}}\Pi_-(D)\varphi^2_{d_0}Q^-_{\vec{f}_1}P_{\vec{b}_1}\Pi_-(D)(\varphi^1_1-\varphi^2_1)\nonumber\\
&&\cdots Q^-_{\vec{f}_{d_0}}P_{\vec{b}_{d_0}}\Pi_-(D)(\varphi^1_{d_0}-\varphi^2_{d_0})\cdot Q^-_\vartheta P_\varepsilon\Pi_-(D)(\varphi_i^1-\varphi_i^2)\Big|\\
&\lesssim&\Big|\Big|\tilde{P}_{j+1}\phi\Big|\Big|_{S_{j+1}^{\mp}}\sum\limits_{|\mathbf{p}|=1}^{\infty}\sum\limits_{i=1}^{d_0}\sum\limits_{\mathbf{m}_i+\mathbf{n}_i=(\mathbf{p}-e_i)^+}3^{1+|(\mathbf{p}-e_i)^+|}B^i_{|(\mathbf{p}-e_i)^+|}C^{|(\mathbf{p}-e_i)^+|/4-3d/2+d|(\mathbf{p}-e_i)^+|/2}\\
&&2^{2j-3dj/2+d|(\mathbf{p}-e_i)^+|/2+j(2\max\{\lfloor d/2-d/b\rfloor,\lfloor d/2-d/a\rfloor\}+2)}\\
&&\sum\limits_{\mathbf{k}_i+\mathbf{l}_i=\mathbf{m}_i}\sum\limits_{\mathbf{r}_i+\mathbf{s}_i=\mathbf{n}_i}\Big|\Big|\Pi_+(D)(\psi^1-\psi^2)\Big|\Big|^{|\mathbf{k}_i|}_{S^{+,d/2}}\Big|\Big|\Pi_-(D)(\varphi^1-\varphi^2)\Big|\Big|^{|\mathbf{l}_i|+1}_{S^{-,d/2}}\\
&&\Big|\Big|\Pi_+(D)\psi^2\Big|\Big|^{|\mathbf{r}_i|}_{S^{+,d/2}}\Big|\Big|\Pi_-(D)\psi^2\Big|\Big|^{|\mathbf{s}_i|}_{S^{-,d/2}}\\
&\lesssim&\Big|\Big|\tilde{P}_{j+1}\phi\Big|\Big|_{S_{j+1}^{\mp}}\sum\limits_{|\mathbf{p}|=1}^{\infty}\sum\limits_{i=1}^{d_0}\sum\limits_{|\mathbf{k}_i|+|\mathbf{l}_i|+|\mathbf{r}_i|+|\mathbf{s}_i|=|(\mathbf{p}-e_i)^+|}3^{1+|(\mathbf{p}-e_i)^+|}B^i_{|(\mathbf{p}-e_i)^+|}\\
&&C^{|(\mathbf{p}-e_i)^+|/4-3d/2+d|(\mathbf{p}-e_i)^+|/2}2^{2j-3dj/2+d|(\mathbf{p}-e_i)^+|/2+j(2\max\{\lfloor d/2-d/b\rfloor,\lfloor d/2-d/a\rfloor\}+2)}\\
&&\Big|\Big|\Pi_+(D)(\psi^1-\psi^2)\Big|\Big|^{|\mathbf{k}_i|}_{S^{+,d/2}}\Big|\Big|\Pi_-(D)(\varphi^1-\varphi^2)\Big|\Big|^{|\mathbf{l}_i|}_{S^{-,d/2}}\Big|\Big|\Pi_+(D)\psi^2\Big|\Big|^{|\mathbf{r}_i|}_{S^{+,d/2}}\\
&&\Big|\Big|\Pi_-(D)\psi^2\Big|\Big|^{|\mathbf{s}_i|}_{S^{-,d/2}}\left\{||\psi^1-\psi^2||_{S^{+,d/2}}+||\varphi^1-\varphi^2||_{S^{-,d/2}}\right\},
\end{eqnarray*}
where we have used Theorem \ref{thm0}. Since the same method also works for (\ref{11}), we obtain
\begin{eqnarray*}
&&\Big|\int_{\mathbb{R}}\int_{\mathbb{T}^d}\langle F(\Pi_+(D)\psi^1+\Pi_-(D)\varphi^1)-F(\Pi_+(D)\psi^2+\Pi_-(D)\varphi^2),\beta\Pi_{\mp}(D)P_j\phi \rangle_{\mathbb{C}^{d_0}}\,dxdt\Big|\\
&\lesssim&2^{2j-3dj/2+j(2\max\{\lfloor d/2-d/b\rfloor,\lfloor d/2-d/a\rfloor\}+2)}\Big|\Big|\tilde{P}_{j+1}\phi\Big|\Big|_{S_{j+1}^{\mp}}\\
&&\sum\limits_{|\mathbf{p}|=1}^{\infty}\sum\limits_{i=1}^{d_0}\sum\limits_{|\mathbf{k}_i|+|\mathbf{l}_i|+|\mathbf{r}_i|+|\mathbf{s}_i|=|(\mathbf{p}-e_i)^+|}3^{1+|(\mathbf{p}-e_i)^+|}B^i_{|(\mathbf{p}-e_i)^+|}C^{|(\mathbf{p}-e_i)^+|/4-3d/2+d|(\mathbf{p}-e_i)^+|/2}\\
&&2^{d|(\mathbf{p}-e_i)^+|/2}\Big|\Big|\Pi_+(D)(\psi^1-\psi^2)\Big|\Big|^{|\mathbf{k}_i|}_{S^{+,d/2}}\Big|\Big|\Pi_-(D)(\varphi^1-\varphi^2)\Big|\Big|^{|\mathbf{l}_i|}_{S^{-,d/2}}\Big|\Big|\Pi_+(D)\psi^2\Big|\Big|^{|\mathbf{r}_i|}_{S^{+,d/2}}\\
&&\Big|\Big|\Pi_-(D)\psi^2\Big|\Big|^{|\mathbf{s}_i|}_{S^{-,d/2}}\left\{||\psi^1-\psi^2||_{S^{+,d/2}}+||\varphi^1-\varphi^2||_{S^{-,d/2}}\right\}\\
&\leq&2^{2j-3dj/2+j(2\max\{\lfloor d/2-d/b\rfloor,\lfloor d/2-d/a\rfloor\}+2)}\Big|\Big|\tilde{P}_{j+1}\phi\Big|\Big|_{S_{j+1}^{\mp}}\\
&&\sum\limits_{|\mathbf{p}|=1}^{\infty}\sum\limits_{|\mathbf{k}|+|\mathbf{l}|+|\mathbf{r}|+|\mathbf{s}|=|\mathbf{p}|}3^{1+|\mathbf{p}|}A_{|\mathbf{p}|}C^{|\mathbf{p}|/4-3d/2+d|\mathbf{p}|/2}2^{d|\mathbf{p}|/2}\\
&&\Big|\Big|\Pi_+(D)(\psi^1-\psi^2)\Big|\Big|^{|\mathbf{k}|}_{S^{+,d/2}}\Big|\Big|\Pi_-(D)(\varphi^1-\varphi^2)\Big|\Big|^{|\mathbf{l}|}_{S^{-,d/2}}\Big|\Big|\Pi_+(D)\psi^2\Big|\Big|^{|\mathbf{r}|}_{S^{+,d/2}}\\
&&\Big|\Big|\Pi_-(D)\psi^2\Big|\Big|^{|\mathbf{s}|}_{S^{-,d/2}}\left\{||\psi^1-\psi^2||_{S^{+,d/2}}+||\varphi^1-\varphi^2||_{S^{-,d/2}}\right\}
\end{eqnarray*}
where we have defined $A_r:=\max\{B^i_r|1\leq i\leq d_0\}$ and used the fact $|(\mathbf{p}-e_i)^+|\leq|\mathbf{p}|$. Recalling $d\geq9$ and $a=b=4$, it is not difficult to get
$$
2j-3dj/2+j(2\max\{\lfloor d/2-d/b\rfloor,\lfloor d/2-d/a\rfloor\}+2)<-dj/2.
$$
Define
$$
g(x,y,z,w):=\sum\limits_{p=1}^{\infty}3^{1+p}A_pC^{p/4-3d/2+dp/2}2^{dp/2}\sum\limits_{k+l+r+s=p}x^ky^lz^rw^s
$$
for $x,y,z,w\in[0,\epsilon]$. Obviously, the following inequality holds true
\begin{eqnarray*}
g(x,y,z,w)&\leq&\sum\limits_{p=1}^{\infty}3^{1+p}A_pC^{p/4-3d/2+dp/2}2^{dp/2}\left(\sum\limits_{k=0}^px^k\right)\left(\sum\limits_{l=0}^py^l\right)\left(\sum\limits_{r=0}^pz^r\right)\left(\sum\limits_{s=0}^pw^s\right)\\
&\leq&\epsilon^4\sum\limits_{p=1}^{\infty}3^{1+p}A_pC^{p/4-3d/2+dp/2}2^{dp/2}(p+1)^4.
\end{eqnarray*}
Thanks to
$$
\limsup\limits_{p\rightarrow\infty}A_p^{1/p}<C^{-1/4-d/2}2^{-d/2}/3
$$
we get
$$
\sum\limits_{p=1}^{\infty}3^{1+p}A_pC^{p/4-3d/2+dp/2}2^{dp/2}(p+1)^4<\infty.
$$
This completes the proof.
\endproof

\section{Appendix}
\begin{thm}\label{thm0}
The  following linear operator is bounded
$$
\Pi_{\pm}(D): S^{\pm,d/2}\longrightarrow S^{\pm,d/2}.
$$
\end{thm}
\textbf{Proof.} From the definition it follows that
\begin{eqnarray*}
&&||\Pi_{\pm}(D)P_jf||_{S^{\pm}_{j+1}}:=||\Pi_{\pm}(D)P_jf||_{L_t^{\infty}L_x^2}+||\Pi_{\pm}(D)P_jf||_{\dot{X}^{\pm,1/2,\infty}}\\
&&+\sup\limits_{\substack{0\leq k'\leq j+1\\(d+2)(j+1)/(2d-2)\leq l\leq j+1}}\Bigg\{2^{-\frac{k'+j+1}{a}}||\Pi_{\pm}(D)P_jf||_{L_t^aL_x^b[l,k']}\\
&&+2^{-\frac{k'+j+1}{b}}||\Pi_{\pm}(D)P_jf||_{L_t^bL_x^a[l,k']}\Bigg\}.
\end{eqnarray*}
It is easy to get
\begin{eqnarray*}
&&||\Pi_{\pm}(D)P_jf||_{L_t^{\infty}L_x^2}=\Bigg|\Bigg|\frac{1}{2}\left\{I_{d_0}\pm\left(\sum\limits_{\dot{k}=1}^d\alpha^{\dot{k}}\xi_{\dot{k}}+\beta\right)/\langle\xi\rangle\right\}\varphi(2^{-j}\cdot)\mathcal{F}_{x\mapsto\xi}f\Bigg|\Bigg|_{L_t^{\infty}l_{\xi}^2}\\
&\lesssim&||\varphi(2^{-j}\cdot)\mathcal{F}_{x\mapsto\xi}f||_{L_t^{\infty}l_{\xi}^2}=||P_jf||_{L_t^{\infty}L_x^2}.
\end{eqnarray*}
Similarly we have
$$
||\Pi_{\pm}(D)P_jf||_{\dot{X}^{\pm,1/2,\infty}}\lesssim||P_jf||_{\dot{X}^{\pm,1/2,\infty}}.
$$
Now we will control
\begin{eqnarray*}
||\Pi_{\pm}(D)P_jf||_{L_t^aL_x^b[l,k']}=\sum\limits_{\kappa\in\mathcal{K}_l}\sum\limits_{n\in\Xi_{k'}}||\Gamma_{k',n}P_{\kappa}\Pi_{\pm}(D)P_jf||_{L_t^aL_x^b}.
\end{eqnarray*}
By Gagliardo-Nirenberg-Moser estimates we can get
\begin{eqnarray*}
&&||\Gamma_{k',n}P_{\kappa}\Pi_{\pm}(D)P_jf||_{L_t^aL_x^b}\\
&\lesssim&\sum\limits_{|\alpha|=\lfloor d/2-d/b\rfloor+1}\Big|\Big|\Gamma_{k',n}P_{\kappa}\Pi_{\pm}(D)P_jf\Big|\Big|^{\theta}_{L_t^aL_x^2}\Big|\Big|D^{\alpha}\Gamma_{k',n}P_{\kappa}\Pi_{\pm}(D)P_jf\Big|\Big|^{1-\theta}_{L_t^aL_x^2}\\
&\lesssim&\sum\limits_{|\alpha|=\lfloor d/2-d/b\rfloor+1}\Big|\Big|\Gamma_{k',n}P_{\kappa}P_jf\Big|\Big|^{\theta}_{L_t^aL_x^2}\Big|\Big|D^{\alpha}\Gamma_{k',n}P_{\kappa}P_jf\Big|\Big|^{1-\theta}_{L_t^aL_x^2}
\end{eqnarray*}
with $\theta:=d/\{2(\lfloor d/2-d/b\rfloor+1)\}-d/\{b(\lfloor d/2-d/b\rfloor+1)\}$. Since of
$$
\mbox{supp}\{\mathcal{F}_{x\mapsto\xi}(\Gamma_{k',n}P_{\kappa}P_jf)\}\subset\left\{\xi\Big|n_{\dot{k}}-2^{k'}\leq\xi_{\dot{k}}\leq n_{\dot{k}}+2^{k'},1\leq\dot{k}\leq d, 2^j\leq|\xi|\leq 2^{j+2}\right\}\bigcap 2\dot{\Gamma}_{\kappa}
$$
with $n\in\Xi_{k'}:=2^{k'}\cdot\mathbb{Z}^d$, we have
\begin{eqnarray*}
&&\Big|\Big|D^{\alpha}\Gamma_{k',n}P_{\kappa}P_jf\Big|\Big|_{L_t^aL_x^2}=\Big|\Big|D^{\alpha}(\hat{\eta}_{\kappa}\ast\Gamma_{k',n}P_{\kappa}P_jf)\Big|\Big|_{L_t^aL_x^2}\\
&=&\Big|\Big|D^{\alpha}\hat{\eta}_{\kappa}\ast\Gamma_{k',n}P_{\kappa}P_jf\Big|\Big|_{L_t^aL_x^2}\leq\Big|\Big|\Gamma_{k',n}P_{\kappa}P_jf\Big|\Big|_{L_t^aL_x^2}||D^{\alpha}\hat{\eta}_{\kappa}||_{L^{\infty}_x}\\
&\lesssim&\Big|\Big|\Gamma_{k',n}P_{\kappa}P_jf\Big|\Big|_{L_t^aL_x^b}||D^{\alpha}\hat{\eta}_{\kappa}||_{L^{\infty}_x}
\end{eqnarray*}
with $\hat{\eta}_{\kappa}:=\mathcal{F}^{-1}_{\xi\mapsto x}\tilde{\eta}_{\kappa}$ and $\tilde{\eta}_{\kappa}$ being a function such that $0\leq\tilde{\eta}_{\kappa}\leq1$, $\mbox{supp}\{\tilde{\eta}_{\kappa}\}\subset 3\dot{\Gamma}_{\kappa}\cap\{\xi|2^j\leq|\xi|\leq 2^{j+2}\}$ and $\tilde{\eta}_{\kappa}(\xi)=1$ for $\xi\in 2\dot{\Gamma}_{\kappa}\cap\{\xi|2^j\leq|\xi|\leq 2^{j+2}\}$, where we have applied Young's inequality and the fact that flat torus has finite measure. Direct computation gives
\begin{eqnarray*}
||D^{\alpha}\hat{\eta}_{\kappa}||_{L^{\infty}_x}=||\xi^{\alpha}\tilde{\eta}_{\kappa}||_{l^1_{\xi}}\lesssim 2^{j|\alpha|}||\tilde{\eta}_{\kappa}||_{l^1_{\xi}}\leq 2^{j|\alpha|}|3\dot{\Gamma}_{\kappa}\cap\{\xi|2^j\leq|\xi|\leq 2^{j+2}\}|.
\end{eqnarray*}
Because the diameter of $\kappa$ is $2^{-l}$ and $l\geq j(d+2)/(2d-2)$,
$$
|3\dot{\Gamma}_{\kappa}\cap\{\xi|2^j\leq|\xi|\leq 2^{j+2}\}|\lesssim 2^{-(d-1)l+j(1+d/2)}\leq 1.
$$
So it is easy to get
$$
||D^{\alpha}\hat{\eta}_{\kappa}||_{L^{\infty}_x}\lesssim1
$$
which means
\begin{eqnarray*}
||\Pi_{\pm}(D)P_jf||_{L_t^aL_x^b[l,k']}\lesssim||P_jf||_{L_t^aL_x^b[l,k']}.
\end{eqnarray*}
Then the desired theorem holds true.
\endproof

\textbf{Acknowledgement}
The author is supported by Fundamental Research Funds for the Central Universities (No. 2021MS045).

\end{document}